\documentclass[12pt]{amsart}
\usepackage{amsmath,amsthm, amssymb, epsf, amscd}
\usepackage{times}

\numberwithin{equation}{section}

\newtheorem{Thm}{Theorem}[section]
\newtheorem{Prop}[Thm]{Proposition}
\newtheorem{Lem}[Thm]{Lemma}
\newtheorem{Cor}[Thm]{Corollary}
\newtheorem{Conj}[Thm]{Conjecture}

\theoremstyle{remark}
\newtheorem{Rem}[Thm]{Remark}
\newtheorem{Exa}[Thm]{Example}

\theoremstyle{definition}

\newtheorem{Prob}[Thm]{Problem}
\newtheorem{Ass}[Thm]{Assumption}

\newcommand{\fig}[1]
        {\raisebox{-0.5\height}%
                 {\epsfbox{#1}}%
        }

\newcommand{\cfig}[1]
        {\begin{center}\fig{#1}\end{center}%
        }
\def\Z{{\mathbb Z}}
\def\R{{\mathbb R}}
\def\Q{{\mathbb Q}}

\def\calA{\mathcal{A}}
\def\calB{\mathcal{B}}

\def\calD{\mathcal{D}}

\def\calG{\mathcal{G}}

\def\calI{\mathcal{I}}
\def\calJ{\mathcal{J}}

\def\calL{\mathcal{L}}

\def\calP{\mathcal{P}}

\def\calS{\mathcal{S}}

\def\sfA{\mathsf{A}}
\def\sfD{\mathsf{D}}
\def\sfL{\mathsf{L}}
\def\sfM{\mathsf{M}}
\def\sfS{\mathsf{S}}
\def\sfX{\mathsf{X}}
\def\rma{\mathrm{a}}
\def\rmb{\mathrm{b}}
\def\rmc{\mathrm{c}}
\def\rmd{\mathrm{d}}
\def\rme{\mathrm{e}}
\def\Im{\mathrm{Im}\,}
\def\deg{\mathrm{deg}\,}
\def\span{\mathrm{span}}

\def\eqdef{\stackrel{\rm def}{=}}
\def\Aut{\mathrm{Aut}\,}
\def\pAut{\mathrm{Aut}'\,}
\def\nequiv{\equiv\kern-.9em /\kern.4em}
\newcommand{\mapright}[1]{%
	\smash{\mathop{%
		\hbox to 1cm{\rightarrowfill}}\limits^{#1}}}

\newcommand{\tworows}[2]{\stackrel{#1}{\scriptscriptstyle #2}}

\begin{document}
\title[Configuration space integral for long $n$-knots]{Configuration space integral for long $n$-knots, the Alexander polynomial and knot space cohomology}

\author[T. Watanabe]{Tadayuki Watanabe}
\address{Research Institute for Mathematical Sciences, Kyoto University, Kyoto, Japan}
\email{tadayuki@kurims.kyoto-u.ac.jp}

\date{\today}

\begin{abstract}
There is a higher dimensional analogue of the perturbative Chern-Simons theory in the sense that a similar perturbative series as in 3-dimension, which is computed via configuration space integral, yields an invariant of higher dimensional knots (Bott-Cattaneo-Rossi invariant), which is constructed by Bott for degree 2 and by Cattaneo-Rossi for higher degrees. However, its feature is yet unknown. In this paper we restrict the study to long ribbon $n$-knots and characterize the Bott-Cattaneo-Rossi invariant as a finite type invariant of long ribbon $n$-knots introduced in \cite{HKS}. As a consequence, we obtain a non-trivial description of the Bott-Cattaneo-Rossi invariant in terms of the Alexander polynomial. 

The results for higher codimension knots are also given. In those cases similar differential forms to define Bott-Cattaneo-Rossi invariant yields infinitely many cohomology classes of $\mathrm{Emb}(\R^n, \R^m)$ if $m,n\geq 3$ odd and $m>n+2$. We observe that half of these classes are non-trivial, along a line similar to Cattaneo-CottaRamusino-Longoni \cite{CCL}.
\end{abstract}
\maketitle
\section{Introduction}

Witten gave in \cite{Wi} a path integral formulation of the Chern-Simons field theory which gives a framework for understanding many invariants of knots in a 3-manifold. But path integral is not yet well defined mathematically. Though, the perturbative expansion method for the path integral gives a mathematical definition of knot invariants. It is defined by configuration space integral with certain uni-trivalent graphs called Jacobi diagrams. This construction of invariants has been studied by Axelrod-Singer, Bar-Natan, Guadagnini-Martellini-Mintchev, Kohno, Kontsevich, Bott-Taubes, Dylan Thurston and others \cite{AS, BN, GMM, Koh, Kon, BT, AF, T}. It is known that the configuration space integral invariant of knots in $S^3$ is a universal finite type (or Vassiliev) invariant \cite{AF,T}. 

There is a higher dimensional analogue of these theory with a certain variant of Jacobi diagrams. The suitable diagrams for higher dimensions are graphs with two kinds of edges corresponding to the ``angular forms" in $\R^{n+2}$ and in $\R^n$ respectively. The diagrams have degrees given by half the number of vertices. In his seminal paper, Bott constructed an invariant for higher dimensional knots \cite{Bot} associated to degree 2 diagrams. After that Cattaneo and Rossi gave a path integral formula for invariants of higher dimensional embedded manifolds into a manifold and as a perturbative expansion of it, they obtained invariants for higher dimensional knots, which will be denoted by ${z}_k$ in this paper, associated to degree $k$ diagrams \cite{CR, R}. It may be a generalization of the 3-dimensional perturbative Chern-Simons theory and Bott's invariant is the degree 2 part of Cattaneo-Rossi's invariant. So we call $z_k$ the Bott-Cattaneo-Rossi invariant, or the BCR invariant for short. More precisely, Cattaneo-Rossi proved that $z_k$ (or its framing correction) is an isotopy invariant if $n$ is odd (or $(n,k)=(2,3)$). The BCR invariant currently seems to be a new invariant of long $n$-knots. However, it has not been known whether $z_k$ is non-trivial or not. So we try to understand its features, in particular, whether there is some connection with known invariants. 

In this paper we study the invariant ${z}_k$ restricting to a certain class of long $n$-knots called long ribbon $n$-knots. The class of long ribbon $n$-knots is a rather familiar class known to have similarities to the classical knot theory in 3-dimension. In particular, in \cite{HKS}, Habiro, Kanenobu and Shima introduced the notion of finite type invariant of ribbon 2-knots, which is straightforwardly generalizable to long ribbon $n$-knots, and showed that the coefficients of the Alexander polynomial expanded as a power series in $t-1$, are finite type invariants. Finite type invariant of type $k$ is defined by the condition that a given invariant vanishes at any `$(k+1)$-th order differential' of a ribbon 2-knot. Moreover, Habiro and Shima obtained in \cite{HS} a remarkable result saying that the set of all finite type invariants is isomorphic as a graded algebra to the polynomial algebra in the coefficients of the Alexander polynomial. 

Restricting to long ribbon $n$-knots, we characterize the BCR invariant $z_k$ as a finite type invariant (Theorem~\ref{thm:Zuniversal}). In particular, we determine its `highest order term' explicitly. As a consequence, we obtain a description in terms of the coefficients of the logarithm of the Alexander polynomial at $t=e^h$ with an explicit highest order term. In particular, it shows that the BCR invariant is non-trivial. The proof is done by choosing an embedding for a long ribbon $n$-knot in some extreme situation and then computing $z_k$ explicitly at the $k$-th order differentials of long ribbon $n$-knots. In the limit, the computation of the integral is highly reduced and we can compute explicitly. The presence of two kinds of edges in the diagrams appearing in the definition of the BCR invariant makes the computation slightly complicated. But the spirit in the computation is similar to that of \cite{AF,T}.

We consider also the cases of long $n$-knots of codimension $>2$. Recently, Budney proved in \cite{Bud} that the space $\mathrm{Emb}(\R^n,\R^m)$ of long embeddings is $(2m-3n-4)$-connected if $2m-3n-3\geq 0$, namely, all homotopy groups in dimensions $\leq 2m-3n-4$ vanish. In this paper we study the cohomology of $\mathrm{Emb}(\R^n,\R^m)$ in dimensions higher than Budney's bound $2m-3n-4$. A similar proof to Theorem~\ref{thm:Zuniversal} works to prove the non-triviality of some cohomology classes of $\mathrm{Emb}(\R^n,\R^m)$ in this range. This study is inspired mainly by the result of Cattaneo, Cotta-Ramusino and Longoni, which concerns the case $n=1$, proving that the Chern-Simons perturbative series classes are non-trivial in cohomology of $\mathrm{Emb}(\R^1,\R^m)$ when $m>3$ \cite{CCL}. In the case of higher dimensions $n>1$, $m-n>2$, the analogous construction to the BCR invariant also yields cohomology classes of $\mathrm{Emb}(\R^n,\R^m)$ if $m,n$ are odd. We prove that half of these classes are non-trivial and thus give an estimate below for $H^*(\mathrm{Emb}(\R^n,\R^m);\R)$. The strategy of the proof is similar to \cite{CCL}, that is to construct a certain cycle in $\mathrm{Emb}(\R^n,\R^m)$ and to evaluate the cohomology class on the cycle.

We also study the BCR invariant for some classes, which we call long handle knots, other than ribbon. We observe that for long handle knots, the BCR invariant is expressed non-trivially in terms of the coefficients of Levine's generalizations of the Alexander polynomial in \cite{Lev}.

This paper is organized as follows. In \S2, the definition of the BCR invariant $z_k$ is given and the Cattaneo-Rossi's result concerning the invariance of their invariant is stated. The BCR invariant is defined as a linear combination of the integrals over the configuration spaces associated with certain graphs.  In \S3, we will explain the generality on the space of graphs from which each coefficient in $z_k$ at a configuration space integral of a graph are determined. This section is somewhat digressive but needed in the invariance proof of $z_k$ in the appendix. In \S4, we restrict the study to long ribbon $n$-knots. The notion of Habiro-Kanenobu-Shima's finite type invariant is recalled and the main theorem (Theorem~\ref{thm:Zuniversal}) about the characterization of $z_k$ as a finite type invariant is stated and proved. The relation between $z_k$ and the Alexander polynomial is explained at the end of \S4. In \S5, higher codimension results are given. We prove that higher BCR invariant classes are non-trivial. The proof is somewhat parallel to that of the Theorem~\ref{thm:Zuniversal} and sometimes we refer some words from the former. So if the reader has read the proof of Theorem~\ref{thm:Zuniversal}, then the proof of non-triviality may be more understandable. In \S\ref{s:remark}, we remark some result about the value of ${z}_k$ for some long $n$-knots other than ribbon and state a few problems. Appendix~A is devoted to the self-contained proof of invariance of the BCR invariant ${z}_k$, which is originally described in \cite{R}, filling the details not explicitly described in \cite{R}. 
\par\vspace{2mm}

\noindent{\bf Acknowledgments.}\ 
The author would like to thank his adviser Professor T.~Ohtsuki for helpful comments, careful reading of a manuscript and encouragements. The author would also like to thank Professors/Doctors A.~Cattaneo, K.~Habiro, T.~Kohno, C.~Lescop, R.~Budney, C.~Rossi and referee for reading an earlier version of this paper, pointing out some errors, and for helpful suggestions. Especially, C.~Rossi helped the author to understand their proof of invariance by correcting the author's misunderstandings and by providing his thesis.

\section{Invariants of long $n$-knots}

We shall review the definition of the BCR invariant. Roughly in this section, we review the definitions of:
\begin{itemize}
\item Jacobi diagrams and weights $w_k$ of them,
\item configuration space integral associated to a Jacobi diagram,
\end{itemize}
and define the BCR invariant as a linear combination of configuration space integrals for diagrams whose coefficients are the weights of them.

Let $n\geq 2$. A {\it long $n$-knot} is the image of a smooth long embedding $\psi$ of $\R^n$ into $\R^{n+2}$ that is standard near $\infty$, i.e., we assume that there exists an $(n+2)$-ball $D\subset\R^{n+2}$ such that $\Im \psi\cap(\R^{n+2}\setminus D)=(\R^n\times\{0\}\times\{0\})\cap(\R^{n+2}\setminus D)$.

\subsection{Jacobi diagrams}

A {\it Jacobi diagram} is an oriented graph $\Gamma$ with valence at most 3 and a choice of a vertex orientation where $\Gamma$ has two kinds of edges, {\it $\theta$-edges} (depicted by directed dashed lines) and {\it $\eta$-edges} (depicted by directed solid lines), such that :
\begin{itemize}
\item The admissible combinations of incident edges to a vertex are:
\cfig{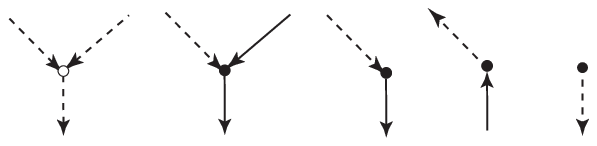}
We will call a trivalent vertex with three incident $\theta$-edges an {\it internal vertex} and in pictures we depict such a vertex by a white dot. We will call a non-internal vertex an {\it external vertex} and in pictures we denote it by a black dot. 
\item The vertex orientation of $\Gamma$ is a choice of ordering of two ingoing $\theta$-edges incident to each internal vertex modulo even number of swappings. In order to represent vertex orientations, we assume that diagrams may be depicted so that the three incident $\theta$-edges $(\mbox{the outgoing edge of $v$}, \mbox{1st incident $\theta$-edge}, \mbox{2nd incident $\theta$-edge})$ are arranged in the anti-clockwise order in a plane diagram.
\end{itemize}

\begin{figure}
\cfig{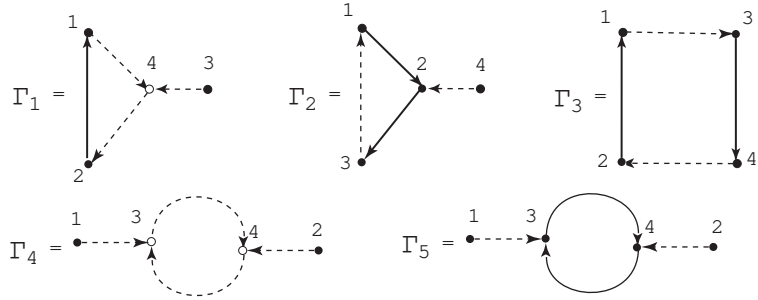}
\caption{All degree 2 connected Jacobi diagrams}\label{fig:deg2d}
\end{figure}

The {\it degree} of a Jacobi diagram is defined to be half the number of vertices. The complete list of connected Jacobi diagrams of degree 2 up to vertex orientations is shown in Figure~\ref{fig:deg2d}.

An automorphism of a Jacobi diagram is an automorphism $\varphi$ on the set of vertices sending external vertices to external vertices (and thus internal vertices to internal vertices) and inducing a bijection on the set of edges which sends each oriented edge $e=(i,j)$ to an oriented edge $(\varphi(i),\varphi(j))$. We call a $\theta$-edge connecting two external vertices a {\it chord}. We denote the group of automorphisms of $\Gamma$ by $\Aut\Gamma$.

\subsubsection{The weight function $w_k$ on Jacobi diagrams}
Let $\calG_k^0$ denote the the set all degree $k$ connected Jacobi diagrams. We will define the invariants of higher dimensional knots in \S\ref{ss:def-integral} as $\R$-linear combinations of the weight function $w_k$ on Jacobi diagrams, namely, a certain map $w_k:\calG_k^0\to\R$, so that their values are in $\R$. $w_k$ is defined as follows.

Let $\Gamma\in\calG_k^0$ be any Jacobi diagram. If $\Gamma$ has a subgraph of one of the following form:
\begin{equation}\label{eq:y-l}
 \fig{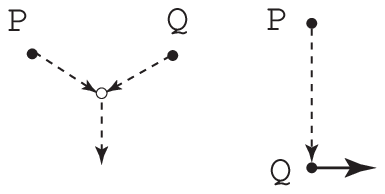} 
\end{equation}
where we assume that  no other edges are ingoing to both $P$ and $Q$, we set $w_k(\Gamma)=0$.

If $\Gamma$ does not have such subgraphs, then the form of $\Gamma$ is rather restricted: it must be a cyclic alternating sequence of the following two paths with $\theta$-edges stuck into:
\begin{equation}\label{eq:hairy-strut}
 \fig{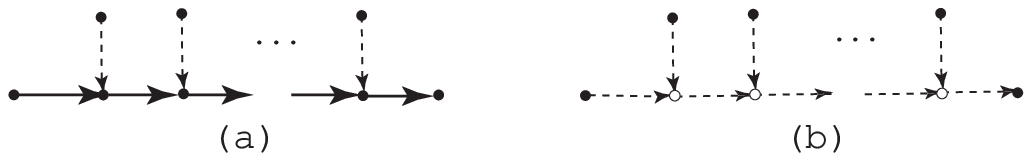} 
 \end{equation}
In particular, any such diagram includes just one cycle. The following picture is a typical example of this observation:
\[ \fig{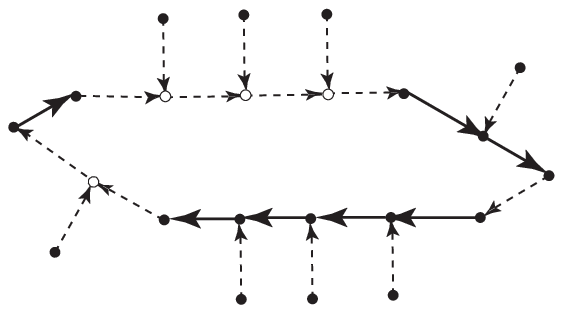} \]

Now we define the weight function $w_k$ for such diagrams. Then we associate a sign with each trivalent vertex (not only internal vertices) by the following rule:
\[ \fig{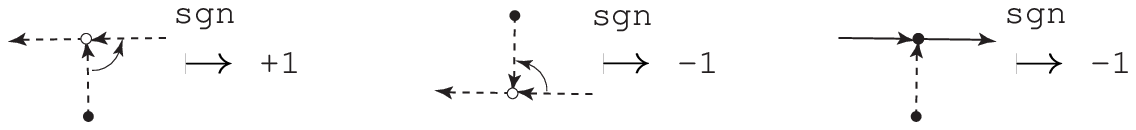} \]
and set
\[ w_k(\Gamma)=\prod_{v:\mbox{\tiny trivalent vertex}}\mathrm{sgn}(v) \]
for any $\Gamma$ without subgraphs as in (\ref{eq:y-l}). Now the map $w_k:\calG_k^0\to\R$ is defined.

\subsection{Configuration space}\label{ss:configuration-space}
Let $\psi:\R^n\to\R^{n+2}$ be a smooth embedding for a long $n$-knot which is standard near $\infty$. We consider that $\R^{n+2}$ is the complement of a fixed point $\infty$ in $S^{n+2}$. 

Let $\Gamma$ be a Jacobi diagram with $q$ external and $s$ internal vertices. The {\it configuration space associated with $(\Gamma, \psi)$} denoted by $C_{\Gamma}^0(\psi)$ is defined by
\[ \begin{split}
	C_{\Gamma}^0(\psi)\eqdef \{&(x_1,\ldots,x_q;x_{q+1},\ldots x_{q+s}),\\
		&x_1,\ldots,x_q\in \R^n, x_{q+1},\ldots,x_{q+s}\in\R^{n+2}\times S^{n-1}\times S^{n-1}\,|\,\\
		&p_1(x_i)\neq p_1(x_j)\quad \mbox{if $i,j\geq q+1$, $i\neq j$, $(i,j)$: edge of $\Gamma$},\\
		&\psi(x_i)\neq\psi(x_j)\quad \mbox{if $i,j\leq q$, $i\neq j$, $(i,j)$: edge of $\Gamma$},\\
		&\psi(x_i)\neq p_1(x_j)\quad \mbox{if $i\leq q$, $j\geq q+1$, $(i,j)$: edge of $\Gamma$}\},
	\end{split} \]
where $p_1:\R^{n+2}\times S^{n-1}\times S^{n-1}\to\R^{n+2}$ is the projection onto the first factor. Note that $C^0_{\Gamma}=\overline{C}^0_{\Gamma}\times (S^{n-1})^{\times a}$ for some non negative integer $a$ where
\[ \begin{split}
	\overline{C}_{\Gamma}^0(\psi)\eqdef \{&(x_1,\ldots,x_q;x_{q+1},\ldots x_{q+s}),\\
		&x_1,\ldots,x_q\in \R^n, x_{q+1},\ldots,x_{q+s}\in\R^{n+2}\,|\,\\
		&x_i\neq x_j\quad \mbox{if $i,j\geq q+1$, $i\neq j$, $(i,j)$: edge of $\Gamma$},\\
		&\psi(x_i)\neq\psi(x_j)\quad \mbox{if $i,j\leq q$, $i\neq j$, $(i,j)$: edge of $\Gamma$},\\
		&\psi(x_i)\neq x_j\quad \mbox{if $i\leq q$, $j\geq q+1$, $(i,j)$: edge of $\Gamma$}\}.
	\end{split} \]
We include the factor $S^{n-1}\times S^{n-1}$ in $C_{\Gamma}^0$ because we want to make the target of the Gauss map for a $\theta$-edge, which will be defined in the next subsection, uniformly $S^{n+1}\times S^{n-1}$. By this trick, treatments of signs and degrees of forms on $C_{\Gamma}^0$ may become easier since the volume form on $S^{n+1}\times S^{n-1}$ has even degree $2n$ and we will later assign each $\theta$-edge a $2n$ degree form. We will call this additional factor a dummy factor. There is a bijective correspondence with the set of $\theta$-edges directed to internal vertices with the spheres in the dummy factor. (See Remark~\ref{rem:dummy} below).

To see the convergence of the integrals over $C_{\Gamma}^0$, we use the compactification $C_{\Gamma}(\psi)$ of $C_{\Gamma}^0(\psi)$ used in \cite{AS, BT}, in analogy of Fulton-MacPherson \cite{FM}, which is a smooth manifold with corners. It is obtained by a sequence of blow-ups along most of the diagonals (see \S\ref{ss:face}) and the infinity. Detailed descriptions of the compactification used in this paper is found in \cite{R}. We use a slightly modified version of the compactification as the one used in \cite{T} so that we do not blow-up along the diagonal 
\[ \{(x_1,\ldots,x_{k})\,|\,x_i=x_j, \mbox{$i$ and $j$ are not connected by an edge}\}\]
while we do blow-up along its lower dimensional sub-diagonals.

Note that the addition of the dummy factor does not affect the blow-ups since the addition changes normal bundles over each diagonal just by direct product of spheres (see \S\ref{ss:face}), i.e., taking the direct product with the dummy factor commutes with the blow-ups.

\subsection{Integral over configuration space}\label{ss:def-integral}
\begin{figure}
\cfig{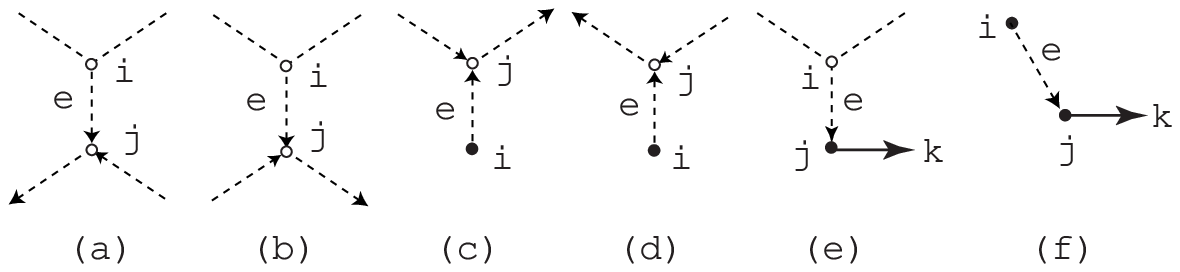}
\caption{}\label{fig:typeedges}
\end{figure}

Let $E_{\mathrm{\theta}}(\Gamma)$ be the set of $\theta$-edges in $\Gamma$. We define a form on $C_{\Gamma}(\psi)$ of degree $2n|E_{\mathrm{\theta}}(\Gamma)|$ by
\[ \omega(\Gamma)\eqdef \bigwedge_{e\in E_{\mathrm{\theta}}(\Gamma)}\phi_e^*(\omega_{n+1}\wedge\omega_{n-1}), \]
where $\omega_{p}$ is the $SO(p+1)$-invariant unit volume form on $S^p$. The Gauss map $\phi_e:C_{\Gamma}(\psi)\to S^{n+1}\times S^{n-1}$ is defined by
\[ \left\{
	\begin{array}{ll}
		(u(p_1(x_j)-p_1(x_i)), p_3(x_j))&\mbox{if $e$ is as in Figure~\ref{fig:typeedges}(a)},\\
		(u(p_1(x_j)-p_1(x_i)), p_2(x_j))&\mbox{if $e$ is as in Figure~\ref{fig:typeedges}(b)},\\
		(u(p_1(x_j)-\psi(x_i)),p_3(x_j))&\mbox{if $e$ is as in Figure~\ref{fig:typeedges}(c)},\\
		(u(p_1(x_j)-\psi(x_i)),p_2(x_j))&\mbox{if $e$ is as in Figure~\ref{fig:typeedges}(d)},\\
		(u(\psi(x_j)-p_1(x_i)), u'(x_k-x_j))&\mbox{if $e$ is as in Figure~\ref{fig:typeedges}(e)},\\
		(u(\psi(x_j)-\psi(x_i)), u'(x_k-x_j))&\mbox{if $e$ is as in Figure~\ref{fig:typeedges}(f)},\\
	\end{array}\right. \]
where $p_2, p_3:\R^{n+2}\times S^{n-1}\times S^{n-1}\to S^{n-1}$ is the projections onto the second and the third factor respectively, $u:\R^{n+2}\setminus\{0\}\to S^{n+1}, u':\R^n\setminus\{0\}\to S^{n-1}$ are defined by
\begin{equation*}
 u(x)=\frac{x}{\|x\|},\quad u'(x)=\frac{x}{\|x\|}. 
\end{equation*}

Let $\Omega$ be an orientation on $C_{\Gamma}(\psi)$. We define
\[ I(\Gamma,\Omega)\eqdef \int_{C_{\Gamma}(\psi)}\omega(\Gamma), \]
in which the orientation on $C_{\Gamma}$ is given by $\Omega$. Let
\[ I(\Gamma)\eqdef I(\Gamma,\Omega(\Gamma)), \]
where $\Omega(\Gamma)$ is defined in Lemma~\ref{lem:ori} below.
\begin{Lem}\label{lem:ori}
Let $\Gamma$ be a Jacobi diagram and let $\overline\Gamma$ be $\Gamma$ with its vertex orientation reversed. Then there exists a choice of orientation $\Omega=\Omega(\Gamma)$ on $C_{\Gamma}(\psi)$ such that
\begin{equation}\label{eq:omega}
\Omega(\overline\Gamma)=(-1)^n\Omega(\Gamma).
\end{equation}
\end{Lem}
\begin{proof}
Let $(X^1,\ldots, X^{n+2}, Y^1, Y^2)\in\R^{n+2}\times S^{n-1}\times S^{n-1}$ be the coordinate of an internal vertex of $\Gamma$ and $(X^1,\ldots, X^n)\in \R^n$ be the coordinate of an external vertex of $\Gamma$. Decompose each $\theta$-edge into two half edges. To each half edge $\bar{e}$, we associate an $n$-form $\Omega_{\bar{e}}$ which are combined into an orientation on $C_{\Gamma}$ as follows.

If $\bar{e}$ is one of three half edges meeting at an internal vertex $v$ with the vertex orientation given by a bijection $o_v:\{\mbox{the two ingoing half edges}\}\to\{1,2\}$, set
\[ \Omega_{\bar{e}}\eqdef \left\{
	\begin{array}{ll}
		dX_v^{o_v}\wedge\omega_{n-1}(v^{(o_v)}) & \mbox{if $\bar{e}$ is ingoing and $o_v(\bar{e})=1, 2$,}\\
		dX_v^3\wedge \cdots\wedge dX_v^{n+2} & \mbox{if $\bar{e}$ is outgoing,}\\
	\end{array}\right.
\]
where $(v^{(1)},v^{(2)})$ denote the coordinate of $v$ on $S^{n-1}\times S^{n-1}$. If $\bar{e}$ is a half edge meeting an external vertex, then set
\[ \Omega_{\bar{e}}\eqdef dX_v^1\wedge\cdots\wedge dX_v^n. \]

In terms of this, we define $\Omega_e\eqdef \Omega_{\bar{e}_+}\wedge\Omega_{\bar{e}_-}$ for each $\theta$-edge $e$ with $e=(\bar{e}_-,\bar{e}_+)$ and define
\[	\Omega(\Gamma)\eqdef \bigwedge_{e\in E_{\mathrm{\theta}}(\Gamma)}\Omega_e.\]
Then the property (\ref{eq:omega}) is because
\[\begin{split}
	 &(dX_v^1\wedge\omega_{n-1}(v^{(1)})\wedge\Omega_{\bar{e}})\wedge (dX_v^2\wedge\omega_{n-1}(v^{(2)})\wedge\Omega_{\bar{e}'})\\
	 =&(-1)^n(dX_v^1\wedge\omega_{n-1}(v^{(1)})\wedge\Omega_{\bar{e}'})\wedge (dX_v^2\wedge\omega_{n-1}(v^{(2)})\wedge\Omega_{\bar{e}})
	 \end{split} \]
for some half-edges $\bar{e}$ and $\bar{e}'$.
\end{proof}

Further, $\omega(\overline{\Gamma})=(-1)^{n+1}\omega(\Gamma)$ because
\[ \theta_{kj}\omega_{n-1}^{(1)}\theta_{ij}\omega_{n-1}^{(2)}=(-1)^{n+1}\theta_{ij}\omega_{n-1}^{(1)}\theta_{kj}\omega_{n-1}^{(2)} \]
where $\theta_{ij}\eqdef u(x_j-x_i)^*\omega_{n+1}$ and $\omega_{n-1}^{(i)}$ is a copy of $\omega_{n-1}$. This together with Lemma~\ref{lem:ori} implies
\[ I(\overline\Gamma,\Omega(\overline\Gamma))=-I(\Gamma,\Omega(\Gamma)). \]
It follows from the definition of $w_k$ that $w_k(\overline{\Gamma})=-w_k(\Gamma)$. Therefore the product $I(\Gamma)w_k(\Gamma)$ does not depend on the vertex orientation of $\Gamma$. But it is still sensitive under a reversion of the orientation of an edge. This is the reason why we fix the edge-orientation.

\begin{Rem}\label{rem:dummy}
The integral $I(\Gamma)$ and the integral for $\Gamma$ defined in \cite{CR} differ by the integration along the dummy factor with a suitable orientation. The results of the two integrals coincide. We will check this in Example~\ref{exa:bott-inv}.

The reason to introduce the dummy factor is as follows. If we use $\overline{C}_{\Gamma}^0$, the suitable form in the integrand for a $\theta$-edge $(i,j)$ with internal target may be $\theta_{ij}\eqdef u(\psi(x_j)-\psi(x_i))^*\omega_{n+1}$. So if there are two $\theta$-edges with a common target, the total form on $\overline{C}_{\Gamma}^0$ may include a part like $\theta_{ij}\wedge \theta_{kl}$. Since $\theta_{ij}\wedge \theta_{kl}=(-1)^{n+1}\theta_{kl}\wedge \theta_{ij}$, the sign of the integrand form depends on the choice of the order on such $\theta$-edges if $n$ is even. The addition of the dummy factor allows us to avoid this problem and we can give a formula for the configuration space integral for general diagrams without choosing particular orders on $\theta$-edges, although in this paper even dimensional case is considered only for $n=2$ and low degree graphs. However, we consider it may be worth for future research to give such a formulation.
\end{Rem}

\subsection{Invariants of long $n$-knots}\label{ss:invariant}
For $k\geq 2$, let 
\begin{equation}\label{eq:Z}
z_k(\psi)\eqdef \frac{1}{2}\sum_{\Gamma\in\calG_k^0}\frac{I(\Gamma)(\psi)w_k(\Gamma)}{|\Aut\Gamma|}\in\R
\end{equation}
where $|\Aut\Gamma|$ denotes the order of the group $\Aut\Gamma$.

In \cite{CR, R}, Cattaneo and Rossi remark that they proved the following theorem.
\begin{Thm}\label{thm:invariant}
If $n$ is an odd integer $\geq 3$ and $k\geq 2$, then $z_k$ is an isotopy invariant of long $n$-knots.

If $n=2$, then there exists a 2-form $\rho$ on $C_1\eqdef Bl(S^2,\{\infty\})$ such that $\hat{z}_3\eqdef z_3+\int_{C_1}\rho$ is an isotopy invariant of long $2$-knots, where $Bl(S^2,\{\infty\})$ denotes the blow-up of $S^2$ along $\infty\in S^2$.
\end{Thm}
We give a self-contained proof of Theorem~\ref{thm:invariant} in the appendix.

By Theorem~\ref{thm:invariant}, we may write $z_k(K)$ (or $\hat{z}_3(K)$) instead of $z_k(\psi)$ (resp. $\hat{z}_3(\psi)$) if $K=\Im\psi$.
\begin{Exa}\label{exa:bott-inv}
From Figure~\ref{fig:deg2d}, the degree 2 term $z_2$ may be given as follows:
\[ z_2=\frac{1}{2}I(\Gamma_1)w_2(\Gamma_1)+\frac{1}{2}I(\Gamma_2)w_2(\Gamma_2)+\frac{1}{4}I(\Gamma_3)w_2(\Gamma_3)
	+\frac{1}{4}I(\Gamma_4)w_2(\Gamma_4)+\frac{1}{4}I(\Gamma_5)w_2(\Gamma_5). \]
We have
\[ w_2(\Gamma_1)=-w_2(\Gamma_2)=w_2(\Gamma_3)=w_2(\Gamma_4)=w_2(\Gamma_5)=1. \]
Moreover, $I(\Gamma_4)=I(\Gamma_5)=0$ because $I(\Gamma_4)$ is the integral of $8n$-form over $8n-(n-1)$ dimensional submanifold of $(S^{n+1}\times S^{n-1})^{4}$ (since the targets of two edges by the map $u$ coincide) and similarly for $I(\Gamma_5)$. Hence
\[ z_2=\frac{1}{2}I(\Gamma_1)-\frac{1}{2}I(\Gamma_2)+\frac{1}{4}I(\Gamma_3). \]

Now we shall rewrite this formula in terms of the integrals over $C_{q,s}$, the compactification of $\overline{C}_{\Gamma}$ defined as $C_{\Gamma}$.

$I(\Gamma_1)$ can be reduced to the integral over $C_{3,1}$ as follows. Assume that the orientation on $C_{3,1}$ is the one naturally induced from $\R^n\times \R^n\times \R^n\times \R^{n+2}$. Let $(x_1,x_2,x_3,x_4)$ denote the coordinate on $C_{3,1}$ and $x_4=(y_1,\cdots, y_{n+2})\in\R^{n+2}$ be the coordinate for $x_4$ component.

By Lemma~\ref{lem:ori}, 
\[ \begin{split}
	\Omega(\Gamma_1)=&(d^nx_1\wedge dy_3\wedge\cdots\wedge dy_{n+2})
		\wedge (dy_1\wedge\omega_{n-1}(v^{(1)})\wedge d^nx_2)\\
		&\wedge (dy_2\wedge\omega_{n-1}(v^{(2)})\wedge d^nx_3)\\
	=-&d^nx_1\wedge d^nx_2\wedge d^nx_3\wedge d^{n+2}x_4
		\wedge\omega_{n-1}(v^{(1)})\wedge\omega_{n-1}(v^{(2)})
	\end{split}\]
where $d^nx$ denotes a volume element for the $x$ component. Let $\theta_{ij}\eqdef u(\psi(x_j)-\psi(x_i))^*\omega_{n+1}$ and $\eta_{ij}\eqdef u'(x_j-x_i)^*\omega_{n-1}$. Then by definition,
\[ 	\begin{split}
	\omega(\Gamma_1)&=\theta_{41}\eta_{12}\theta_{24}\omega_{n-1}^{(1)}\theta_{34}\omega_{n-1}^{(2)}
		= (-1)^{(n-1)(n+1)}\theta_{41}\eta_{12}\theta_{24}\theta_{34}\omega_{n-1}^{(1)}\omega_{n-1}^{(2)}\\
	&=(-1)^{n-1}\theta_{41}\theta_{24}\theta_{34}\eta_{12}\omega_{n-1}^{(1)}\omega_{n-1}^{(2)}
		=(-1)^{n-1}(-1)^n\theta_{14}\theta_{24}\theta_{34}\eta_{12}\omega_{n-1}^{(1)}\omega_{n-1}^{(2)}\\
	&=-\theta_{14}\theta_{24}\theta_{34}\eta_{12}\omega_{n-1}^{(1)}\omega_{n-1}^{(2)}.
	\end{split}\]
Thus 
\[ I(\Gamma_1)=\int_{C_{3,1}\times S^{n-1}\times S^{n-1}}\theta_{14}\theta_{24}\theta_{34}\eta_{12}\omega_{n-1}^{(1)}\omega_{n-1}^{(2)}
	=\int_{C_{3,1}}\theta_{14}\theta_{24}\theta_{34}\eta_{12}. \]

A similar computation yields
\[ 	I(\Gamma_2)=\int_{C_{4,0}}\theta_{13}\theta_{24}\eta_{12}\eta_{34}
	\quad\mbox{and}\quad
	I(\Gamma_3)=\int_{C_{4,0}}\theta_{13}\theta_{24}\eta_{12}\eta_{23}
\]
where the orientation on $C_{4,0}$ is assumed induced from the one on $\R^n\times \R^n\times \R^n\times \R^n$. Therefore we have
\[ z_2=\frac{1}{2}\int_{C_{3,1}}\theta_{14}\theta_{24}\theta_{34}\eta_{12}
	-\frac{1}{2}\int_{C_{4,0}}\theta_{13}\theta_{24}\eta_{12}\eta_{34}
	+\frac{1}{4}\int_{C_{4,0}}\theta_{13}\theta_{24}\eta_{12}\eta_{23}. \]	
This is precisely $\frac{1}{2}$ times the invariant in \cite{Bot} written in the notations of \cite{CR}.
\qed
\end{Exa}

\subsection{Vanishing of some terms of $z_k$}

By the following proposition, most terms of $z_k$ in lower odd degrees vanish in the case $n$ is odd.
\begin{Prop}\label{prop:vanish-axial-symm}
Suppose $n$ is odd and $\Gamma$ is a Jacobi diagram of odd degree with an axial symmetry, not having a subgraph as in (\ref{eq:y-l}). Then $I(\Gamma)=0$.
\end{Prop}
\begin{proof}
For each odd degree Jacobi diagram $\Gamma$, we consider the Jacobi diagram ${\Gamma}^*$ obtained by reversing the orientations of all the edges involved in the cycle in $\Gamma$. 

We consider two integrals $I({\Gamma}^*)$ and $I(\Gamma)$. By the assumption, ${\Gamma}^*$ and $\Gamma$ are equal in $\calG_k$ up to vertex orientation and ${\Gamma}^*$ may be obtained from $\Gamma$ by the axial symmetry inducing an automorphism $S:C_{\Gamma}\to C_{\Gamma}$ on the associated configuration space. One of the following axial symmetries may occur:
\[ \fig{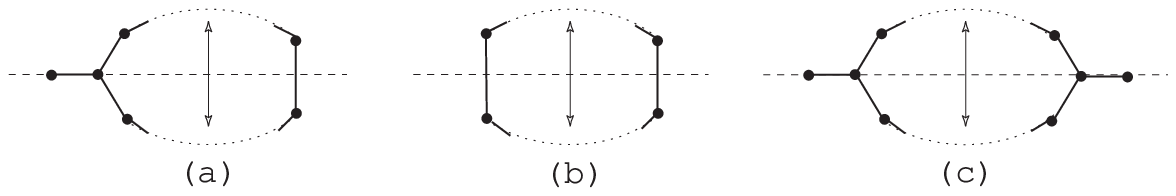} \]

In the case (a), $S$ transposes even number of pairs. So it does not change the orientation of $C_{\Gamma}$. Since the length of the cycle in (a) is odd, $S$ reverses the sign of the integrand form. One has
\[ I(\Gamma)=\int_{C_{\Gamma}}\omega(\Gamma)=+\int_{C_{\Gamma}}S^*\omega(\Gamma)=-\int_{C_{\Gamma}}\omega(\Gamma) \]
and hence $I(\Gamma)=0$. 

In the case (b) and (c), $S$ reverses the orientation of $C_{\Gamma}$ since $S$ transposes odd number of pairs. Moreover, $S$ preserves the sign of the integrand form since the length of the cycle is even. By the same reason as in (a), one concludes $I(\Gamma)=0$.
\end{proof}

If $\Gamma$ without internal vertices does not have a subgraph as in (\ref{eq:y-l}), then a chord of $\Gamma$ is one of the following form:
\begin{equation}\label{eq:type-chord}
\fig{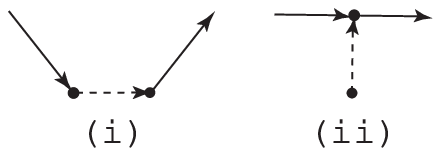}
\end{equation}
Let $k_1$ and $k_2\ (k=k_1+k_2)$ denote the numbers of chords of type (i) and (ii) respectively.
\begin{Prop}\label{prop:odd-chord-diag}
Suppose $n$ is odd and $\Gamma$ is a Jacobi diagram without internal vertices, not having a subgraph as in (\ref{eq:y-l}). Then
\begin{equation}\label{eq:Iw=Iw}
 I(\Gamma^*)w_k(\Gamma^*)=(-1)^kI(\Gamma)w_k(\Gamma) 
 \end{equation}
where $\Gamma^*$ denotes $\Gamma$ with the orientations of all edges involved in its cycle reversed.
\end{Prop}
\begin{proof}
One can show that
\begin{enumerate}
\item $w_k(\Gamma^*)=w_k(\Gamma)=(-1)^{k_2}$.
\item $\Omega(\Gamma^*)=(-1)^{k_1}\Omega(\Gamma)$.
\item $\omega(\Gamma^*)=(-1)^{k_2}\omega(\Gamma)$.
\end{enumerate}
(1) is by definition. (2) is because the number of chords involved in the cycle of $\Gamma$ is $k_1$. (3) is because the number of edges involved in the cycle of $\Gamma$ is $2k_1+k_2$. Hence (\ref{eq:Iw=Iw}) holds.
\end{proof}

\section{The space of Jacobi diagrams}

We shall see that the weight function $w_k:\calG_k^0\to\R$ satisfies some axioms, which are needed in the proof of Theorem~\ref{thm:invariant} and see that $w_k$ may arise naturally from the generality on some space of Jacobi diagrams defined by the axioms. The reader who does not need to read the proof of Theorem~\ref{thm:invariant} may skip this section.

Let $\R\calG_k^0$ denote the vector space spanned by elements of $\calG_k^0$. We define {\it ST, SU, STU, C relations} on $\R\calG_k^0$ as follows\footnote{If the external chains are replaced by 1-handles and the $\theta$-part is replaced by ``oriented tree claspers" \cite{W}, then these relations are essentially those appeared in \cite{HS}.}:
\begin{equation}\label{eq:STU}
\fig{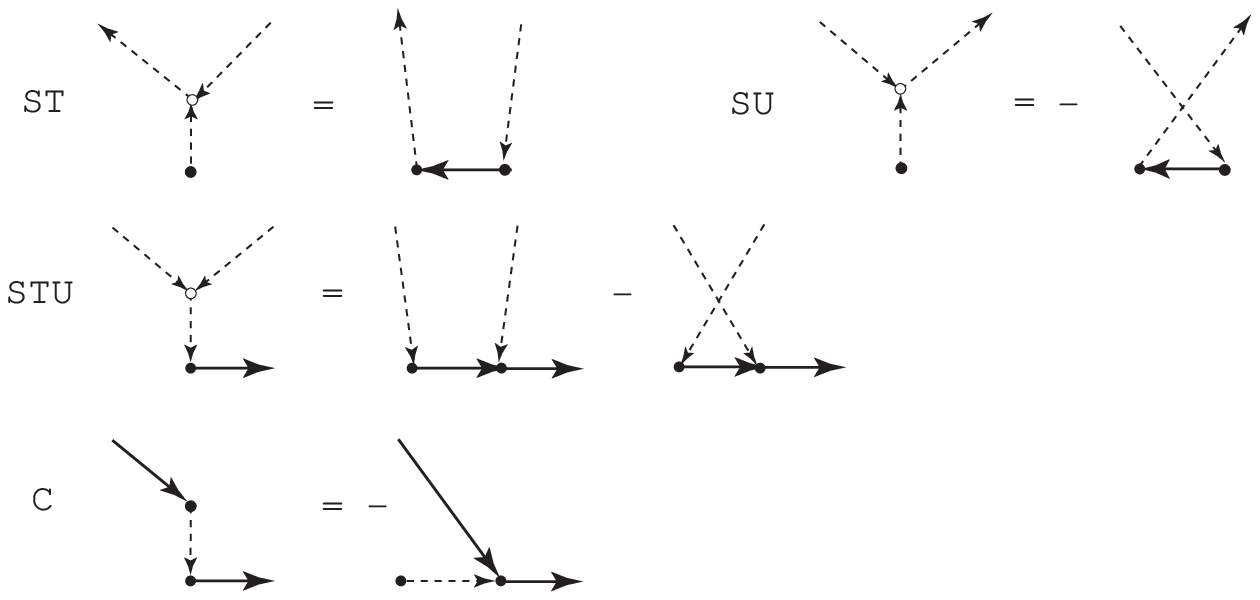}
\end{equation}
Here the graphs in (\ref{eq:STU}) are subgraphs of Jacobi diagrams. Other edges may incident to vertices in the subgraphs. Let 
\begin{equation*}
	\calA_k\eqdef \R\calG_k^0/\mathrm{relations\ in\ (\ref{eq:STU})}.
\end{equation*}
Then as in the 3-dimensional perturbative Chern-Simons theory, one can define a graph valued series:
\[ z_k^*(\psi)\eqdef\frac{1}{2}\sum_{\Gamma\in\calG_k^0}\frac{I(\Gamma)(\psi)[\Gamma]}{|\Aut\Gamma|}\in\calA_k. \]
In 3-dimension, the graph valued perturbative invariant is the strongest one in the sense that any scalar valued perturbative invariant can be obtained from the graph valued one via some weight function on graphs. But in higher dimensions, the graph valued one is not so strong than a scalar valued one, as observed in the rest of this section.
\begin{Prop}\label{prop:IHX}
The following relations are satisfied in $\calA_k$:
\[\fig{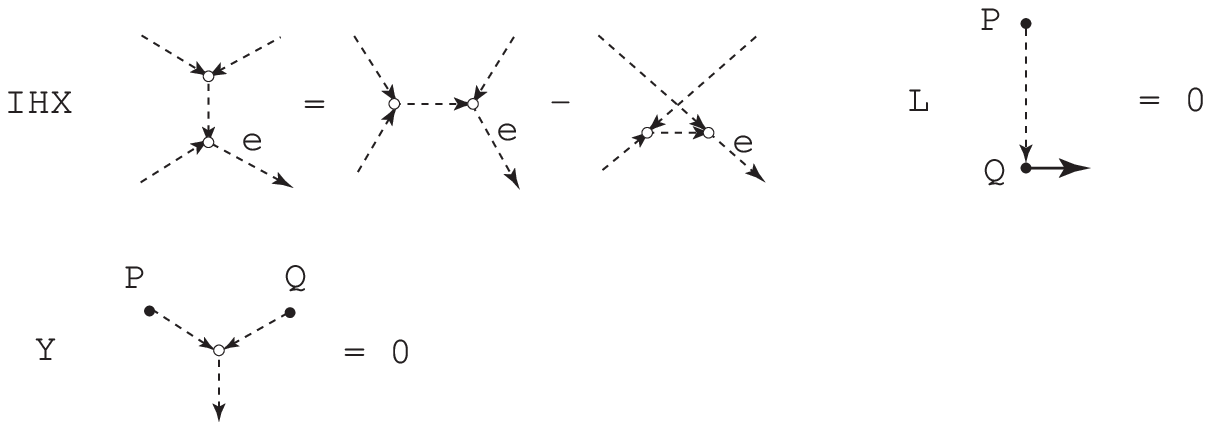}\]
Here we assume that no other edges are ingoing to both $P$ and $Q$.
\end{Prop}
\begin{proof}
For IHX relation, we apply some STU relations to the component connected to the distinguished edge $e$ in the relation just in the same way for the three terms until $e$ touches an external vertex. Then IHX relation follows from the following expansion by STU relations:
\[ \fig{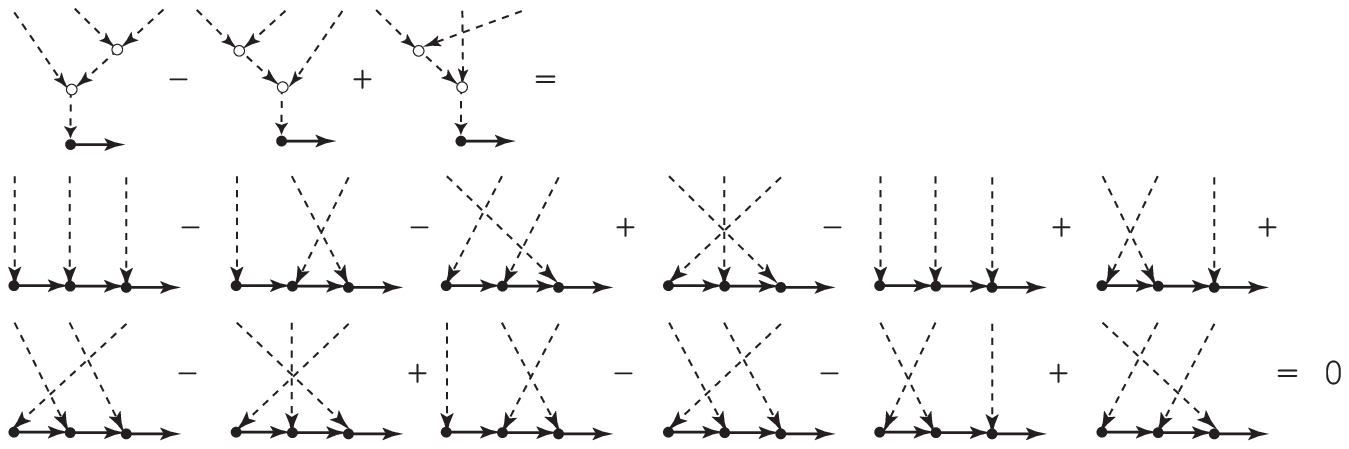} \]
Similarly, Y relation follows from the following expansion by an STU relation:
\[ \fig{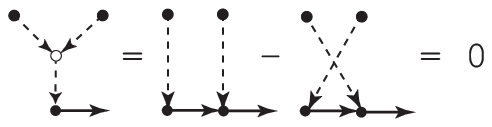} \]
A sequence of applications of ST relations to the LHS of L relation yields the LHS of Y relation as follows:
\cfig{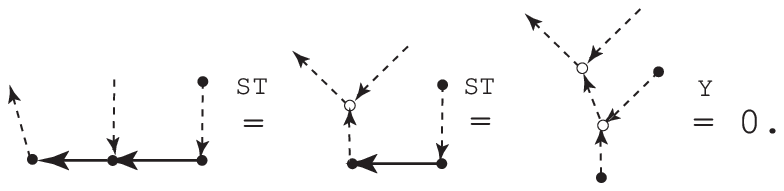}
\end{proof}
The vector space structure on $\calA_k$ is determined as in the following propositions.
\begin{Prop}\label{prop:A_n}
\begin{enumerate}
\item $\theta$ part of a Jacobi diagram forms a disjoint union of chords, trees and wheels. Here wheels are Jacobi diagrams of the following forms:
\cfig{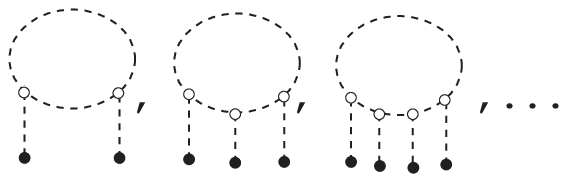}
\item Extend $w_k$ linearly to a map $w_k:\R\calG_k^0\to \R$. Then the map $w_k:\R\calG_k^0\to\R$ descends to a well defined $\R$-linear map $w_k:\calA_k\to\R$ which spans $\calA_k^*$.
\item The dimension of $\calA_k$ is 1 for each $k\geq 2$.
\end{enumerate}
\end{Prop}
\begin{proof}
(1) immediately follows from the edge-orientation condition of internal vertices. 

(2) follows from a direct check of the fact that $w_k$ satisfies the relations in (\ref{eq:STU}).

For (3), observe first that any Jacobi diagram can be transformed into a Jacobi diagram without $\eta$-edges by ST, SU, C and Y relation. The resulting graph must be a wheel-like graph plus terms annihilated by Y relation. So the dimension of $\calA_k$ is at most one. By (2), $w_k$ descends to a non-trivial linear map $w_k:\calA_k\to\R$, so the dimension of $\calA_k$ is at least one. Therefore (3) follows.
\end{proof}

Since the dimension of $\calA_k$ is one, any linear map $\calA_k\to\R$ is a scalar multiple of $w_k$. Further, we obtain the following identity:
\[ z_k^*(\psi)=z_k(\psi)[\Gamma],\ z_k(\psi)=w_k\circ z_k^*(\psi) \]
for some $[\Gamma]\in\calA_k$. Therefore, the BCR invariant is equivalent to the graph valued one.


\section{BCR invariant for long ribbon $n$-knots}\label{s:proof}


We restrict our study to a certain class of long $n$-knots called long ribbon $n$-knots and characterize the BCR invariant as a finite type invariant in the sense of \cite{HKS}. By this characterization, we obtain a description of the BCR invariant in terms of the Alexander polynomial.

\subsection{Long ribbon $n$-knots and ribbon presentations}

\subsubsection{Long ribbon $n$-knots}

A {\it long ribbon $(n+1)$-disk} is the image of an immersion $f$ of a lower half $(n+1)$-plane $D_-^{n+1}=\{(x_1,\cdots,x_{n+1})\in\R^{n+1}\,|\,x_{n+1}\leq 0 \}$ into $\R^{n+2}$ such that the singularity of $f$ consists of finitely many ribbon singularities and such that $f$ is standard outside a sufficiently large $(n+1)$-ball in $\R^{n+1}$ whose center is the origin. Here a {\it ribbon singularity} is an $n$-disk consisting of transverse double points and whose pre-image consists of a proper $n$-disk in $D_-^{n+1}$ and an $n$-disk in the interior of $D_-^{n+1}$. A {\it long ribbon $n$-knot} is a long $n$-knot bounding a long ribbon $(n+1)$-disk. 

\subsubsection{Ribbon presentations}

We use ribbon presentations to present long ribbon $n$-knots. 
A {\it ribbon presentation} $P=\calD\cup \calB$ in $\R^3$ is an immersed oriented 2-disk into $\R^3$ with a base point on its boundary having a decomposition consists of $p+1$ disjoint embedded 2-disks $\calD=D_0\sqcup D_1\sqcup D_2\sqcup\ldots\sqcup D_{p}$ and $p$ disjoint embedded bands $\calB=B_1\sqcup B_2\sqcup\ldots\sqcup B_p$ satisfying the following conditions.

\begin{enumerate}
\item The base point is on the boundary of $D_0$.
\item Each band $B_i$ transversely intersects the interiors of the disks in $\calD$ (as in Figure~\ref{fig:crossing}(i)).
\item Each end of a band is attached to the boundary of some disk in $\calD$.
\end{enumerate}
An example of a ribbon presentation is depicted in Figure~\ref{fig:ex-starlike}.

\subsubsection{Associating a long ribbon $n$-knot to a ribbon presentation}
We can construct a ribbon $(n+1)$-disk $V_P$ associated to each ribbon presentation $P$ by 
\[ V_P=N(\calD\times [-2,2]^{n-1}\cup \calB\times [-1,1]^{n-1})\subset\R^3\times\R^{n-1},\]
where $N(\,)$ denotes a smoothing of corners. Then we obtain a long ribbon $n$-knot associated to $P$ by taking the boundary of $V_P$ followed by connect summing the standardly embedded plane $\R^n$ at the point associated with the base point (the base point may be thickened to form $D^{n-1}$. Then define the base point of the long $n$-knot to be its center). Note that the connect summing of $\R^n$ is not unique depending on the position of the base point. We denote by $K_P$ the long ribbon $n$-knot associated to a (based) ribbon presentation $P$ and by $[K_P]$ its isotopy class. It is known that any (long) ribbon $n$-knot is isotopic to the one associated to some ribbon presentation.


\subsubsection{Crossings}

\begin{figure}
\fig{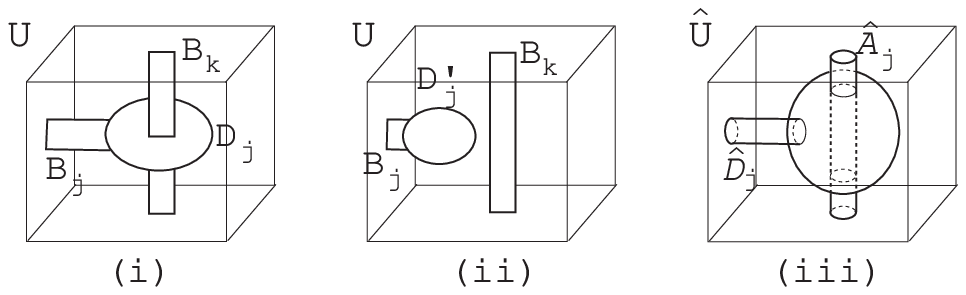}
\caption{}\label{fig:crossing}
\end{figure}

Now we give a definition of crossings of both ribbon presentations and long ribbon $n$-knots. For crossings of ribbon presentations, consider a regular neighborhood $U$ of a disk $D_j\subset \calD$ (of $P=\calD\cup\calB$) including a ribbon singularity inside so that $U\cap P$ consists of 
\begin{itemize}
\item the disk $D_j$,
\item a part of the band $B_j\subset\calB$ incident to $D_j$,
\item a part of some band $B_k$ intersecting $D_j$.
\end{itemize}
We call the triple $[U,U\cap (D_j\cup B_j),U\cap B_k]$ a {\it crossing of $P$} (see Figure~\ref{fig:crossing}(i)).

In the associating process of $P\leadsto K_P$ above, $U$ may also be thickened by taking direct product with $[-3,3]^{n-1}$ to include the associated pieces thickened from $U\cap(D_j\cup B_j)$ and $U\cap B_k$. We call the triple $(\widehat{U}, \widehat{D}_j,\widehat{A}_j)$ a crossing of $K_P$ associated to $[U,U\cap (D_j\cup B_j),U\cap B_k]$ where
\begin{itemize}
\item $\widehat{U}=U\times [-3,3]^{n-1}\subset \R^{n+1}$,
\item $\widehat{D}_j$ is one of the components of $\widehat{U}\cap K_P$ associated to $U\cap(D_j\cup B_j)$,
\item $\widehat{A}_j$ is one of the components of $\widehat{U}\cap K_P$ associated to $U\cap B_k$.
\end{itemize}
Observe that $\widehat{D}_j\cong D^n$ and $\widehat{A}_j\cong S^{n-1}\times I$. $(\widehat{U},\widehat{D}_j,\widehat{A}_j)$ looks like Figure~\ref{fig:crossing}(iii).

A crossing of a long ribbon $n$-knot isotopic to $K_P$ is a triple $(\widetilde{U}, \widetilde{D}_j,\widetilde{A}_j), (\widetilde{D}_j\cup\widetilde{A}_j\subset\widetilde{U})$ which is obtained from an associated crossing $(\widehat{U}, \widehat{D}_j,\widehat{A}_j)$ on $K_P$ as above by an isotopy deforming $K_P\cup\widehat{U}$ such that 
\begin{itemize}
\item it induces an isotopy deforming the trivial link $K_P\cap \partial\widehat{U}$ in $\partial\widehat{U}$,
\item the number of the components in $K_P\cap\partial\widehat{U}$ is preserved during the isotopy,
\item it sends $\widehat{D}_j$ and $\widehat{A}_j$ into $\widetilde{D}_j$ and $\widetilde{A}_j$ respectively.
\end{itemize}

\subsubsection{Unclasping of a crossing}

An {\it unclasping of a crossing} on a ribbon presentation is defined as a replacement of a crossing $[U,U\cap (D_j\cup B_j),U\cap B_k]$ in Figure~\ref{fig:crossing}(i) with another triple $[U,U\cap (D_j'\cup B_j),U\cap B_k]$ in Figure~\ref{fig:crossing}(ii) keeping near $\partial U$ unchanged. It is easy to see that any ribbon presentation can be made into the one without ribbon singularities by a sequence of unclaspings. We can also define an unclasping of a crossing on a long ribbon $n$-knot by applying the associating procedure to the both sides in Figure~\ref{fig:crossing} keeping near the boundary unchanged.

\subsubsection{Connected sum}
A connected sum of two long $n$-knots $K\#L$ is defined up to isotopy by arranging $K$ and $L$ along the standard plane $\R^n\subset\R^{n+2}$.

\subsection{Finite type invariants of long ribbon $n$-knots}

First we recall the notion of $k$-schemes defined in \cite{HKS} which is in some sense a higher dimensional analogue of singular knots in Vassiliev's theory of knot invariants. 

Let $P=\calB\cup\calD$ be a based ribbon presentation and $\{c_1,\ldots, c_k\}$ be a set of disjoint crossings on $P$. A {\it $k$-scheme} $[P;c_1,\ldots,c_k]$ is defined by 
\[ [P;c_1,\ldots,c_k]=\sum_{S\subset\{c_1,\ldots,c_k\}}(-1)^{|S|}[K_{P^S}]\]
where $|S|$ is the size of $S$ and $P^S$ is the ribbon presentation obtained from $P$ by unclasping at all the crossings whose labels are in $S$.

Let $\calJ_k$ be the subspace of $\calL^n\eqdef\span_\R\{\mbox{isotopy classes of long ribbon $n$-knots}\}$ spanned by all $k$-schemes. It is easy to check that this constitutes a descending filtration on $\calL^n$:
\begin{equation}\label{eq:filt-J}
 \calL^n=\calJ_0\supset \calJ_1\supset \calJ_2\supset\cdots\supset \calJ_k\supset\cdots. 
 \end{equation}

Let $g$ be an $\R$-valued invariant of long ribbon $n$-knots. Then we can extend $g$ naturally to $\calL^n$ by linearity. We say that $g$ is {\it of type $k$} if it vanishes on $\calJ_{k+1}$, or equivalently, if $g$ is an element of $(\calL^n/\calJ_{k+1})^*$.

There is a filtration on the set $\calI$ of all $\R$-valued finite type invariants:
\begin{equation}\label{eq:filt-I}
 \calI_0\subset\calI_1\subset\calI_2\subset\cdots\subset\calI_k\subset\cdots =\calI
\end{equation}
with $\calI_k$ being the set of all type $k$ invariants. The two filtrations (\ref{eq:filt-J}) and (\ref{eq:filt-I}) are dual to each other in the sense that there is an isomorphism
\begin{equation}\label{eq:I=J}
 \calI_k/\calI_{k-1}\cong (\calJ_k/\calJ_{k+1})^*. 
\end{equation}
The dimension of $\calI_k/\calI_{k-1}$ is equal to the number of possible sequences $(m_1,\ldots,m_r)$ of positive integers with $m_1+\cdots+m_r=k, 2\leq m_1\leq\cdots\leq m_r$. These facts immediately follow from the results in \cite{HS}. For $g\in\calI_k$, we call its projected image in $\calI_k/\calI_{k-1}$ a {\it principal term of $g$}. Then (\ref{eq:I=J}) says that the principal term is determined by a linear functional on $J_k/J_{k+1}$.

\begin{figure}
\cfig{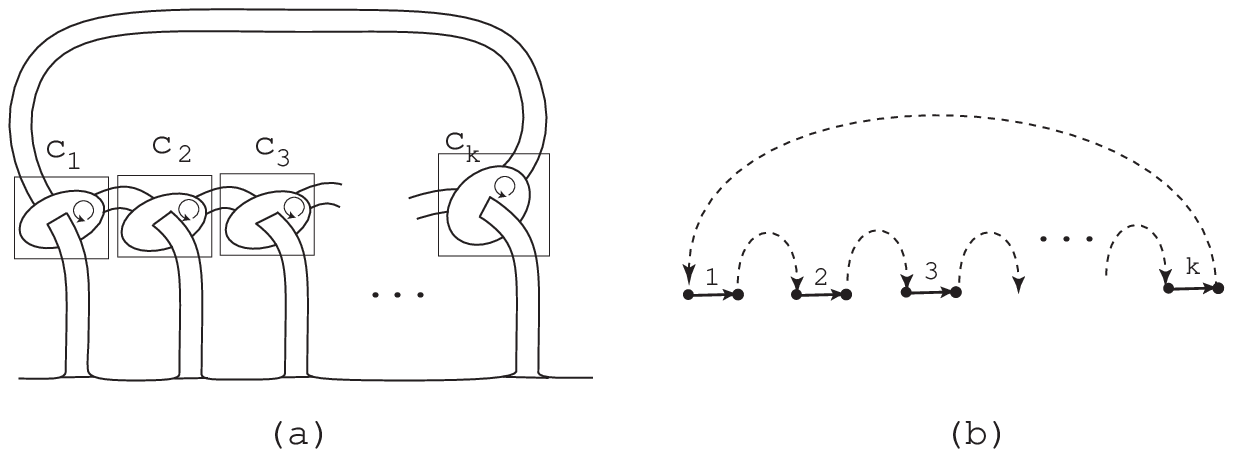}
\caption{}\label{fig:wheel}
\end{figure}

Let $\calP(\calJ_k/\calJ_{k+1})$ be the subspace of $\calJ_k/\calJ_{k+1}$ spanned by elements which cannot be written in $\calJ_k/\calJ_{k+1}$ as a connected sum of two schemes. According to \cite{HS}, $\calJ_k/\calJ_{k+1}$ is generated by connected sums of schemes of the form $[W_i;c_1,\ldots,c_i]$, which is defined in Figure~\ref{fig:wheel}(a). It follows that $\calP(\calJ_k/\calJ_{k+1})$ is one dimensional and spanned by $[W_k;c_1,\ldots,c_k]$.

\begin{Rem}\label{rem:long}
Here, although the set of ribbon $n$-knots and the set of long ribbon $n$-knots are different, it can be shown that the filtrations of finite type invariants of both are the same\footnote{For instance, any ``long" $k$-scheme may be reduced modulo long $(k+1)$-schemes to ``wheel-like" ones (see Figure~\ref{fig:wheel}). So the graded piece $\span_\R\{\mbox{long $k$-schemes}\}/\span_\R\{\mbox{long $(k+1)$-schemes}\}$ may be isomorphic to the one of usual (non-long) $k$-schemes.}.
\end{Rem}

We set $\hat{z}_k=z_k$ when $n$ is odd and let $\hat{z}_3$ be as in Theorem~\ref{thm:invariant} when $n=2$. The main theorem of this paper is stated as follows.
\begin{Thm}\label{thm:Zuniversal}
If $n$ is an odd integer $\geq 3$ and $k\geq 2$, or $(n,k)=(2,3)$, then $\hat{z}_k$ is a finite type invariant of type $k$. Its principal term corresponds to the dual of $[W_k;c_1,\ldots,c_k]$, the basis of the one dimensional subspace $\calP(\calJ_k/\calJ_{k+1})\subset \calJ_k/\calJ_{k+1}$, if $n\nequiv k\ \mbox{mod 2}$, and is zero if $n\equiv k\ \mbox{mod 2}$.
\end{Thm}

\subsection{Proof of Theorem~\ref{thm:Zuniversal}}\label{ss:proof-thm}
\subsubsection{Star-like ribbon presentations}\label{ss:assoc}

The determination of the principal term of $\hat{z}_k$ is reduced to the computations of its values at star-like $k$-schemes, defined as follows.

\begin{figure}
\fig{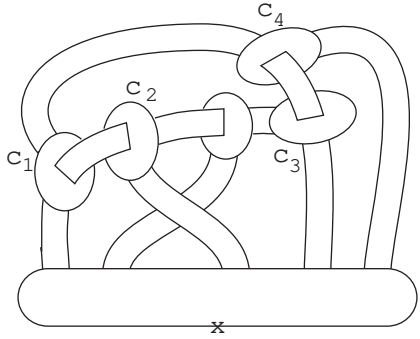}
\caption{star-like $4$-marked ribbon presentation}\label{fig:ex-starlike}
\end{figure}

Let $\{P;c_1,\ldots,c_k\}$ be a combination of a ribbon presentation $P=\calD\cup\calB$ and a set of crossings $\{c_1,\ldots, c_k\}$ on it. We will call such a combination a {\it $k$-marked ribbon presentation}. Let $D_0\in\calD$ be the based disk and $D_j$ be the disk entirely included in $c_j$. We will say that $\{P;c_1,\ldots,c_k\}$ is {\it star-like} if, for each $D_j$, there is a band $B_j$ connecting $D_0$ and $D_j$, and no other band is incident to $D_j$. An example of a star-like marked ribbon presentation is shown in Figure~\ref{fig:ex-starlike}. Also we say that a $k$-scheme associated to a marked star-like ribbon presentation is star-like. The following fact implies that the class of star-like $k$-schemes is general enough to study the space of finite type invariant.
\begin{Prop}[\cite{HS}]\label{prop:HSthm}
Any $k$-scheme can be transformed into a sum of star-like $k$-schemes.
\end{Prop}

We will call the disk $D_j\cup B_j$ a {\it branch} of a marked star-like ribbon presentation. Of course this depends on the choice of a ribbon presentation. 

\subsubsection{Hypotheses}
 
Assume that $n$ is an odd integer $\geq 3$ and $k\geq 2$ or $(n,k)=(2,3)$. Theorem~\ref{thm:Zuniversal} will be proved by computing the values of $\hat{z}_k$ for star-like $m$-schemes $[P;c_1,\ldots,c_m]$ with $m\geq k$. Let 
\[ \psi^S(P;c_1,\ldots,c_m):\R^n\to \R^{n+2}\ (S\subset\{1,\ldots,m\}) \]
be a choice of an embedding for each $[K_{P^S}]$ which is standard near $\infty$. We will often write $\psi^S$ for $\psi^S(P;c_1,\ldots,c_m)$ when the marked ribbon presentation $\{P;c_1,\ldots,c_m\}$ considered is understood from the context.

\begin{figure}
\fig{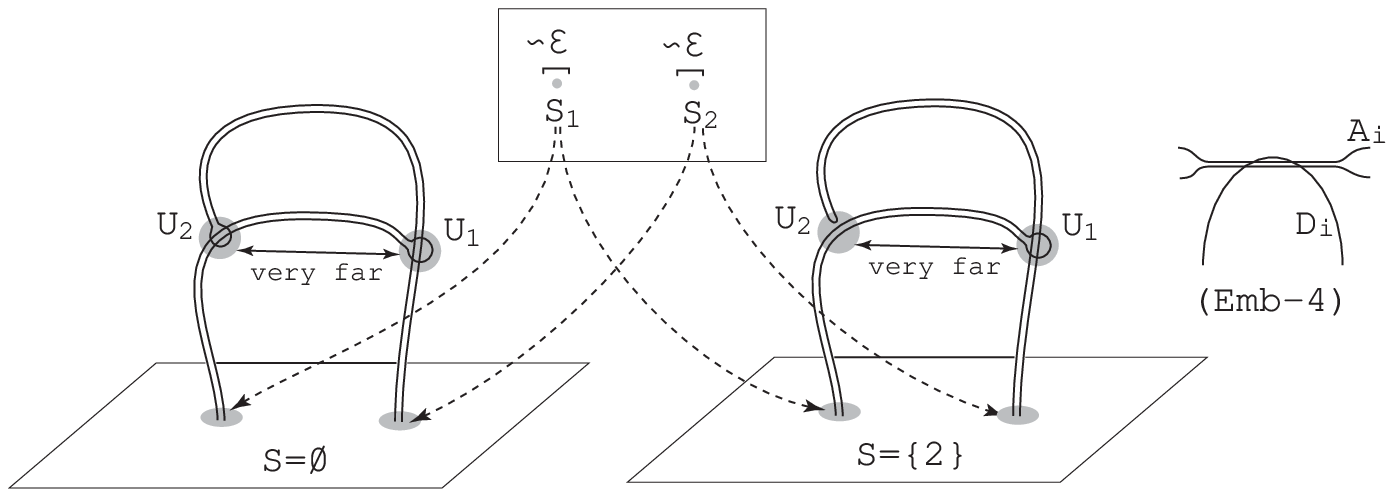}
\caption{}\label{fig:assum}
\end{figure}

Now we make some assumptions throughout the proof, which do not lose the generality. Let $[\breve{U}_i,\breve{D}_i,\breve{A}_i]$ $\ (i=1,\ldots,m)$ denote the crossings on $P$ corresponding to $c_1,\ldots,c_m$ and $(U_i, D_i, A_i)\ (i=1,\ldots,m)$ denote the corresponding crossings on $\Im{\psi^{\emptyset}}$. We assume the conditions in Assumption~\ref{ass:emb} for the embeddings. Some of these assumptions are made depending on small parameters $\varepsilon>0,\ \varepsilon_i>0\ (i=1,\ldots,m)$ with $\varepsilon_i<\varepsilon$. Label the branches on $\{P;c_1,\ldots,c_m\}$ by positive integers $\{1,\ldots,m\}$. 
\begin{Ass}\label{ass:emb}
There are disjoint subsets $\sfD_j$ and $\sfA_j\ (j=1,\ldots,m)$ of $\R^n$ each diffeomorphic to $D^n$ and $S^{n-1}\times I$ respectively such that
\begin{description}
\item[(Emb-0)] $\psi^S=\psi^{\emptyset}$ on $\R^n\setminus(\bigcup_j\sfD_j\cup\bigcup_j\sfA_j)$.
\item[(Emb-1)] $U_i\cap\Im{\psi^S}=\psi^S(\sfD_i\cup \sfA_i)$.
\item[(Emb-2)] If the intersection of the $j$-th branch of $P$ and $\breve{U}_1\cup\cdots\cup\breve{U}_m$ is $\breve{A}_{j_1}\cup\cdots\cup\breve{A}_{j_{r-1}}\cup\breve{D}_{j}$, then 
\[ \sfA_{j_1}\cup\cdots\cup \sfA_{j_{r-1}}\cup \sfD_{j}\subset \sfS_j \]
where
\[ \sfS_j\eqdef\{(x_1,\ldots,x_n)\in\R^n\,|\,\|x_1-j\|^2+\|x_2\|^2+\cdots+\|x_n\|^2\leq\varepsilon^2\}\subset\R^n. \]
\item[(Emb-3)] The distance between the crossings $U_i$ and $U_j$ for $i\neq j$ is very large relative to the diameters of both $U_i$ and $U_j$. More precisely, the distance is assumed larger than $\frac{1}{\varepsilon}\max\{\mathrm{diam}\ U_i,\mathrm{diam}\ U_j\}$. 
\item[(Emb-4)] For any $S$ and $i\in\{1,\ldots,m\}\setminus S$, $\psi^S$ and $\psi^{S\cup\{i\}}$ is chosen so that they coincide outside an $(n+1)$-ball in $U_i$ with radius $\varepsilon_i<\varepsilon$. Indeed, such $\psi^S$ and $\psi^{S\cup\{i\}}$ may be obtained by contracting around the center $\{\frac{1}{2}\}\times S^{n-1}\subset I\times S^{n-1}$ of $A_i$ into a very thin cylinder with small $S^{n-1}$ component and let them approach near $D_i$.
\end{description}
\end{Ass}
Figure~\ref{fig:assum} is a picture explaining these assumptions. One may check that these assumptions are compatible. Note that (Emb-0), (Emb-3) and (Emb-4) are the assumptions concerning the image of embeddings and (Emb-1) and (Emb-2) are those concerning the choice of parameterization of $\R^n$. Note also that we can not define a limit of the embedding at $\varepsilon=0$ while we can define for the limit $\varepsilon_i=0$.

\begin{figure}
\fig{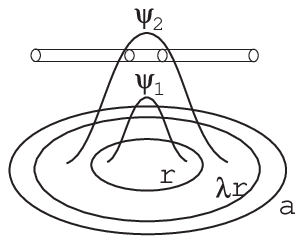}
\caption{}\label{fig:frmchange}
\end{figure}

With all these assumptions, we additionally make the assumption in the following lemma in the case $n=2$ to make the correction term uniformly over the terms of a $k$-scheme. From the procedure of framing correction in this case (see \S\ref{ss:anomaly-n=2}), the correction term depends only on the map 
\[ f= d\psi:\R^2\to I_2(\R^4) \]
giving the tangent 2-framing in $\R^4$ of $\psi(x)\in\R^4\ (x\in\R^2)$ where $I_2(\R^4)$ is the space of all linear injective maps $\R^2\to\R^4$, or the set of 2-framings in $\R^4$. The choice of the embeddings in the following lemma can be made compatible to (Emb-0)--(Emb-4).
\begin{Lem}\label{lem:frmchange}
Let $P_1$ and $P_2$ be two ribbon presentations which differ only by one unclasping of a crossing. For long ribbon 2-knots, we can choose certain embeddings $\psi_i:\R^2\to\R^4,\ i=1,2$ for the isotopy classes $[K_{P_1}]$ and $[K_{P_2}]$ respectively, so that their tangent 2-frames $f_i$ coincide.
\end{Lem}
\begin{proof}
We can assume without losing generality that there exist some positive real numbers $r$ and $a$ with $0<r<a$ such that $f_1(x)$ is equal to a constant $g\in I_2(\R^4)$ in $r\leq ||x||\leq a$ and such that $f_1(x)$ is not a constant function in $||x||<r$. Then we can take $\psi_2$ as
\[ \psi_2(x)=\left\{\begin{array}{ll}
	\lambda\psi_1(\lambda^{-1}x) & ||x||\leq \lambda r \\
	\psi_1(x) & ||x||>\lambda r
	\end{array}\right. \]
for some constant $\lambda$ with $0<\lambda<a/r$. See Figure.~\ref{fig:frmchange}(a) for an explanation of this condition. The naturally induced framing $f_2$ for $\psi_2$ from $f_1$ is
\[ f_2(x)=\left\{\begin{array}{ll}
	f_1(\lambda^{-1}x) & ||x||\leq \lambda r \\
	f_1(x) & ||x||> \lambda r
	\end{array}\right. \]
After some change of parametrization in the disk $||x||\leq a$, $f_2$ becomes equal to $f_1$ without changing the image of $\psi_2$.
\end{proof}

We define the formal alternating sum:
\[ \psi(P;c_1,\ldots,c_m)\eqdef\sum_{S\subset\{1,\ldots,m\}}(-1)^{|S|}\psi^S. \]
Then its configuration space integral $I(\Gamma)$ is given by 
\[ I(\Gamma)(\psi(P;c_1,\ldots,c_m))=\sum_{S\subset\{1,\ldots,m\}}(-1)^{|S|}I(\Gamma)(\psi^S). \]

\subsubsection{Proof in four steps}

Theorem~\ref{thm:Zuniversal} is a consequence of the following four claims.
\begin{description}
\item[Claim 1] {\it %
$\hat{z}_k(\psi(P;c_1,\ldots,c_m))=0$ for any $m$-scheme $[P;c_1,\cdots,c_m]$ with $m>k$. Namely, $\hat{z}_k$ is a finite type invariant of type $k$.
}\par\vspace{2mm}
\item[Claim 2] {\it %
$\hat{z}_k(\psi(P_1;c_1,\ldots,c_i)\# \psi(P_2;d_1,\ldots,d_j))=0$ for any $i$-scheme $[P_1;c_1,\ldots,c_i]$ and $j$-scheme $[P_2;d_1,\ldots,d_j]$ with $i+j=k,\ i,j>0$.
}\par\vspace{2mm}
\item[Claim 3] {\it %
If $n\nequiv k\ $mod 2, $\hat{z}_k(\psi(W_k;c_1,\ldots,c_k))=w_k(\Gamma_k)$ where $W_k$ and the crossings $c_1,\ldots,c_k$ on it are defined in Figure~\ref{fig:wheel}(a) and $\Gamma_k$ is defined in Figure~\ref{fig:wheel}(b). 
}\par\vspace{2mm}
\item[Claim 4] {\it %
If $n\equiv k$ mod 2, $\hat{z}_k(\psi(W_k;c_1,\ldots,c_k))=0$. 
}\par\vspace{2mm}
\end{description}

The reason for Claim 2--4 to be sufficient to determine the principal term of $z_k$ is as follows. Claim 2 implies that $\hat{z}_k$ can be non trivial only on the subspace $\calP(\calJ_k/\calJ_{k+1})$. So the principal term of $\hat{z}_k$ reduces to $\calP(\calJ_k/\calJ_{k+1})^*$. Thus Claim 3 and 4 characterize the element of $\calP(\calJ_k/\calJ_{k+1})^*$ for $\hat{z}_k$. 

Proofs of the claims are done by explicit computations of $I(\Gamma)(\psi(P;c_1,\ldots,c_m))$ for $\deg\Gamma=k,\ m\geq k$.\par\vspace{2mm}

\noindent{\bf Proof of Claim 1.}\ 
Let $\pi^S:C_{\Gamma}(\psi^S)\to C_q(\R^n)$ be the natural projection where $C_q(\R^n)$ denotes the configuration space for external vertices in $\Gamma$. Then $\pi^S$ can be considered as a bundle and we can rewrite the integral by the pushforward:
\[ I(\Gamma)(\psi^S)=\int_{C_{\Gamma}(\psi^S)}\omega(\Gamma)(\psi^S)=\int_{C_q(\R^n)}\pi^S_*\omega(\Gamma)(\psi^S) \]
for some $qn$-form $\pi^S_*\omega(\Gamma)(\psi^S)$ on $C_q(\R^n)$. By using this expression, we can rewrite the alternating sum $I(\Gamma)(\psi(P;c_1,\ldots,c_m))$ as an integral of an alternating sum of forms over the common space $C_q(\R)$:
\[ I(\Gamma)(\psi(P;c_1,\ldots,c_m))=\int_{C_q(\R^n)}\sum_{S\subset \{1,\ldots,m\}}(-1)^{|S|}\pi^S_*\omega(\Gamma)(\psi^S). \]

Let $V_1(D_i)$ be the subset of $C_q(\R^n)$ consisting of configurations such that no points are mapped into $D_i\subset U_i$ and $V_1(A_i)$ is similarly defined for $A_i\subset U_i$. Then we have
\begin{equation}\label{eq:I(V2)}
 \int_{V_1(D_i)\cup V_1(A_i)}\sum_{S\subset \{1,\ldots,m\}}(-1)^{|S|}\pi^S_*\omega(\Gamma)(\psi^S)=O(\varepsilon_i).
 \end{equation}
Here $\varepsilon_i<\varepsilon$ is the number taken in (Emb-4) and $O(\varepsilon_i)$ denotes a term which vanishes at the limit $\varepsilon_i\to 0$. This is because the integral of the LHS is equal to 
\[ \sum_{S\subset\{1,\ldots,\hat{i},\ldots,m\}}(-1)^{|S|}\int_{V_1(D_i)\cup V_1(A_i)}(\pi^S_*\omega(\Gamma)(\psi^S)-\pi^{S\cup\{i\}}_*\omega(\Gamma)(\psi^{S\cup\{i\}}))
\]
and by (Emb-4), this contributes as $O(\varepsilon_i)$. 

By (\ref{eq:I(V2)}), the computation reduces to the one for configurations such that both $D_i$ and $A_i$ includes at least one external vertices. If $m>k$, this is impossible. So we have $\hat{z}_k(\psi(P;c_1,\ldots,c_m))=O(\varepsilon)$ and this must be zero because $\hat{z}_k$ is an invariant by Theorem~\ref{thm:invariant}. Claim 1 is proved.
\qed\par\vspace{2mm}

Assume $m=k$ to prove Claim 2--4.\par\vspace{2mm}

\noindent{\bf Proof of Claim 2.}\ 
According to the observations in the proof of Claim 1, the Jacobi diagrams giving non-zero contribution of $I(\Gamma)(\psi(P;c_1,\ldots,c_k))$ are those without internal vertices and so we restrict to such Jacobi diagrams in the following. Note that $C_{\Gamma}=C_q(\R^n)$ for such $\Gamma$. Let $C_q'(\R^n)=C_q(\R^n)\setminus\bigcup_i (V_1(D_i)\cup V_1(A_i))$. 

\begin{Lem}\label{lem:long-theta-edge}
Let $V_2(e)$ ($e$: $\theta$-edge of $\Gamma$) be the set of configurations for $\Gamma$ such that the two ends of $e$ are mapped by the embeddings into some pair of two different balls $U_i$ and $U_j$ respectively. Then 
\[ \int_{V_2(e)}\sum_{S\subset\{1,\ldots,m\}}(-1)^{|S|}\omega(\Gamma)(\psi^S)=O(\varepsilon). \]
\end{Lem}
\begin{proof}
By the assumption (Emb-3), the locus of the image of the relative vector connecting a point in $U_i$ and another point in $U_j\ (i\neq j)$ via the Gauss map $u:\R^{n+2}\setminus\{0\}\to S^{n+1}$ (defined in \S\ref{ss:def-integral}) for the $\theta$-edge $e$, is included in an arbitrarily small ball embedded into $S^{n+1}$, as $\varepsilon$ tends to $0$. So the integral is $O(\varepsilon)$.
\end{proof}
According to Lemma~\ref{lem:long-theta-edge}, only the Jacobi diagrams each $\theta$-edge of which is entirely mapped into some $U_i$ contributes. Since $\Gamma$ is of degree $k$, it follows that such diagrams can not have internal vertices.

\begin{Lem}\label{lem:long-eta-edge}
Let $V_3(e)$ ($e$: $\eta$-edge of $\Gamma$) be the subset of $C_q'(\R^n)$ consisting of configurations for $\Gamma$ such that two successive external vertices connected by $e$ split into $\sfS_i$ and $\sfS_j\ (i\neq j)$. Then 
\[ \int_{V_3(e)}\sum_{S\subset\{1,\ldots,m\}}(-1)^{|S|}\omega(\Gamma)(\psi^S)=O(\varepsilon). \]
\end{Lem}
\begin{proof}
By (Emb-2), the locus of the image of the relative vector connecting the point in $\sfS_i$ and another point in some $\sfS_j\ (i\neq j)$ via the Gauss map $u':\R^{n}\setminus\{0\}\to S^{n-1}$ (defined in \S\ref{ss:def-integral}) for the $\eta$-edge $e$, is included in an arbitrarily small ball in $S^{n-1}$ as $\varepsilon$ tends to 0, by the definition of $\sfS_i$'s. So the integral is $O(\varepsilon)$.
\end{proof}

For $k$-schemes of the form of connected sum $K_i\# K_j$ of $i$- and $j$-schemes with $i+j=k,\ i,j>1$, any configurations not in $V_2(e)$'s associated to a connected Jacobi diagram, i.e., each $\theta$-edge is mapped entirely into some $U_i$, must be in some $V_3(e)$ because the image of any $\theta$-edge cannot connect $K_i$ and $K_j$ and then some $\eta$-edge must connect them. Hence Claim 2 is proved.
\qed\par\vspace{2mm}

\noindent{\bf Proof of Claim 3.}
Now we compute the precise value of $\hat{z}_k(\psi(W_k;c_1,\ldots,c_k))$. By Lemma~\ref{lem:long-theta-edge} and \ref{lem:long-eta-edge}, the only contributing term in the sum (\ref{eq:Z}) is the term:
\[ \frac{1}{2}\frac{I(\Gamma_k)(\psi(W_k;c_1,\cdots,c_k))w_k(\Gamma_k)}{|\Aut \Gamma_k|}. \]
By Lemma~\ref{lem:IG=Aut} below, this equals 
\[ \frac{|\pAut{\Gamma_k}|w_k(\Gamma_k)}{2|\Aut{\Gamma_k}|}+O(\varepsilon)=\frac{2|\Aut{\Gamma_k}|w_k(\Gamma_k)}{2|\Aut{\Gamma_k}|}+O(\varepsilon)=w_k(\Gamma_k)+O(\varepsilon) \]
where $\pAut \Gamma$ denotes the group of automorphisms of $\Gamma$ considered as an unoriented graph. Roughly, by the above observations, the computation in Lemma~\ref{lem:IG=Aut} reduces to the integral over a direct product of some simple spaces and then it is computed as a product of `mapping degrees'. This completes the proof of Claim 3. \qed\par\vspace{6mm}

\noindent{\bf Proof of Claim 4.}
In the case $n\equiv k$ mod 2, we always have $I(\Gamma_k)=0$ by Proposition~\ref{prop:vanish-axial-symm} and thus Claim 4 is proved.\qed
 
\begin{Lem}\label{lem:IG=Aut}
Under the assumptions (Emb-0)--(Emb-4), we have $ I(\Gamma_k)(\psi(W_k;c_1,\cdots,c_k))=|\pAut\Gamma_k|+O(\varepsilon).$
\end{Lem}
\begin{proof}
Then there are exactly $|\pAut{\Gamma_k}|$ connected components in the reduced configuration space $C_q'(\R^n)\setminus (\bigcup_e V_2(e)\cup \bigcup_f V_3(f))$ for $\Gamma_k$ each of which is equal up to permutations to
\[ \sfM_k\eqdef \sfA_1\times \sfD_1\times \sfA_2\times \sfD_2\times \cdots \times \sfA_k\times \sfD_k \]
By symmetry, all $|\pAut \Gamma_k|$ components contributes as a common value up to sign.
 Note that in the computation of $I(\Gamma_k)(\psi(W_k;c_1,\cdots,c_k))$, the integration domain $\sfM_k$ may be assumed to be common for all $S$  because of (Emb-0) and (Emb-1) while the integrand form varies depending on $\psi^S$'s.

After a suitable $\mathrm{Diff}_+(\R^n)$ action on $\R^n$ fixing 
\begin{itemize}
\item near infinity,
\item outside $\sfS_1\cup\cdots\cup \sfS_k\subset\R^n$ and 
\item the image of the embeddings, 
\end{itemize}
we may assume in addition to (Emb-0)--(Emb-4) that\par\vspace{2mm}
\hspace{3mm}{\bf (Emb-5):}
\[\begin{split}
	 \sfA_i&=\Bigl\{(x_1,\ldots,x_n)\in\R^n\,\Bigl|\Bigr.\,\Bigl(\frac{\varepsilon}{2}\Bigr)^2\leq \|x_1-(i+1)\|^2+\|x_2\|^2+\cdots+\|x_n\|^2\leq \Bigl(\frac{2\varepsilon}{3}\Bigr)^2\Bigr\}\\
	 &\mbox{for }1\leq i<k,\\
	 \sfA_k&=\Bigl\{(x_1,\ldots,x_n)\in\R^n\,\Bigl|\Bigr.\,\Bigl(\frac{\varepsilon}{2}\Bigr)^2\leq \|x_1-1\|^2+\|x_2\|^2+\cdots+\|x_n\|^2\leq \Bigl(\frac{2\varepsilon}{3}\Bigr)^2\Bigr\},\\
	 \sfD_i&=\Bigl\{(x_1,\ldots,x_n)\in\R^n\,\Bigl|\Bigr.\,\|x_1-i\|^2+\|x_2\|^2+\cdots+\|x_n\|^2\leq \varepsilon^4\Bigr\}.
	 \end{split}
 \]
This additional assumption can be made so as not to affect the assumptions (Emb-0)--(Emb-4) by deforming all $\psi^S$'s for $S\subset\{1,\ldots, k\}$ simultaneously. In the case of $n=2$, $\sfS_j, \sfD_j, \sfA_j$ look like the following picture.
\cfig{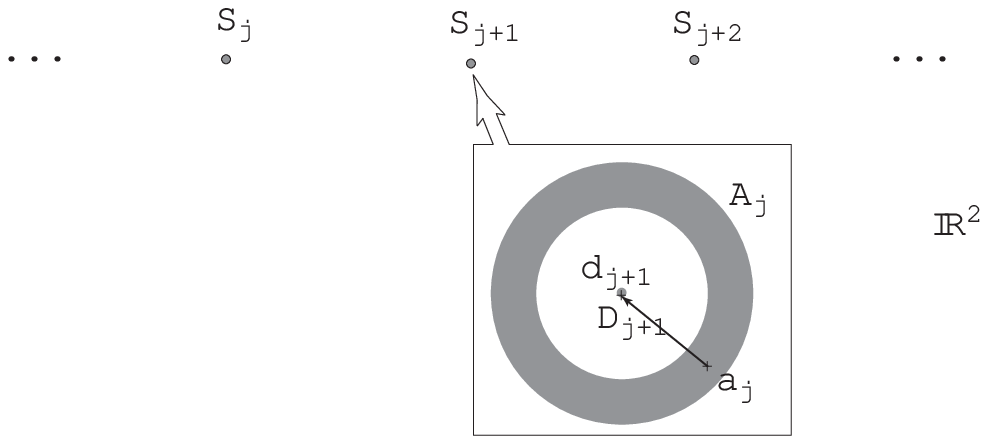}

Let $(a_1,d_1,\cdots,a_k,d_k)\ (a_i\in \sfA_i, d_i\in \sfD_i)$ be the coordinate of $\sfM_k$. The explicit expression for the LHS integral restricted to $\sfM_k$ is
\begin{equation}\label{eq:I(GW)}
 \sum_{S\subset\{1,\ldots,k\}}(-1)^{|S|}\int_{\sfM_k}\bigwedge_i
	\Bigl(u(\psi^S a_i-\psi^S d_i)^*\omega_{n+1}\wedge u'(d_{i+1}-a_i)^*\omega_{n-1}\Bigr).
\end{equation}
It suffices to show that this quantity equals $1+O(\varepsilon)$ because in the case $n\nequiv k$ mod 2, any automorphism on $\sfM_k$ induced from an automorphism of $\Gamma_k$ does not change the orientation of the configuration space and the integrand form. So all $|\pAut\Gamma_k|$ connected components contribute as a common value up to $O(\varepsilon)$.

From the fact: 
\begin{quote}
	varying $d_i$ fixing all other points only affects (modulo $O(\varepsilon)$) to $u(\psi^S a_i-\psi^S d_i)^*\omega_{n+1}$ in $\omega(\Gamma_k)$ because the other form $u'(d_i-a_{i-1})^*\omega_{n-1}$, which may also depend on $d_i$, contributes as $O(\varepsilon)$ by (Emb-5),
\end{quote}
we can write (\ref{eq:I(GW)}) as
\[
 	\sum_{S\subset\{1,\ldots,k\}}
	(-1)^{|S|}\int_{\sfM_k}\bigwedge_i \Bigl(u(\psi^S a_i-\psi^S d_i)^*\omega_{n+1}\wedge u'(d_{i+1}^0-a_i)^*\omega_{n-1}\Bigr)+O(\varepsilon)
\]
where $d_{i+1}^0=(i+1,0,\cdots,0)\in\R^n$ for $1\leq i<k$, $d_{k+1}^0=(1,0,\cdots,0)$. 
We will write $\Phi_i(\psi)\eqdef u(\psi a_i-\psi d_i)^*\omega_{n+1}\wedge u'(d_{i+1}^0-a_i)^*\omega_{n-1}$ for simplicity bearing in mind that it depends on $a_i$ and $d_i$. We have
\[ \int_{\sfM_k}\bigwedge_i\Phi_i(\psi^S)-\int_{\sfM_k}\bigwedge_i\Phi_i(\psi^{S\cup\{j\}})
	=\int_{\sfM_k(j)}\bigwedge_{i\neq j}\Phi_i(\psi^S)
		\int_{\sfA_j\times \sfD_j}\Bigl(\Phi_j(\psi^S)-\Phi_j(\psi^{S\cup \{j\}})\Bigr)
 \]
$(S\subset\{1,\ldots,\hat{j},\ldots,k\})$ where 
\[ \sfM_k(j)\eqdef \sfA_1\times \sfD_1\times\cdots\times\widehat{\sfA_j\times \sfD_j}\times\cdots\times \sfA_k\times \sfD_k, \]
because $\Phi_i(\psi^S)=\Phi_i(\psi^{S\cup\{j\}})$ on $\sfM_k(j)$ for $i\neq j$. Thus (\ref{eq:I(GW)}) equals
\[ 
	\sum_{S\subset\{1,\ldots,\hat{j},\ldots,k\}}(-1)^{|S|}
	\int_{\sfM_k(j)}\bigwedge_{i\neq j}\Phi_i(\psi^S)
		\int_{\sfA_j\times \sfD_j}\Bigl(\Phi_j(\psi^S)-\Phi_j(\psi^{S\cup \{j\}})\Bigr)+O(\varepsilon).
 \]
By using Lemma~\ref{lem:part-I(GW)} below iteratively, we obtain the desired result.
\end{proof}

\begin{Lem}\label{lem:part-I(GW)}
\[ \int_{\sfA_j\times \sfD_j}\Bigl(\Phi_j(\psi^S)-\Phi_j(\psi^{S\cup \{j\}})\Bigr)=1.\]
\end{Lem}
\begin{proof}
For each crossing $(U_i,D_i,A_i)$, consider the three unclaspings of it defined as follows.
\begin{description}
\item[$(U_i,D_i',A_i)$] The result of an unclasping of $(U_i,D_i,A_i)$ keeping $A_i$ fixed.
\item[$(U_i,D_i,A_i')$] The result of an unclasping of $(U_i,D_i,A_i)$ keeping $D_i$ fixed.
\item[$(U_i,D_i',A_i')$] Obtained by mixing the $D_i'$ and $A_i'$ in the previous two.
\end{description}
Note that these three are isotopic each other.

Denote by $\psi^S(D_j'), \psi^S(A_j'), \psi^S(D_j',A_j')$ the corresponding embeddings for the three which are all isotopic to $\psi^{S\cup \{j\}}$. Since the three are all isotopic relative to $\partial(D_j\times A_j)$ and the integrand form is closed, 
\[
	\int_{\sfA_j\times \sfD_j}\Phi_j(\psi^S(D_j'))
	=\int_{\sfA_j\times \sfD_j}\Phi_j(\psi^S(A_j'))
	=\int_{\sfA_j\times \sfD_j}\Phi_j(\psi^S(D_j',A_j'))
	=\int_{\sfA_j\times \sfD_j}\Phi_j(\psi^{S\cup\{j\}}).
	 \]
Thus the LHS of the lemma is rewritten as
\[ 
	\int_{\sfA_j\times \sfD_j}\Bigl(\Phi_j(\psi^S)-\Phi_j(\psi^{S}(D_j'))\\
	-\Phi_j(\psi^{S}(A_j'))+\Phi_j(\psi^{S}(D_j',A_j'))\Bigr)\\
	=\int_{\sfA_j^{\circ}\times \sfD_j^{\circ}}\Phi_j(\psi_j^{S\circ\circ}). 
\]
Here $\sfX^{\circ}$ denotes $\sfX\cup_{\partial}(-\sfX)$ and $\psi_j^{S\circ\circ}$ is the embedding of $\sfA_j^{\circ}\cup \sfD_j^{\circ}$ into $U_j\cup_{\partial}(-U_j)\cong S^{n+2}$ obtained by gluing $\psi^S(*)$'s at boundaries.
 
The computation of the last expression decomposes with respect to the splitting $\sfA_j^{\circ}\cong S^{n-1}\times S^1$ as follows:
\[\begin{split}
  \int_{(v,l,d_j)\in S^{n-1}\times S^1\times \sfD_j^{\circ}}&
	u(\psi^{S\circ\circ}_jl-\psi^{S\circ\circ}_jd_j)^*\omega_{n+1}
		\wedge u'(d_{j+1}^0-(v,l))^*\omega_{n-1}\\
	=\int_{\sfL_j^{\circ}(v)\times \sfD_j^{\circ}}&
	u(\psi_j^{S\circ\circ}l-\psi_j^{S\circ\circ}d_j)^*\omega_{n+1}
		\int_{v\in S^{n-1}}\omega_{n-1}(v).
	\end{split}
 \]
where we use the fact:
\begin{quote}
 varying $a_i$ along the line
\begin{equation}\label{eq:def-L}
 \sfL_i(v)\eqdef\{(i+1,0,\ldots,0)+tv\,|\,t>0\}\cap \sfA_i\subset\R^n,\ v\in S^{n-1} 
 \end{equation}
fixing all other points only affects (modulo $O(\varepsilon)$) to the form $u(\psi^S a_i-\psi^S d_i)^*\omega_{n+1}$ by (Emb-5).
\end{quote}
 
To compute the last integral, we choose an orientation on $\sfL_j(v)$ given by the direction of $t$ (in the definition of $\sfL_j(v)$) increases. Then the induced orientation on $S^{n-1}$ from that of $S^{n-1}\times \sfL_j(\cdot)=\sfA_j$ and that of $\sfL_j(\cdot)$ is opposite to the one naturally induced from the orientation on $\R^n$ by the outward-normal first convention. 

We consider first the integral along $\sfL_j^{\circ}(v)\times \sfD_j^{\circ}$. Here the induced orientation on $S^{n+1}$ by mapping $\mbox{(orientation on $\sfL_j$)}\wedge\mbox{(orientation on $\sfD_{i}$)}$ via $u$ coincides with the one naturally induced from $\R^{n+2}$. To see this, consider the orientation on a point $(l_j(v),d_j)\in \sfL_j(v)\times \sfD_j$ where
\[ \begin{split}
	l_j(v)&=(p,y_1,y_2,\cdots,y_{n-1})\in\R^3\times[-1,1]^{n-1}\subset\R^{n+2},\\
	d_j&=(q,y_1',y_2',\cdots,y_{n-1}')\in\R^3\times[-1,1]^{n-1}\subset\R^{n+2}
	\end{split}
\]
where $(p,q)$ are points on the boundary of the two components in a crossing of a ribbon presentation. By the definition of the ribbon presentation $W_k$ (Figure~\ref{fig:wheel}(a)) and the orientation on $\sfL_j(v)$, $dp\wedge dq$ is mapped by the Gauss map into the natural orientation on $S^2\subset\R^3$ because the linking number of the two arcs is $1$. The orientation on $\sfL_j(v)\times \sfD_j$ corresponding to varying $p$ fixing $y_1,\ldots,y_{n-1}$ (varying along the orientation on $\sfL_j(v)$) and varying $d_j$, is $dp\wedge dq\wedge dy_1'\wedge\cdots\wedge dy_{n-1}'$. This is mapped by the map $u$ into the natural orientation on $S^{n+1}\subset\R^{n+2}$. Further, the linking number of $\Im\psi_j^{S\circ\circ}|\sfL_j^\circ(v)$ and $D_j^\circ$ is $\pm 1$ because it is equal to the intersection number of $\Im\psi_j^{S\circ\circ}|\sfL_j^\circ(v)$ and an $(n+1)$-disk bounded by $D_j^\circ$. Thus the integral along $\sfL_j^{\circ}(v)\times \sfD_j^{\circ}$ contributes by $1$. 

For the integral along $S^{n-1}$, since the direction of the vector from $v_j$ to $d_{j+1}^0$ coincides with the ingoing normal vector on $S^{n-1}$, the integral contributes by $(-1)^n\times (-1)^n=1$ and the result follows.
\end{proof}

\subsection{Relation with the Alexander polynomial}

The Alexander polynomial ${\Delta}_K(t)\in\Z[t,t^{-1}]$ for a (long) $n$-knot $K$ is defined by using the Fox calculus for the knot group. See  \cite{HKS} for detailed definition. The Alexander polynomial is uniquely determined by the conditions: ${\Delta}_K(1)=1$ and $(d{\Delta}_K/dt)(1)=0$, and we use such a normalized one. Then a series of invariants $\alpha_2, \alpha_3,\ldots$ of (long) $n$-knots are defined by
\begin{equation*}
\log{\Delta}_K(t)|_{t=e^h}=\alpha_2(K)h^2+\alpha_3(K)h^3+\ldots\in\\\Q[[h]].
\end{equation*}

The following result gives a complete correspondence between finite type invariant of ribbon 2-knots and $\alpha_j$-invariants.
\begin{Prop}\label{prop:HS}
Let $k>1$ and $n>1$. Then 
\begin{enumerate}
\item (Habiro-Kanenobu-Shima) $\alpha_k$ is a primitive ($=$additive) invariant of type $k$ of ribbon $n$-knots. 
\item (Habiro-Shima) We have the identification:
\[ \{\mbox{\rm $\Q$-valued finite type invariants}\}=\Q[\alpha_2,\alpha_3,\ldots] \]
of graded Hopf algebras.\footnote{In \cite{HS}, the results are stated for $n=2$ and it is remarked there that their result is generalized to ribbon $n$-knots. Habiro says that the proof for $n>2$ is exactly the same as for $n=2$}. 
\end{enumerate}
\end{Prop}
By Remark~\ref{rem:long}, Proposition~\ref{prop:HS} also holds for long ribbon $n$-knots.

\begin{Cor}\label{cor:z_k-a_k}
Let $n$ be an odd integer $>1$. For long ribbon $n$-knots, 
$\hat{z}_k$ is a degree at most $k$ polynomial in $\alpha_2,\alpha_3,\ldots,\alpha_k\ (\deg\alpha_j=j)$ and we have
\[\hat{z}_k\equiv \left\{
 	\begin{array}{ll}
 		\alpha_k & \mbox{in $\calI_k/\calI_{k-1}$ if $k$ is even},\\
 		0 & \mbox{in $\calI_k/\calI_{k-1}$ if $k$ is odd}.
	\end{array}\right.
\]
For long ribbon 2-knots, $\hat{z}_3=\alpha_3+\lambda\alpha_2\ \ \mbox{for some constant $\lambda\in\R$}$.
\end{Cor}
In particular, Corollary~\ref{cor:z_k-a_k} implies
\begin{Cor}
The BCR invariant is non-trivial.
\end{Cor}

\begin{proof}[Proof of Corollary~\ref{cor:z_k-a_k}]
That $\hat{z}_k$ is a degree at most $k$ polynomial in $\alpha_2,\ldots,\alpha_k$ follows from the fact that $\hat{z}_k$ is of type $k$ by Theorem~\ref{thm:Zuniversal} and the result of \cite{HS}.

The principal term of $\hat{z}_k$ is determined by Theorem~\ref{thm:Zuniversal} and Lemma~\ref{lem:alpha_k} below.
\end{proof}

\begin{Lem}\label{lem:alpha_k}
The type $k$ invariant $\alpha_k$ satisfies the following conditions.
\begin{enumerate}
\item $\alpha_k(\calJ_i\#\calJ_j)=0\ \mbox{ for } i+j=k,\,\, i,j>0.$
\item $\alpha_k([W_k;c_1,\ldots,c_k])=1,$\\ where $W_k$ and the crossings $c_1,\ldots,c_k$ on it are defined in Figure~\ref{fig:wheel}(a).
\end{enumerate}
\end{Lem}
\begin{proof}
(1) is because $\alpha_k$ is additive.

For (2), it suffices to compute the value of $\alpha_k$ at $[K_{W_k}]-1$ because $[K_{W_k^{c_j}}]=1$ and
\[ \begin{split} 
	[W_k;c_1,\ldots,c_k]&=[W_k;c_2,c_3,\ldots,c_k]-[W_k^{c_1};c_2,c_3,\ldots,c_k]\\
	&=[W_k;c_2,c_3,\ldots,c_k]=\cdots=[W_k;c_k]=[K_{W_k}]-1.
	\end{split} \]
From $\Delta_{[K_{W_k}]}(t)=1+(t-1)^k$, we have
\[ \log{\Delta_{[K_{W_k}]}(e^h)}=\log{(1+(e^h-1)^k)}=h^k+o(h^{k}). \]
Hence $\alpha_k([K_{W_k}])=1$ and $\alpha_k([W_k;c_1,\ldots,c_k])=\alpha_k([K_{W_k}-1])=1$.
\end{proof}

\section{Bott-Cattaneo-Rossi classes in $H^*(\mathrm{Emb}(\R^n,\R^m);\R)$}

In this section, along a similar line to \cite{CCL} we show that certain cocycles of $\mathrm{Emb}(\R^n,\R^m)$ for $m,n\geq 3$ odd, defined by configuration space integral are non-trivial in cohomology. Here $\mathrm{Emb}(\R^n,\R^m)$ is the space of long embeddings $\R^n\to\R^m$ which are standard near $\infty$, equipped with the Whitney $C^{\infty}$-topology. We construct a morphism of cohomology groups from $H^{k(m-(n+2))}(\mathrm{Emb}(\R^n,\R^m);\R)$ to the 0-th cohomology group of a certain space, which is identified with the dual of the real vector space spanned by chord diagrams and show that the configuration space integral classes correspond non-trivial 0-dimensional cohomology class on that space. Here a {\it chord diagram} means a Jacobi diagram without internal vertices and without a cycle consisting only of $\eta$-edges. 

\subsection{Singular disk}
We will identify the 0-th cohomology group of the space of singular disks with the dual of the space of chord diagrams below. A singular disk is an immersion of the lower half plane $\R^2_-$ into $\R^l\ (l>3)$ obtained from a marked star-like ribbon presentation as follows. Let $\{P;c_1,\ldots,c_k\}$ be a $k$-marked star-like ribbon presentation without unmarked crossings. Replace the based disk of $P$ with the lower half plane $\R^2_-$. Then modify each marked crossing as follows.
\begin{equation}\label{eq:sing-cross}
 \fig{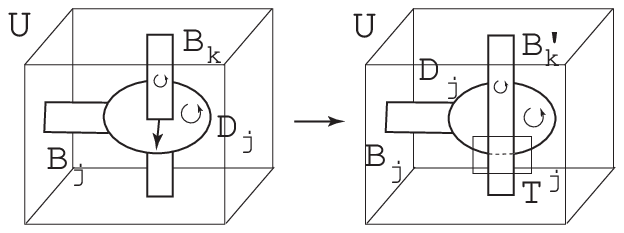}
\end{equation}
Here we assume that $\partial D_j$ and $\partial B_k'$ in the RHS are tangent to some 2-plane $T^j$ near a neighborhood of their intersection. Then embed the resulting immersed disk into $\R^l$ via the inclusion $\R^3\hookrightarrow\R^l$. We call an immersed half-plane into $\R^l$ isotopic to the resulting object obtained now such that 
\begin{itemize}
\item each singularity and boundary arcs near it are included in some $T^j$,
\item two arcs from a band are parallel,
\end{itemize}
a {\it $k$-singular disk}. We denote by $\mathrm{Imm}^r_k(\R^2_-,\R^l)$ the space of all $k$-singular disks. We equip $\mathrm{Imm}(\R^2_-,\R^l)$ with the quotient topology with respect to the projection
\[ \mathrm{Imm}(\R^2,\R^l)\times\mathrm{Emb}(\R^1,\R^2)\to \mathrm{Imm}(\R^2_-,\R^l) \]
and equip $\mathrm{Imm}^r_k(\R^2_-,\R^l)\subset \mathrm{Imm}(\R^2_-,\R^l)$ with the induced topology. Here $\mathrm{Imm}(A,\R^l)\ (A=\R^2\mbox{ or }\R^2_-)$ denotes the space of (long) immersions $A\to \R^l$ which are standard near $\infty$, equipped with the Whitney $C^{\infty}$-topology. However, the particular topology is not so important below. Only the connectivity is important.

\begin{Rem}
One can also define a singular disk by singularizing from the other side of the band instead of (\ref{eq:sing-cross}). However, one can show that the resulting cycle in $\mathrm{Emb}(\R^n,\R^m)$, which will be constructed later, is homologous to the one obtained from the above definition.
\end{Rem}

\subsection{Associating a chord diagram to a singular disk}
Now we see the relationship between $\mathrm{Imm}_k^r(\R^2_-,\R^l)$ and chord diagrams. We define a map 
\[\Gamma(\cdot):\{\mbox{$k$-singular disk}\}\to\calG_k^0\]
associating a chord diagram with each $k$-singular disk $S$ as follows. Choose a star-like $k$-marked ribbon presentation $\{P;c_1,\ldots,c_k\}$ yielding the given $k$-singular disk. Let $D_0$ be the based disk of $P$, $D_1,\ldots,D_k$ be disks included in the crossings $c_1,\ldots,c_k$ respectively, $B_1, \ldots, B_k$ be the bands connecting $D_0$ and $D_1,\ldots,D_k$. Label the marked crossings intersecting $B_j$ by $U_{j1},\ldots,U_{jm_j}$ in the order from $D_0$ on $B_j$. Note that $m_j$ can be zero.

\begin{figure}
\fig{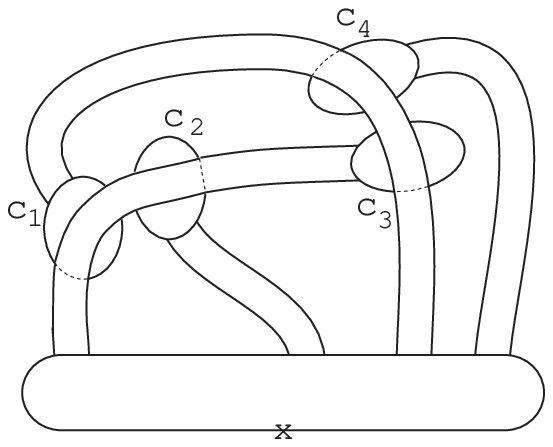}
\caption{}\label{fig:ex-sing}
\end{figure}

With these data, form a graph as follows. 
\begin{enumerate}
\item[\bf Step 1] Put $k$ distinct paths of $\eta$-edges $C_1,\ldots,C_k$ of lengths $m_1,\ldots,m_k$ respectively. Then label each vertex in the head of $C_j$ by $D_j$ and label each of the other vertices by $U_{jp}$ if it is on $C_j$ and is on the $p$-th position from the tail. For example, the result for Figure~\ref{fig:ex-sing} is shown in the following picture:
\[ \fig{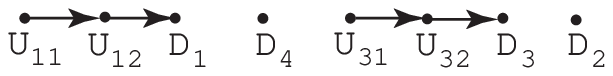} \]
\item[\bf Step 2] Connect a point labeled $D_i$ with a point labeled $U_{jp}$ by a $\theta$-edge with the orientation $(D_i,U_{jp})$ if and only if the intersection of the disk $D_i$ and $P\setminus D_i$ is included in $U_{jp}$.  An example for Figure~\ref{fig:ex-sing} is shown in the following picture:
\[ \fig{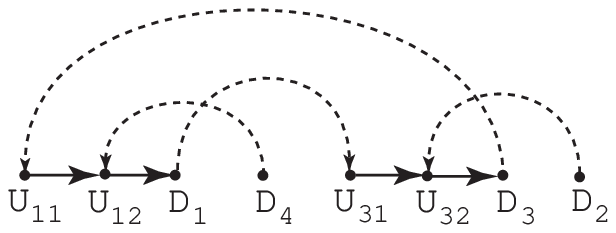} \]
Define $\Gamma(S)$ to be the resulting chord diagram.
\end{enumerate}
By definition, the result does not depend on the choice of a marked ribbon presentation.

\begin{Prop}
If $l>3$, then the connected components of $\mathrm{Imm}_k^r(\R^2_-,\R^l)$ are in one-to-one correspondence with the set of chord diagrams with $k$ chords via the map $\Gamma(\cdot)$.
\end{Prop}
Proof of this proposition is just the same as \cite[Proposition~2.5]{CCL}.

\subsection{Construction of a cycle in $\mathrm{Emb}(\R^n,\R^m)$}
We will construct a morphism
\[ i_k^r:H^{k(m-(n+2))}(\mathrm{Emb}(\R^n,\R^m);\R)\to H^0(\mathrm{Imm}^r_k(\R^2_-,\R^l);\R) \]
for $l-3=m-(n+2)$.

Let $\widetilde{G}_2(\R^l)$ denote the Grassmann manifold of oriented 2-planes in $\R^l$ and let
\[ \rho_k:\mathrm{Imm}_k^r(\R^2_-,\R^l)\to \widetilde{G}_2(\R^l)^{\times k} \]
be the map associating the tangent 2-planes $T^j$ in $\R^l$ for each $j$. Then we define the pullback bundle $\mathfrak{D}_k^r\eqdef \rho_k^*Q_l^{\times k}$ with fiber $(S^{l-3})^{\times k}$ so that the following diagram is commutative.
\[ \begin{CD}
	\mathfrak{D}_k^r @>>> Q_l^{\times k}\\
	@VVV @VVV\\
	\mathrm{Imm}_k^r(\R^2_-,\R^l) @>\rho_{k}>> \widetilde{G}_2(\R^l)^{\times k}
	\end{CD}
\]
Here $Q_l\eqdef SO(l)\times_{SO(2)\times SO(l-2)}S^{l-3}$, the unit normal sphere bundle over $\widetilde{G}_2(\R^l)$. For a fiber $\mathfrak{F}_k^r(\gamma)$ over $\gamma\in\mathrm{Imm}_k^r(\R^2_-,\R^l)$, we construct a `blowing-up' map $\mathfrak{F}_k^r(\gamma)\to \mathrm{Imm}(\R^2_-,\R^l)$ as follows.

\begin{figure}
\fig{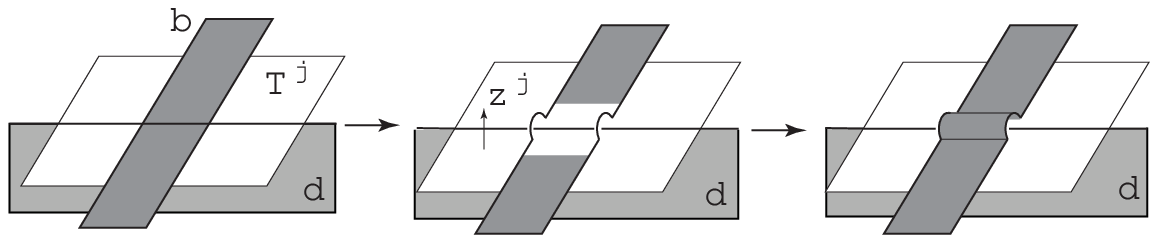}
\caption{}\label{fig:push}
\end{figure}

A point in $\mathfrak{F}_k^r(\gamma)$ can be considered as a sequence of unit normal vectors $(\mathbf{z}^1,\ldots,\mathbf{z}^k)$ to the tangent 2-planes $T^1,\ldots,T^k$. For two transversely intersecting arcs $(c_1,c_2)$ lying on $T^j$ in $\R^l$, a canonical way to blow-up the singularity is given in \cite{CCL} as a family of deformations $((c_1,c_2),\mathbf{z}^j)\mapsto (c_1+\alpha_1^j(\mathbf{z}^j),c_2+\alpha_2^j(\mathbf{z}^j))\ (\mathbf{z}^j\in S^{l-3})$ where $\alpha_{a_j}^j(\mathbf{z}^j)\ (a_j=1,2)$ is a little push along $(-1)^{a_j+1}\mathbf{z}^j$. Here we use a slightly modified version of it, which is equivalent to the original one. Namely, we use a family of deformations $((c_1,c_2),\mathbf{z}^j)\mapsto (c_1+\alpha_1^j(\mathbf{z}^j),c_2)\ (\mathbf{z}^j\in S^{l-3})$. We apply this blowing-up to the three arcs involved in $T^j$, which are parts of boundaries of a band $b$ and a disk $d$ as follows. Here two arcs from $b$ and an arc from $d$ intersect transversely at two double points. We push both of the two arcs from $b$ by $\alpha_1^j(\mathbf{z}^j)$. Then the two arcs from $b$ are still parallel and thus we can refill the interior of the band in a canonical way (see Figure~\ref{fig:push}). 

\begin{figure}
\fig{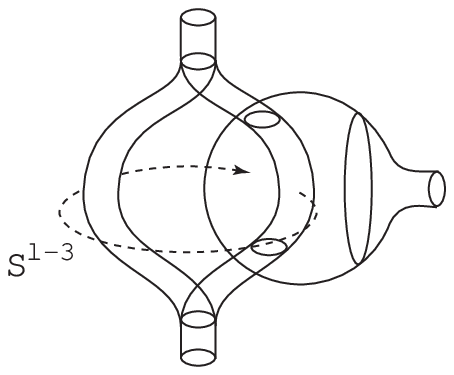}
\caption{}\label{fig:fam-cross}
\end{figure}

Now we have obtained a continuous family $\{P_t\}_t$ of immersed half planes parametrized by $t\in\mathfrak{F}_k^r(\gamma)\cong (S^{l-3})^{\times k}$, which are embeddings on $(S^{l-3}\setminus\{\mathrm{point}\})^{\times k}$ and otherwise non-embedding immersions which are ribbon 2-disks.\footnote{This is because any arc $c$ parallel to the band $b$ may be blown-up to form a $(l-2)$-sphere with two arcs stuck into, such that it intersects the disk $d$ at a point, which is collected into a ribbon singularity if the parallel arc $c$ varies.} By fixing a continuous family of decompositions of ribbon disks into disks and bands as for ribbon presentations, we may associate a family $K_{P_t}$ of long ribbon $n$-knots embedded into $\R^m$. Namely, if $P_t=\calD_t\cup\calB_t$ be the continuous family of decompositions into disks and bands, then we associate an immersed disk
\[ V_{P_t}\eqdef N(\calD_t\times[-2,2]^{n-1}\cup\calB_t\times [-1,1]^{n-1})\subset \R^l\times\R^{n-1}=\R^m \]
where $N(\ )$ denotes a canonical way of smoothing of corners. Then we obtain a family of long ribbon $n$-knots $K_{P_t}\eqdef \partial V_{P_t}$. Figure~\ref{fig:fam-cross} is a picture trying to explain the associated blown-up family of a crossing by using a 2-dimensional knot.

We may also realize the family of long ribbon $n$-knots in $\R^m$ as a family of embeddings $\psi^t:\R^n\to \R^m\ (t\in\mathfrak{F}_k^r(\gamma))$ that are all coincide outside crossings. Here and after a {\it crossing} means the part associated to one of the blown-up family of ribbon 2-disks from the RHS of (\ref{eq:sing-cross}), unlike in \S\ref{s:proof}. We will show in Proposition~\ref{prop:psi-exist} below the existence of such a family of embeddings. Hence we have obtained a $k(l-3)=k(m-(n+2))$-cycle $\sigma^{k(l-3)}=\sigma^{k(m-(n+2))}$ embedded into $\mathrm{Emb}(\R^n,\R^m)$ from the choice of a point in $\mathrm{Imm}_k^r(\R^2_-,\R^l)$, the choice of a disk-band decomposition, and the choice of a family of embeddings $\{\psi^t\}_t$.

Given a $k(m-(n+2))$-cocycle $\zeta^{k(m-(n+2))}$, if the evaluation
\begin{equation}\label{eq:eval}
 \langle (\psi^\bullet)^*\zeta^{k(m-(n+2))},[\mathfrak{F}_k^r(\gamma)]\rangle
	\equiv\langle\zeta^{k(m-(n+2))},[\sigma^{k(m-(n+2))}]\rangle\in\R 
\end{equation}
where $\psi^\bullet\eqdef\{\psi^t\}_t$ is considered as a map $\mathfrak{F}_k^r(\gamma)\to\mathrm{Emb}(\R^n,\R^m)$, does not depend on the choices made and is constant on the connected component of $\mathrm{Imm}_k^r(\R^2_-,\R^l)$, then the correspondence
\[ i_k^r: \zeta^{k(m-(n+2))}\mapsto (\gamma \mapsto \langle (\psi^\bullet)^*\zeta^{k(m-(n+2))},[\mathfrak{F}_k^r(\gamma)]\rangle) \]
descends to a morphism of cohomology groups:
\[ i_k^r:H^{k(m-(n+2))}(\mathrm{Emb}(\R^n,\R^m);\R)\to H^0(\mathrm{Imm}^r_k(\R^2_-,\R^l);\R). \]
Note that a simultaneous isotopy over $\sigma^{k(m-(n+2))}$ does not change its homology class. So it does not change the value (\ref{eq:eval}) too.

Let $[\breve{U}_i,\breve{D}_i,\breve{A}_i]\ (i=1,\ldots,m)$ denote the crossings on the singular star-like ribbon presentation presenting the $k$-singular disk $\gamma$ and $(U_i,D_i,A_i)\ (i=1,\ldots, m)$ denote the corresponding crossings on $K_{P_{t=*}}, *\in\mathfrak{F}_k^r(\gamma)$.

\begin{Prop}\label{prop:psi-exist}
There exits a smooth family of embeddings $\psi^t\ (t\in\mathfrak{F}_k^t(\gamma))$ such that 
\begin{itemize}
\item $\Im{\psi^t}=K_{P_t}$,
\item $\psi^t$'s coincide outside crossings.
\end{itemize}
\end{Prop}
\begin{proof}
First choose a long embedding $\psi^*=\psi^{t=*}:\R^n\to\R^m$ for $K_{P_*}$. Since all $K_{P_t}$'s coincide outside crossings, one can choose a family of embeddings
\[ \psi^t_{\mathrm{out}}:\R^n\setminus\bigcup_i(\sfD_i\cup \sfA_i)\hookrightarrow \R^m \]
such that they coincide for all $t\in\mathfrak{F}_k^r(\gamma)$ and $\Im\psi^t_{\mathrm{out}}=K_{P_t}\cap (\R^m\setminus \cup_i U_i)$, where $\sfD_i\eqdef (\psi^*)^{-1}D_i, \sfA_i\eqdef (\psi^*)^{-1}A_i$. 

Then we extend $\{\psi^t_{\mathrm{out}}\}$ over $\cup_i\sfA_i$. To do so, we ignore the part $\cup_i D_i$ in $K_{P_t}$. Then we can choose a simultaneous isotopy deforming the family of embedded images $U_i\cap K_{P_t}$ into one particular image $U_i\cap K_{P_*}\ *\in\mathfrak{F}_k^r(\gamma)$. By pulling-back $\psi^*$ by the simultaneous isotopy, one obtains a family of embeddings extended to $\R^n\setminus \cup_i\sfD_i$.

Finally we can extend the previous family of embeddings to $\R^n$ by extending in the same way as $\psi^*$ for all $t\in\mathfrak{F}_k^r(\gamma)$ because the embedded images from $\sfD_i$ is constant in $t\in\mathfrak{F}_k^r(\gamma)$ by the construction.
\end{proof}

\subsection{Higher BCR classes}
Let $\zeta_{k(m(n+2))}\in H^{k(m-(n+2))}(\mathrm{Emb}(\R^n,\R^m);\R)$ be defined similarly as $z_k$ by replacing $\R^{n+2}$ with $\R^m$ in $C_{\Gamma}$, $\omega_{n+1}$ with $\omega_{m-1}$, integral over $C_{\Gamma}$ with fiber integration along the $C_{\Gamma}$-fiber. Namely, for $k\geq 2$, 
\begin{equation}\label{eq:zeta}
\zeta_{k(m-(n+2))}(\psi)\eqdef\frac{1}{2}\sum_{\Gamma\in\calG_k^0}\frac{I(\Gamma)(\psi)w_k(\Gamma)}{|\Aut\Gamma|}\in\Omega^{k(m-(n+2))}(\mathrm{Emb}(\R^n,\R^m)).
\end{equation}
where we use the same symbol for both the form and the cohomology class. That $\zeta_{k(m-(n+2))}$ is a cocycle can be proved by exactly the same argument as for the proof of Theorem~\ref{thm:invariant}, only replacing $\R^{n+2}$ with $\R^m$.

As mentioned in the introduction, Budney showed in \cite{Bud} that $\mathrm{Emb}(\R^n,\R^m)$ is $(2m-3n-4)$-connected. If $k\geq 2$, then
\[ k(m-(n+2))\geq 2(m-(n+2)) > 2m-3n-4. \]
So it makes sense to consider $k(m-(n+2))$-cocycles.

\begin{Thm}\label{thm:BCR-nontrivial}
Let $k\geq 2$. If $m, n$ are odd integers $\geq 3$, $m>n+2$, then the BCR class 
\[ \zeta_{k(m-(n+2))}\in H^{k(m-(n+2))}(\mathrm{Emb}(\R^n,\R^m);\R) \]
is non-trivial. Moreover $i_k^r(\zeta_{k(m-(n+2))})$ is a cocycle in $H^0(\mathrm{Imm}^r_k(\R^2_-,\R^l);\R)$. If $k$ is even, then it coincides with the weight function $w_k$ restricted to chord diagrams under the identification $H^0(\mathrm{Imm}^r_k(\R^2_-,\R^l);\R)\cong \span_{\R}\{\mbox{chord diagrams with $k$ chords}\}^*$ and if $k$ is odd, it is zero.
\end{Thm}

\begin{proof}
Let $\psi^t:\R^n\to \R^m\ (t\in\mathfrak{F}_k^r(\gamma))$ be the family of embeddings representing the cycle $\sigma^{k(m-(n+2))}$. As in the proof of Theorem~\ref{thm:Zuniversal}, we may assume (the analogues of) (Emb-0)--(Emb-4) after a suitable simultaneous isotopy over $\sigma^{k(m-(n+2))}$:
\begin{description}
\item[(Emb-0)] $\psi^t=\psi^{*}$ on $\R^n\setminus(\bigcup_j\sfD_j\cup\bigcup_j\sfA_j)$.
\item[(Emb-1)] $U_i\cap\Im{\psi^t}=\psi^t(\sfD_i\cup \sfA_i)$.
\item[(Emb-2)] If the intersection of the $j$-th branch of $P_*$ and $\breve{U}_1\cup\cdots\cup\breve{U}_m$ is $\breve{A}_{j_1}\cup\cdots\cup\breve{A}_{j_{r-1}}\cup\breve{D}_{j}$, then 
\[ \sfA_{j_1}\cup\cdots\cup \sfA_{j_{r-1}}\cup \sfD_{j}\subset \sfS_j \]
where
\[ \sfS_j\eqdef\{(x_1,\ldots,x_n)\in\R^n\,|\,\|x_1-j\|^2+\|x_2\|^2+\cdots+\|x_n\|^2\leq\varepsilon^2\}\subset\R^n. \]
\item[(Emb-3)] The distance between the crossings $U_i$ and $U_j$ for $i\neq j$ is very large relative to the diameters of both $U_i$ and $U_j$. More precisely, the distance is assumed larger than $\frac{1}{\varepsilon}\max\{\mathrm{diam}\ U_i,\mathrm{diam}\ U_j\}$. 
\item[(Emb-4)] For any $t, t'\in\mathfrak{F}_k^r(\gamma)$, $\psi^t$ and $\psi^{t'}$ is chosen so that they coincide outside an $(n+1)$-ball in $U_i$ with radius $\varepsilon_i<\varepsilon$. Indeed, such family of embeddings $\psi^t$ may be obtained by contracting around the center $\{\frac{1}{2}\}\times S^{n-1}\subset I\times S^{n-1}$ of $A_i$ (for all $t\in\mathfrak{F}_k^r(\gamma)$ simultaneously) into a very thin cylinder with small $S^{n-1}$ component and let them approach near $D_i$.
\end{description}
Moreover we can make the following additional assumption after a simultaneous isotopy over $\mathfrak{F}_k^r(\gamma)$:
\begin{description}
\item[(Emb-4')] For any $j$, $\psi^t$ restricted to $\sfD_j$ is constant in $t\in\mathfrak{F}_k^r(\gamma)$, namely, $\psi^t|\sfD_j=\psi^{t'}|\sfD_j$ for any $t, t'\in \mathfrak{F}_k^r(\gamma)$.
\end{description}
Theorem~\ref{thm:BCR-nontrivial} is a consequence of the following identity:
\begin{equation}\label{eq:eval=w}
 \langle (\psi^\bullet)^*\zeta_{k(m-(n+2))},[\mathfrak{F}_k^r(\gamma)]\rangle=w_k(\Gamma(\gamma)). 
\end{equation}

First we see that the contributing subdomain for the integral may be reduced as in the proof of Theorem~\ref{thm:Zuniversal}. Let $\pi^t:C_{\Gamma}(\psi^t)\to C_q(\R^n)$ be the natural projection. Then the form $I(\Gamma)(\psi^t)$ obtained by the integral along the fiber $C_{\Gamma}(\psi^t)$ may be rewritten as 
\[ I(\Gamma)(\psi^t)=\int_{C_{\Gamma}(\psi^t)}\omega(\Gamma)(\psi^t)=\int_{C_q(\R^n)}\pi^t_*\omega(\Gamma)(\psi^t) \]
for some $(qn+k(m-(n+2)))$-form on $C_q(\R^n)$-bundle over $\sigma^{k(m-(n+2))}\subset \mathrm{Emb}(\R^n,\R^m)$. Since $C_q(\R^n)$ is common over the cycle $\sigma^{k(m-(n+2))}$, the integral along $\mathfrak{F}_k^r(\gamma)$ equals
\[ \int_{t\in \mathfrak{F}_k^r(\gamma)}\int_{C_q(\R^n)}\pi^t_*\omega(\Gamma)(\psi^t)
	=\int_{C_q(\R^n)}\int_{t\in \mathfrak{F}_k^r(\gamma)}\pi^t_*\omega(\Gamma)(\psi^t). \]

Now we can show the analogue of (\ref{eq:I(V2)}). Let $V_1(D_i)$ be the subset of $C_q(\R^n)$ consisting of configurations such that no points are mapped into $D_i\subset U_i$ and $V_1(A_i)$ is similarly defined for $A_i\subset U_i$. Then we have
\[ \int_{V_1(D_i)\cup V_1(A_i)}\int_{\mathfrak{F}_k^r(\gamma)}\pi^t_*\omega(\Gamma)(\psi^t)=O(\varepsilon_i). \]
This is because at the limit $\varepsilon_i\to 0$ the integrand form is still well defined and then the $S^{l-3}$-variation of embeddings inside $U_i$ degenerates to the one limiting embedding. So the integral vanishes by a dimensional reason. 

We can also prove that the integrals along $V_2(e)$ and $V_3(e)$ (defined in Lemma~\ref{lem:long-theta-edge}, \ref{lem:long-eta-edge}) are $O(\varepsilon)$ by exactly the same arguments as Lemma~\ref{lem:long-theta-edge}, \ref{lem:long-eta-edge} using (Emb-3) and (Emb-2) respectively. So the contributing domain in $C_q(\R^n)$ reduces to a disjoint union $\widehat{\sfM}_k$ of the spaces of the form
\[ \sfM_k\eqdef \sfA_1\times \sfD_1\times \sfA_2\times \sfD_2\times \cdots \times \sfA_k\times \sfD_k. \]
This time the arrangement of $\sfD_j$'s and $\sfA_j$'s may be different from that for $\Gamma_k$ in Figure~\ref{fig:wheel}(b). But the following conditions are satisfied.
\begin{enumerate}
\item $\sfD_j$ is the only $\sfD_*$ included in $\sfS_j$ .
\item $\sfA_i$'s are radially arranged around $\sfD_j$ in each $\sfS_j$.
\end{enumerate}

We show under a suitable assumptions that there are at most two terms in the sum (\ref{eq:zeta}) for which the integral $I(\Gamma)$ is non degenerate. Let $a(j,p)\in\{1,\ldots,k\}$ be the number such that $\sfA_{a(j,p)}\subset\sfS_j$ is the $p$-th from the center of $\sfS_j$. Then we make the following assumption in addition to (Emb-0)--(Emb-4'):\par\vspace{2mm}
\begin{quote}{\bf (Emb-5'):}
\[\begin{split}
	 \sfA_{a(j,p)}=\Bigl\{&(x_1,\ldots,x_n)\in\R^n\,\Bigl|\Bigr.\,\\
	 	&\Bigl(\frac{\varepsilon^{k+1-p}}{2}\Bigr)^2\leq \|x_1-j\|^2+\|x_2\|^2+\cdots+\|x_n\|^2\leq \Bigl(\frac{2\varepsilon^{k+1-p}}{3}\Bigr)^2\Bigr\},\\
	 \sfD_j=\Bigl\{&(x_1,\ldots,x_n)\in\R^n\,\Bigl|\Bigr.\,\|x_1-j\|^2+\|x_2\|^2+\cdots+\|x_n\|^2\leq (\varepsilon^{k+1})^2\Bigr\}.
	 \end{split}
 \]
\end{quote}

\begin{Lem}\label{lem:non-increasing}
Let $q_r, q_{r+1}, q_{r+2}$ be three successive external vertices of a chord diagram $\Gamma$ lying on a path of $\eta$-edges. Let $V_4(q_r,q_{r+1},q_{r+2};j,p)\subset \sfM_k$ be the subset consisting of configurations such that $q_r,q_{r+1},q_{r+2}$ are mapped into $\sfA_{a(j,p)},\sfA_{a(j,p+2)},\sfA_{a(j,p+1)}$ respectively. Then 
\begin{equation}\label{eq:rrr}
 \int_{V_4(q_r,q_{r+1},q_{r+2};j,p)}\int_{\mathfrak{F}_k^r(\gamma)}\pi^t_*\omega(\Gamma)(\psi^t)=O(\varepsilon). 
\end{equation}
\end{Lem}
\begin{proof}
Let $x_r,x_{r+1},x_{r+2}\in \R^n$ be the coordinates for $q_r,q_{r+1},q_{r+2}$ respectively in the configuration. Then $\omega(\Gamma)$ has factors $u'(x_{r+1}-x_r)^*\omega_{n-1}$ and $u'(x_{r+2}-x_{r+1})^*\omega_{n-1}$. Since $\|x_{r+2}-x_r\|$ is very small relative to $\|x_{r+1}-x_r\|$ and $\|x_{r+2}-x_{r+1}\|$ on $V_4(q_r,q_{r+1},q_{r+2};j,p)$ (see Figure~\ref{fig:xrxrxr}), the integral of the form $\omega'(\Gamma)$, obtained from the form $\omega(\Gamma)$ by replacing $u'(x_{r+2}-x_{r+1})^*\omega_{n-1}$ with $u'(x_r-x_{r+1})^*\omega_{n-1}$, differs from the original one by $O(\varepsilon)$. Moreover, the integral of $\omega'(\Gamma)$ vanishes by a dimensional reason. Thus (\ref{eq:rrr}) is proved.
\begin{figure}
\fig{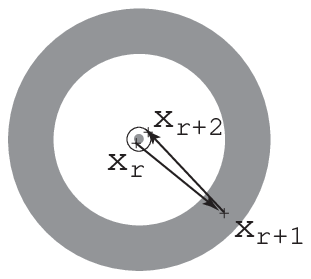}
\caption{}\label{fig:xrxrxr}
\end{figure}
\end{proof}

We can prove the following lemma similarly as Lemma~\ref{lem:non-increasing}.
\begin{Lem}
Let $q_r, q_{r+1}, q_{r+2}$ be three successive external vertices of a chord diagram $\Gamma$ lying on a path of $\eta$-edges. For $V_5(q_r,q_{r+1},q_{r+2};j,p)\subset \widehat{\sfM}_k$ the subset consisting of configurations such that $q_r$, $q_{r+1}$, $q_{r+2}$ are mapped into $\sfA_{a(j,p+1)},\sfA_{a(j,p+2)},\sfA_{a(j,p)}$ respectively, then
\begin{equation}\label{eq:rrr2}
 \int_{V_5(q_r,q_{r+1},q_{r+2};j,p)}\int_{\mathfrak{F}_k^r(\gamma)}\pi^t_*\omega(\Gamma)(\psi^t)=O(\varepsilon). 
\end{equation}
\end{Lem}

By (\ref{eq:rrr}), (\ref{eq:rrr2}) and the observations above, the domains in $\widehat{\sfM}_k$ giving non-degenerate contribution are those consisting of configurations such that the successive external vertices $q_r,q_{r+1},\ldots,q_{r+s}$ lying on a path of $\eta$-edges are mapped bijectively into components
$\sfD_j$, $\sfA_{a(j,1)}$, $\sfA_{a(j,2)}$, $\ldots$, $\sfA_{a(j,s)}$ (or $\sfA_{a(j,s)}$, $\ldots$, $\sfA_{a(j,1)}$, $\sfD_{j}$)
respectively for some $j$. Thus there are at most two chord diagrams for which $\int_{\mathfrak{F}_k^r(\gamma)}I(\Gamma)(\psi^t)$ is non degenerate and if there are two, they are related by a change of all edge orientations involved in the cycle of $\Gamma$. Further, one of them is equal to $\Gamma(\gamma)$.

We can reduce more by the symmetry of graphs. We may assume that $\gamma$ is such that $\Gamma(\gamma)$ does not have a subgraph as in (\ref{eq:y-l}) because otherwise $w_k(\Gamma(\gamma))=0$. Let $\Gamma=\Gamma(\gamma)$ and let $\Gamma^*$ denote $\Gamma$ with the orientations of all edges involved in its cycle reversed. 

If $k$ is odd and $\Gamma=\Gamma^*$, then by Proposition~\ref{prop:odd-chord-diag} extended for higher codimensions, $I(\Gamma)w_k(\Gamma)=0$.

If $k$ is odd and $\Gamma\neq\Gamma^*$, then by Proposition~\ref{prop:odd-chord-diag} extended for higher codimensions, $I(\Gamma)w_k(\Gamma)$ and $I(\Gamma^*)w_k(\Gamma^*)$ cancel each other. Now the result for the case of odd $k$ is proved.

If $k$ is even and $\Gamma=\Gamma^*$, the only contributing term is 
\begin{equation}\label{eq:contrib-term1}
 \frac{1}{2}\int_{\mathfrak{F}_k^r(\gamma)}\frac{I(\Gamma)(\psi^t)w_k(\Gamma)}{|\Aut \Gamma|}
 =\int_{\mathfrak{F}_k^r(\gamma)}\frac{I(\Gamma)(\psi^t)w_k(\Gamma)}{|\pAut \Gamma|}
\end{equation}
where $\pAut \Gamma$ denotes the group of automorphisms of $\Gamma$ considered as an unoriented graph.

If $k$ is even and $\Gamma\neq\Gamma^*$, sum of the contributing terms is equal to 
\begin{equation}\label{eq:contrib-term2}
\frac{1}{2}\int_{\mathfrak{F}_k^r(\gamma)}\frac{I(\Gamma)(\psi^t)w_k(\Gamma)}{|\Aut \Gamma|}
+\frac{1}{2}\int_{\mathfrak{F}_k^r(\gamma)}\frac{I(\Gamma^*)(\psi^t)w_k(\Gamma^*)}{|\Aut \Gamma^*|}
 =\int_{\mathfrak{F}_k^r(\gamma)}\frac{I(\Gamma)(\psi^t)w_k(\Gamma)}{|\pAut \Gamma|}. 
\end{equation}
Hence in either case for even $k$, the computations can be done equally. Indeed, Lemma~\ref{lem:IIG=AutG} below implies that the above two equal $w_k(\Gamma)+O(\varepsilon)$. Therefore (\ref{eq:eval=w}) is proved and it completes the proof of Theorem~\ref{thm:BCR-nontrivial}.
\end{proof}

\begin{Lem}\label{lem:IIG=AutG}Let $k$ be even and $\geq 2$. Then under the assumptions (Emb-0)--(Emb-5'), we have 
\[ \int_{\mathfrak{F}_k^r(\gamma)}I(\Gamma)(\psi^t)=|\pAut\Gamma|+O(\varepsilon). \]
\end{Lem}
\begin{proof}
There are exactly $|\pAut \Gamma|$ connected components in $\widehat{\sfM}_k$ over each of which the integral is non degenerate. Further, we know that if $\Gamma$ does not have a subgraph as in (\ref{eq:y-l}), any element in $\pAut\Gamma$ is a combination of a rotation along the cycle and a reversion of orientations of all the edges involved in the cycle. A rotation does not change the integral and the reversion does not change the integral too as observed above. Hence all $|\pAut\Gamma|$ components contribute as a common value. So we need to compute for one connected component in which external vertices $q_r, q_{r+1},\ldots, q_{r+s}$ on each path of $\eta$-edges with edge orientations coincident with this order are mapped into some $\sfA_{a(j,s)},\ldots, \sfA_{a(j,1)}, \sfD_j$ respectively.

It suffices to prove that
\[ \int_{\sfM_k}\int_{\mathfrak{F}_k^r(\gamma)}
		\bigwedge_i\Bigl(
			u(\psi^t a_i-\psi^t d_i)^*\omega_{m-1}\wedge
			u'(\alpha_i-a_i)^*\omega_{n-1}
		\Bigr)=1+O(\varepsilon)\]
where $\alpha_i$ is some $a_*$ or $d_*$, which is the image of the target vertex of an $\eta$-edge. By (Emb-5'), we can rewrite the LHS as
\[ \int_{\sfM_k}\int_{\mathfrak{F}_k^r(\gamma)}
		\bigwedge_i\Bigl(
			u(\psi^t a_i-\psi^t d_i)^*\omega_{m-1}\wedge
			u'(d_{\tau(i)}^0-a_i)^*\omega_{n-1}
		\Bigr)+O(\varepsilon)
\]
where $\tau(i)=j$ if $\sfA_i\subset\sfS_j$ and $d_p^0\eqdef (p,0,\ldots,0)\in\R^n$. We will write $\Phi_i(\psi)\eqdef u(\psi a_i-\psi d_i)^*\omega_{m-1}\wedge u'(d_{\tau(i)}^0-a_i)^*\omega_{n-1}$ for simplicity. Then the above expression equals
\[ \int_{t\in\prod_{i\neq j}\mathfrak{F}_k^r(\gamma;i)}\int_{\sfM_k(j)}
	\bigwedge_{i\neq j}\Phi_i(\psi^{(t,*)})
	\int_{t_j\in\mathfrak{F}_k^r(\gamma;j)}\int_{\sfA_j\times\sfD_j}\Phi_j(\psi^{(t,t_j)})+O(\varepsilon) \]
where $\mathfrak{F}_k^r(\gamma;i)$ corresponds the fiber of the pullback $S^{l-3}$-bundle from the $i$-th component of $Q_l^{\times k}$ and $*\in\mathfrak{F}_k^r(\gamma;i)$ is the base point. Then the result follows by an iterative use of the following identity.
\begin{equation}\label{eq:int-a-d=1}
 \int_{t_j\in\mathfrak{F}_k^r(\gamma;j)}\int_{\sfA_j\times\sfD_j}\Phi_j(\psi^{(t,t_j)})=1.
\end{equation}

To prove (\ref{eq:int-a-d=1}), we choose some embeddings $\R^m\to \R^n$ independently of $t_j$, which coincide outside $U_j$ to $\psi^{(t,t_j)}$ and whose restrictions to $\sfA_j\cup\sfD_j$ are disjoint from $\psi^{(t,t_j)}(\sfA_j\cup\sfD_j)$ inside $U_i$, and which enclose $\psi^{(t,t_j)}(\sfA_j\cup\sfD_j)$ as in the proof of Lemma~\ref{lem:part-I(GW)}. Since such closing embeddings are chosen constant in $t_j\in\mathfrak{F}_k^r(\gamma;j)$, the LHS integrals of (\ref{eq:int-a-d=1}) with $\psi^{(t,t_j)}$ replaced by those closing embeddings vanish by a dimensional reason. Therefore the LHS of (\ref{eq:int-a-d=1}) equals
\[\begin{split}
	&\int_{\mathfrak{F}_k^r(\gamma;j)}\int_{\sfA_j^{\circ}\times\sfD_j^{\circ}}
	\Phi_j(\psi_j^{(t,t_j)\circ\circ})\\
	=&\int_{\mathfrak{F}_k^r(\gamma;j)}\int_{(v,l,d_j)\in S^{n-1}\times S^1\times \sfD_j^{\circ}}
	u(\psi^{(t,t_j)\circ\circ}_jl-\psi^{(t,t_j)\circ\circ}_jd_j)^*\omega_{m-1}
		\wedge u'(d_{\tau(j)}^0-(v,l))^*\omega_{n-1}\\
	=&\int_{\mathfrak{F}_k^r(\gamma;j)}\int_{\sfL_j^{\circ}(v)\times \sfD_j^{\circ}}
	u(\psi^{(t,t_j)\circ\circ}_jl-\psi^{(t,t_j)\circ\circ}_jd_j)^*\omega_{m-1}
		\int_{v\in S^{n-1}}\omega_{n-1}(v)\\
	=&\int_{\mathfrak{F}_k^r(\gamma;j)}\int_{\sfL_j^{\circ}(v)\times \sfD_j^{\circ}}
	u(\psi^{(t,t_j)\circ\circ}_jl-\psi^{(t,*)\circ\circ}_jd_j)^*\omega_{m-1}
		\int_{v\in S^{n-1}}\omega_{n-1}(v)
	\end{split}
\]
where $\psi_j^{(t,t_j)\circ\circ}$ is the embedding of $\sfA_j^\circ\cup\sfD_j^\circ$ into $U_j\cup_{\partial}(-U_j)$ obtained by closing $\psi^{(t,t_j)}$ with the above chosen embeddings and $\sfL_j$ is defined in (\ref{eq:def-L}). The last equality follows from (Emb-4'). By a similar argument as in the final part of the proof of Lemma~\ref{lem:part-I(GW)}, one may see that the linking number of the suspended arcs $\psi^{(t,t_j)}(\sfL_j(v))$ over $\mathfrak{F}_k^r(\gamma;j)$, which forms $S^{l-2}$, and the $n$-sphere $\psi_j^{(t,*)\circ\circ}(\sfD_j)$ is 1. Hence (\ref{eq:int-a-d=1}) is proved and it completes the proof of Lemma~\ref{lem:IIG=AutG}.
\end{proof}

\section{BCR invariant for long $n$-knots other than ribbon}\label{s:remark}

Now we shall define a certain kind of long $n$-knot which we will call a long $(p,q)$-handle knot. Let $L_1,L_2,\ldots,L_m$ be disjoint copies of the higher dimensional trivially framed Hopf link placed on the trivial long $n$-knot $:S^{p-1}\sqcup S^{q-1}\hookrightarrow\R^{n}\subset\R^{n+2}$ such that $p+q=n+1$. An example of this step for $(p,q)=(1,2)$ is given as follows. We choose a disjoint union of $(0,1)$-dimensional Hopf links $S^0\sqcup S^1\hookrightarrow\R^2$ trivially embedded into $\R^2$:
\cfig{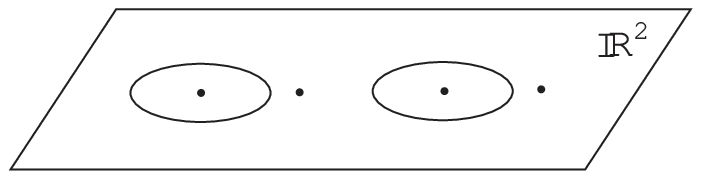}

Consider $\R^n$ is bounded by the lower half space $\R_-^{n+1}\subset\R^{n+2}$. Then attach $(n+1)$-dimensional $p$-handles along $S^{p-1}$'s and $(n+1)$-dimensional $q$-handles along $S^{q-1}$'s to $\R_-^{n+1}$, so that each $p$-handle $D^{q+2}\times D^p$ (resp. $q$-handle $D^{p+2}\times D^q$) is obtained by slight thickening (say, $\varepsilon$-tubular neighborhood for sufficiently small $\varepsilon>0$) the upper hemisphere of $S^p$ (resp. $S^q$). Then the boundary of the resulting handlebody $H_{p,q}$ is again a trivial knot. For $(p,q)=(1,2)$, this step is seen as follows. We attach 2-handles in the place of $S^1$-components and 1-handles in the place of $S^0$-components:
\cfig{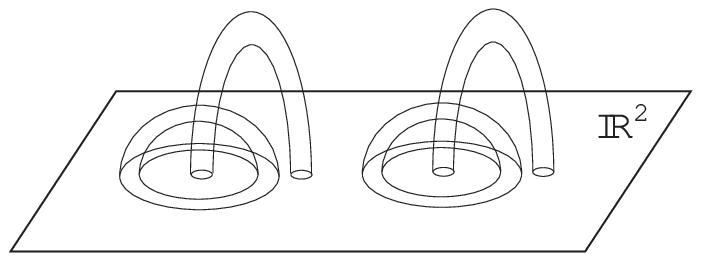}
The resulting handlebody is $H_{1,2}$. Here we assume that for some choice of $(n_1,n_2)$ of an orthonormal 2-framing normal to $\R^2$, all 2-handles are included in the 3-dimensional plane spanned by $\R^2$ and $n_2$, and all 1-handles are included in the 3-dimensional plane spanned by $\R^2$ and $n_1$ so that $H_{1,2}$ does not have a self intersection.

We consider a crossing change between a $p$-handle and a $q$-handle of $H_{p,q}$, which is an analogue of an unclasping of a crossing of a long ribbon $n$-knot, to obtain a non-trivial knot. It is defined as follows: Let $N_a^p$ and $N_b^q$ be $2\varepsilon$-neighborhoods in $\R^{n+2}$ of points on some handles: $(0,a)\in D^{q+2}\times D^p$ and $(0,b)\in D^{p+2}\times D^q$ respectively. Then attach an $(n+2)$-dimensional 1-handle $M$ connecting $N_a^p$ and $N_b^q$ so that it is disjoint from the handlebody and denote $N_a^p\cup M\cup N_b^q$ by $N_{ab}$. Then we can standardly embed in $N_{ab}$ a Hopf link $L':S^{q}\cup S^p$ so that they are disjoint from $H_{p,q}$ and $\mathrm{lk}(S^p,\mbox{core of the $q$-handle})=\mathrm{lk}(S^q,\mbox{core of the $p$-handle})=\pm 1$. We call the modification of the handlebody via a surgery along $L'$ in $N_{ab}$ a {\it crossing change along $N_{ab}$}. Note that this modifies the ambient space $\R^{n+2}$ into another $\R^{n+2}$. Thus the result of the modification of the pair $(\R^{n+2},\partial{H_{p,q}})$ is considered to be another long knot in $\R^{n+2}$. For $(p,q)=(1,2)$, this step is seen as follows. We embed disjointly the Hopf links $S^1\cup S^2$:
\cfig{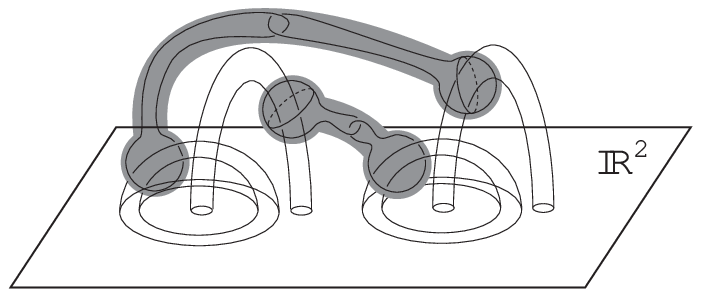}

A {\it long $(p,q)$-handle knot} is defined to be a long $n$-knot obtained from $H_{p,q}$ by a sequence of disjoint crossing changes followed by taking its boundary. Note that if $p,q>1$, isotopy type of a $(p,q)$-handle knot is determined by some equivalence class of the matrix of the bilinear form
\[ L:H_p(V_{p,q};\Q)\otimes H_q(\R^{n+2}-V_{p,q};\Q)\to\Q \]
defined by linking numbers between certain basis elements, where $V_{p,q}$ is a handlebody obtained from $H_{p,q}$ by a sequence of crossing changes.\footnote{That the isotopy types of $(p,q)$-handle knots are determined by $L$ was suggested to the author by K.~Habiro. He also informed me about Levine's generalization of the Alexander invariant in \cite{Lev}.} 

\begin{Exa}
Long $(1,n)$-handle knots are long ribbon $n$-knots and if $n=2q-1$, long $(q,q)$-handle knots are the long versions of simple knots in \cite{Lev}.
\end{Exa}

In \cite{Lev}, a generalization of the Alexander invariant is defined by using the bilinear form like $L$ for general $n$-knots. Denote the invariant with a suitable normalization (by the conditions: $\Delta_p(1)=1, \Delta_p'(1)=0$) by ${\Delta}_p(t)$ and expand it as follows.

\begin{equation*}
\log{\Delta}_p(t)|_{t=e^h}=\alpha^p_2h^2+\alpha^p_3h^3+\ldots\in\\\Q[[h]].
\end{equation*}

Then the following theorem can be proved similarly as Theorem~\ref{thm:Zuniversal}. 
\begin{Thm}
Let $n$ be an odd integer $>1$. For long handle $(p,q)$-knots with $p+q=n-1$, $1<p<q$, we have
\begin{equation}
	\begin{split}
		\hat{z}_{2k}&\in\alpha^p_{2k}+\R[\alpha^p_{2k-1},\ldots,\alpha^p_3,\alpha^p_2]^{(\deg{\leq 2k-1})}\\
		\hat{z}_{2k+1}&\in\R[\alpha^p_{2k},\ldots,\alpha^p_3,\alpha^p_2]^{(\deg{\leq 2k})}
	\end{split}
\end{equation}
for $k\geq 1$.
\end{Thm}
In this case, finite type invariant is also defined in the space of (a certain equivalence classes of) matrices associated to the bilinear form $L$ by some modifications on the matrices representing the crossing changes. For 1-knots, such  a finite type invariant is considered in detail in \cite{MO} (and it may also be generalized to simple knots straightforwardly). That the finite type invariant of matrices reduces to the polynomials in $\alpha_j^p$'s is shown by counting the dimensions of both spaces.
\begin{Conj}
Let $n$ be an odd integer $>1$. For arbitrary long $n$-knots, we have
\begin{equation}
	\begin{split}
		\hat{z}_{2k}&\in\sum_p\alpha^p_{2k}+\R[\alpha^p_{2k-1},\ldots,\alpha^p_3,\alpha^p_2]^{(\deg{\leq 2k-1})}_p\\
		\hat{z}_{2k+1}&\in\R[\alpha^p_{2k},\ldots,\alpha^p_3,\alpha^p_2]^{(\deg{\leq 2k})}_p
	\end{split}
\end{equation}
for $k\geq 1$, where $p$ runs over $1\leq p\leq \lfloor\frac{n}{2}\rfloor$.
\end{Conj}

The Alexander polynomial is known to be a (finite) polynomial with integral coefficients. So it is natural to expect a similar properties for $\hat{z}_k$.
\begin{Prob}
Does $Z=\exp{\Bigl(\sum_{k\geq 2}x_k\hat{z}_k\Bigr)}$ have an integrality and a finiteness property for a suitable choice of $\{x_k\}$?
\end{Prob}
This problem is related to the following realization problem.
\begin{Prob}
Determine what kind of series is realizable as a value of $Z$ of a long $n$-knot. Moreover, determine the complete set of images of $Z$.
\end{Prob}

Other generalization that can be considered is configuration space integral invariant for links $\R^{n_1}\sqcup\cdots\sqcup\R^{n_r}\hookrightarrow\R^{N}$. This is related to the Massey product and the coefficients of Farber's non-commutative invariant \cite{Far} and so on. We will explain this in \cite{W2}.
\begin{appendix}
\section{Invariance proof}
In \cite{CR} and \cite{R}, the invariance of the higher degree $z_k$ is claimed and the idea of the proof is given. But the explicit descriptions of $z_k$ for higher degrees are not given there while the complete definitions and proofs are given for degrees up to 3. So in this section, we give a proof of Theorem~\ref{thm:invariant} filling the details in their idea with our definition of $z_k$, conventions on the diagrams and the orientations on configuration spaces. We also see that $z_k$ can be obtained naturally from a general framework of diagrams and of the geometry of configuration spaces.

\subsection{Descriptions of faces in $\partial C_{\Gamma}$}\label{ss:face}

\subsubsection{Codimension one faces}
It is known that $C_{\Gamma}$ is a smooth manifold with corners and $\partial C_{\Gamma}$ admits a stratification \cite{FM, BT}. We will need only to consider the codimension one faces in $\partial C_{\Gamma}$ for our purpose. So we review here how each codimension one face can be described. 

As observed in \cite{FM,AS}, the set of codimension one faces in $\partial C_{\Gamma}$ is in bijective correspondence with the set of all subsets of the vertices of $\Gamma$, which are obtained by blowing-up along the corresponding diagonals. Denote by $\calS_{\Gamma,A}$ and $\calS_{\Gamma,\infty,i}$ the codimension one faces corresponding to $\Delta_A$ and $\Delta_{\infty,i}$ respectively. 

Let us see how $\calS_{\Gamma,A}$ is presented briefly. Details are found in \cite{R}. Here we shall disregard the dummy factor for simplicity. Let $V(\Gamma)$ denote the set of vertices on $\Gamma$. There are two cases:

\begin{enumerate}
\item[\bf Case 1] If $A\subset V(\Gamma)$ does not have external vertices, then $\calS_{\Gamma,A}$ is the pullback bundle $\varphi^*B_A=B_A\times C_{\Gamma/A}$ in the following commutative diagram.
\begin{equation}\label{cd:faceonspace}
\begin{CD}
\calS_{\Gamma,A} @>\hat{\varphi}>> B_A\\
@V\hat{\pi}^{\partial}VV @VV\pi^{\partial}V\\
C_{\Gamma/A} @>{\varphi}>> *
\end{CD}
\end{equation}
Here $\pi^{\partial}$ and $\varphi$ are the maps to a point and $C_{\Gamma/A}$ is the space of configurations obtained from configurations in $C_{\Gamma}$ by collapsing points in $A$. We consider each piece in (\ref{cd:faceonspace}) as a fiber over a point in $\mathrm{Emb}(\R^n,\R^{n+2})$.

\item[\bf Case 2]If $A\subset V(\Gamma)$ has external vertices, then $\calS_{\Gamma,A}$ is the pullback bundle $\varphi^*\widehat{B}_A$ in the following commutative diagram.
\begin{equation}\label{cd:faceonknot}
  \begin{CD}
	\calS_{\Gamma,A}  @>\hat{\varphi}>>  \widehat{B}_A\\
	@V\hat{\pi}^{\partial}VV @VV\pi^{\partial}V\\
	C_{\Gamma/A} @>\varphi>>  I_n(\R^{n+2})
	\end{CD}
\end{equation}
Here 
\begin{itemize}
\item $\varphi$ is the composition of the embedding and the generalized Gauss map (the tangent map together with an assignment of the $n$-frame in $\R^{n+2}$ determined by the embedding) at the point where the points in $A$ coincide, 
\item $\widehat{B}_A$ is the space of configurations of points on $(\R^{n+2},\iota\R^n)$, ($\iota\in I_n(\R^{n+2})$ is an $n$-frame in $\R^{n+2}$) modulo translations and dilations, together with $\iota\in I_n(\R^{n+2})$, 
\item $\pi^{\partial}$ is the map giving the underlying $n$-frame in $\R^{n+2}$.
\end{itemize}
\end{enumerate}

For the face $\calS_{\Gamma,\infty,i}$, the fiber of the unit normal bundle $SN\Delta_{\infty,i}$ over a point in $\Delta_{\infty,i}\subset (S^{n+2})^{k}$ where the point $x_i$ on knot (resp. in space) coincide at $\infty\in S^{n}\subset S^{n+2}$, is identified with the set of points in $T_{\infty} S^{n}$ (resp. $T_{\infty}S^{n+2}$) modulo overall translations and dilations along $T_{\infty} S^{n}$. Namely, the fiber is identified with the $(n-1)$-dimensional unit sphere (resp. $(n+1)$-dimensional unit sphere). Thus $\calS_{\Gamma,\infty,i}$ may be identified with $C_{\Gamma\setminus \{v_i\}}\times S^{n-1}$ (resp. $C_{\Gamma\setminus \{v_i\}}\times S^{n+1}$). We will call a point in $B_A$ (or $\widehat{B}_A$) a {\it relative configuration}.

The forms on $C_{\Gamma}^0$ defined in \S\ref{ss:def-integral} extend smoothly to $C_{\Gamma}$ as the naturally defined forms by Gauss maps via the coordinates determined by the infinitesimal embedding $\iota$.

\subsubsection{Classification of codimension one faces}

Theorem~\ref{thm:invariant} will be proved by looking at the integrals restricted to the codimension one faces classified as follows.
\begin{description}
\item[Principal faces]The faces corresponding to the diagonal where exactly two of the $|V(\Gamma)|$ points coincide in $S^{n+2}\setminus\{\infty\}$.
\item[Hidden faces]The faces corresponding to the diagonal where at least 3 of $|V(\Gamma)|$ points but not all points coincide in $S^{n+2}\setminus\{\infty\}$.
\item[Infinite faces]The faces corresponding to the diagonal where at least 3 of $|V(\Gamma)|$ points coincide in $\infty\in S^{n}$.
\item[Anomalous faces]The faces corresponding to the diagonal where all the $|V(\Gamma)|$ points coincide in $S^{n}\setminus\{\infty\}$.
\end{description}

Let $\pi:M\to B$ be a bundle with $n$ dimensional fiber $F$. Then the {\it push-forward} (or {\it integral along the fiber}) $\pi_*\omega$ of an $(n+p)$-form $\omega$ on $M$ is a $p$-form on $B$ defined by
\[ \int_c \pi_*\omega=\int_{\pi^{-1}(c)}\omega, \]
where $c$ is a $p$-dimensional chain in $B$. 

Let $\pi^{\partial}:\partial_F M\to B$ be the restriction of $\pi$ to $\partial F$-bundle with the orientation induced from $\mathrm{Int}{(F)}$, i.e., $\Omega(\partial F)=i_n\Omega(F)$ where $n$ is the in-going normal vector field over $\partial F$. Then the generalized Stokes theorem for the pushforward is 
\begin{equation}\label{eq:stokes}
 d\pi_*\omega=\pi_*d\omega+(-1)^{\deg{\pi_*^{\partial}\omega}}\pi^{\partial}_*\omega.
\end{equation}The derivation of (\ref{eq:stokes}) is e.g., in \cite[Appendix]{BT}.
\begin{proof}[Proof of Theorem~\ref{thm:invariant}]
The formula (\ref{eq:stokes}) is used to prove the invariance of $z_k$ (or $\hat{z}_3$) as follows. Consider $C_{\Gamma}(\cdot)$ as a bundle over the space of embeddings $\mathrm{Emb}(\R^n,\R^{n+2})$ with fiber the configuration space $C_{\Gamma}$ and consider $z_k$ as a $0$-form on $\mathrm{Emb}(\R^n,\R^{n+2})$. Then the invariance of $z_k$ relies on the closedness of it because if two embeddings $\psi_0$ and $\psi_1$ are connected by a smoothly parametrized embeddings $\psi_t\ (t\in [0,1])$, then by Stokes' theorem,
\[ z_k(\psi_1)-z_k(\psi_0)=\int_0^1 d z_k(\psi_t). \]

Since $\omega(\Gamma)$ is closed, 
\[ dI(\Gamma)=\int_{\partial C_{\Gamma}}\omega(\Gamma) \]
by (\ref{eq:stokes}) where the integral of the RHS is the pushforward restricted to $\partial C_{\Gamma}$. In particular, only the codimension one faces in $\partial C_{\Gamma}$ contribute to the above integral. Therefore the obstruction to the closedness of $z_k$ is
\begin{equation}\label{eq:int-codim-one}
 dz_k=\sum_{\Gamma}\frac{w_k(\Gamma)}{|\Aut \Gamma|}\sum_{S\subset \partial C_{\Gamma}\mbox{\tiny: codim 1}}\int_S\omega(\Gamma). 
\end{equation}

By Proposition~\ref{prop:principal} below, all the principal face contributions cancel each other in the sum.

By Proposition~\ref{prop:hiddenall}, \ref{prop:infvanish} below, the contributions of the hidden and infinite faces also vanish by kinds of involutive symmetries on faces or otherwise by dimensional reasons.

If $n$ is odd, Proposition~\ref{prop:anomaly-n-odd}, which is proved again by some involutive symmetries, shows that the anomalous faces do not contribute and thus (\ref{eq:int-codim-one}) is proved to be zero. If $n=2$ and $k=3$, Proposition~\ref{prop:anomaly-n=2} shows that by letting $\widetilde{\Theta}(\Gamma)=\int_{C_1}\varphi^*\hat{\rho}(\Gamma)$, the addition of the correction term to $z_3$:
\[ \hat{z}_3(\psi)\eqdef {z}_3(\psi)+\sum_{\Gamma\in\calG_3^0}\frac{\widetilde{\Theta}(\Gamma)(\psi)w_3(\Gamma)}{|\Aut \Gamma|} \]
is proved to be $d$-closed.

\end{proof}

\subsection{Configuration space integral restricted to the faces}\label{ss:integral-faces}
Now we shall see the explicit form of the integral $I(\Gamma)$ extended to the codimension one faces. 

In the diagram (\ref{cd:faceonspace}), the form $\omega(\Gamma)$ extended to $\calS_{\Gamma,A}$ can be written as 
\begin{equation}\label{eq:omega-decomp}
 \omega(\Gamma)=\hat{\varphi}^*\lambda_1(\Gamma_A)\wedge \hat{\pi}^{\partial*}\lambda_2(\Gamma/A)
\end{equation}
for some $\lambda_1(\Gamma_A)\in\Omega^*({B}_{A})$ and $\lambda_2(\Gamma/A)\in\Omega^*(C_{\Gamma/A})$ both determined by graphs and the Gauss maps. Therefore $\hat{\pi}^{\partial}_*\omega(\Gamma)=\hat{\pi}^{\partial}_*\hat{\varphi}^*\lambda_1(\Gamma_A)\wedge\lambda_2(\Gamma/A)=\varphi^*\pi_*^{\partial}\lambda_1(\Gamma_A)\wedge \lambda_2(\Gamma/A)$ by commutativity of (\ref{cd:faceonspace}). 

Similarly, in the diagram (\ref{cd:faceonknot}), the form $\omega(\Gamma)$ extended to $\calS_{\Gamma,A}$ has the decomposition (\ref{eq:omega-decomp}) too for some $\lambda_1(\Gamma_A)\in\Omega^*(\widehat{B}_{A})$ and $\lambda_2(\Gamma/A)\in\Omega^*(C_{\Gamma/A})$ both determined by graphs. $\hat{\pi}^{\partial}_*\omega(\Gamma_A)=\varphi^*\pi_*^{\partial}\lambda_1(\Gamma_A)\wedge \lambda_2(\Gamma/A)$ also holds.

In other words, in both cases the integral along the fiber $C_{\Gamma/A}$:
\[ I(\Gamma,A)\eqdef \int_{C_{\Gamma/A}}\varphi^*\pi_*^{\partial}\lambda_1(\Gamma)\wedge \lambda_2(\Gamma/A) \]
with the orientation induced from $\mathrm{int}(C_{\Gamma})$, is precisely the integral of $\omega(\Gamma)$ along the fiber $\calS_{\Gamma,A}$, as a fiber over the space of embeddings $\mathrm{Emb}(\R^n,\R^{n+2})$.

\subsection{Principal faces}\label{ss:principal}
In this subsection, we will prove the following proposition.
\begin{Prop}\label{prop:principal}For $A\subset V(\Gamma)$ with $|A|=2$, 
\[ \sum_{\Gamma\in\calG_k^0}\sum_A\frac{I(\Gamma,A)w_k(\Gamma)}{|\Aut \Gamma|}=0. \]
\end{Prop}
This proposition shows that the configuration space integrals restricted to principal faces cancel each other.

\subsubsection{Quasi Jacobi diagrams}

Since each principal face corresponds to a collapsing of an edge in a Jacobi diagram, it is represented by using the following graphs obtained by collapsing an edge. We will call such a graph a {\it quasi Jacobi diagram}. We say that the vertex where an edge of $\Gamma$ has been collapsed is {\it exceptional}. Vertex orientation on a quasi Jacobi diagram is also defined. It is defined as a choice of a bijection
\[ o_v:\{\mbox{two ingoing $\theta$-edges incident to $v$}\}\to\{1,2\}  \]
for each non-exceptional internal vertices $v$ and a bijection
\[ \begin{array}{ll}
	o_v:\{\mbox{three ingoing $\theta$-edges incident to $v$}\}\to\{1,2,3\} \quad &\mbox{if $v$ is 4-valent}\\	
	o_v:\{\mbox{two ingoing $\theta$-edges incident to $v$}\}\to\{1,2\} \quad &\mbox{otherwise}
	\end{array} \]
for the exceptional vertex $v$ having incident $\theta$-edges.

Let $\Gamma$ be a Jacobi diagram $\R^n$ and $e$ be a $\theta$-edge of $\Gamma$. The operator $\delta_e$ of $\Gamma$ giving a quasi Jacobi diagram is defined as follows:
\begin{equation}\label{eq:cntr}
\fig{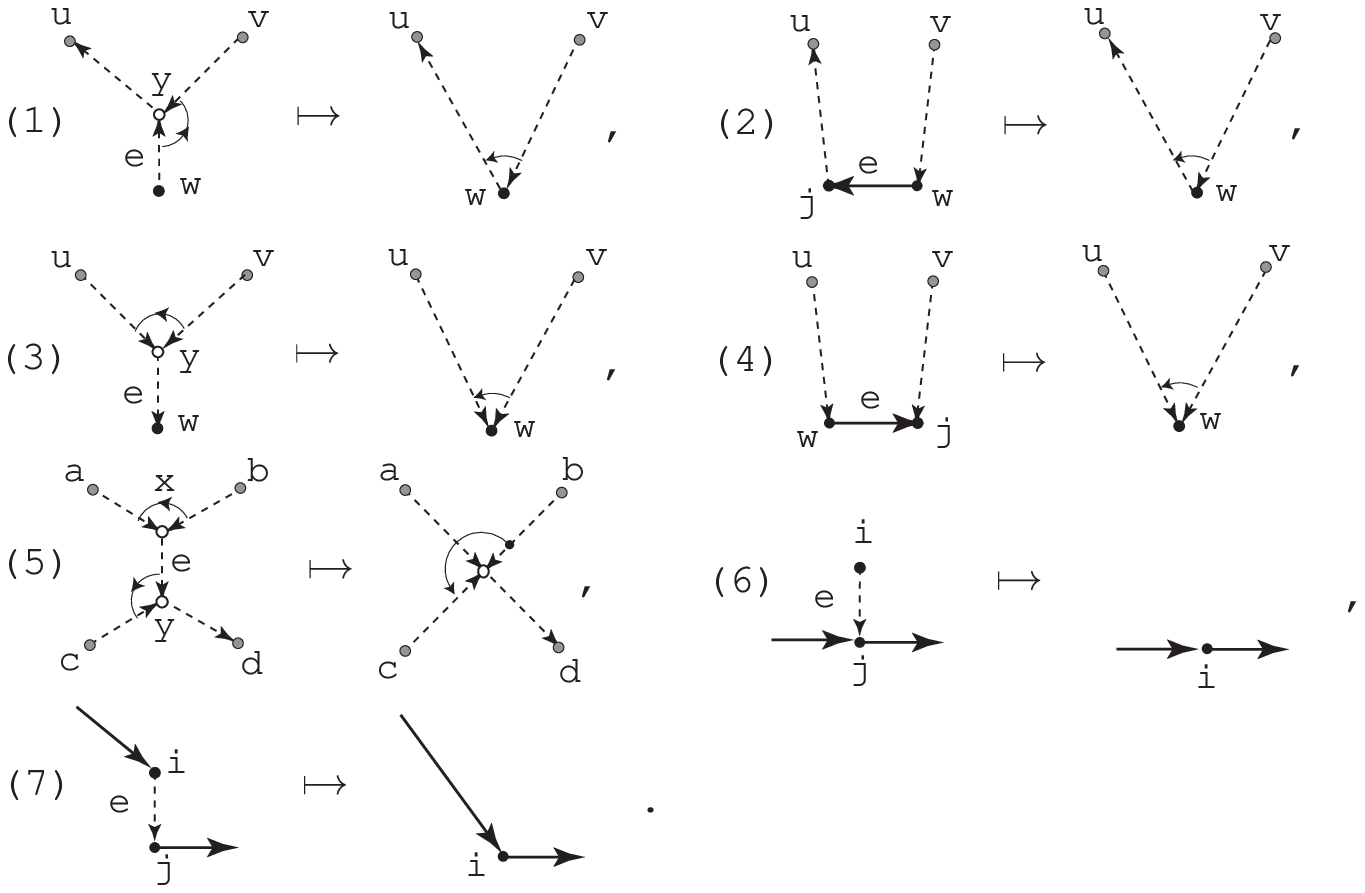}
\end{equation}

\subsubsection{Compactified configuration space for $\delta_e\Gamma$}

We define the configuration space $C^0_{\delta_e\Gamma}(\psi)$ so that a point in $C^0_{\delta_e\Gamma}(\psi)$ represents a position on the diagonal $\Delta_{\{ij\}}$ in $(S^{n+2})^k$ where $e=(i,j)$. Let $C_{\delta_e\Gamma}(\psi)$ be its closure in $C_{\Gamma}(\psi)$.
\begin{itemize}
\item If $e$ is a $\theta$-edge of $\Gamma$ both of whose ends are internal, (the case (5) in (\ref{eq:cntr}))
\[ \begin{split}
	C_{\delta_e\Gamma}^0(\psi)\eqdef \{&(x_1,\ldots,x_q,x_{q+1},\ldots x_{q+s}),\\
		&x_1,\ldots,x_q\in \R^n, x_{q+1},\ldots,x_{q+s-1}\in\R^{n+2}\times S^{n-1}\times S^{n-1},\\
		&x_{q+s}\in\R^{n+2}\times S^{n-1}\times S^{n-1}\times S^{n-1}\,|\,\\
		&p_1(x_i)\neq p_1(x_j)\quad \mbox{if $i,j\geq q+1$, $i\neq j$, $(i,j)$: edge of $\Gamma$},\\
		&\psi(x_i)\neq\psi(x_j)\quad \mbox{if $i,j\leq q$, $i\neq j$, $(i,j)$: edge of $\Gamma$},\\
		&\psi(x_i)\neq p_1(x_j)\quad \mbox{if $i\leq q$, $j\geq q+1$, $(i,j)$: edge of $\Gamma$}\}.
	\end{split} \]
\item If $e$ is a $\theta$-edge on $\Gamma$ one of whose ends is external, or if $e$ is a $\eta$-edge or a chord in $\Gamma$, (the cases (1), (2), (3), (4), (6), (7) in (\ref{eq:cntr}))
\[ \begin{split}
	C_{\delta_e\Gamma}^0(\psi)\eqdef \{&(x_1,\ldots,x_q,x_{q+1},\ldots x_{q+s}),\\
		&x_1,\ldots,x_{q}\in \R^n, x_{q+1},\ldots,x_{q+s}\in\R^{n+2}\times S^{n-1}\times S^{n-1}\,|\,\\
		&p_1(x_i)\neq p_1(x_j)\quad \mbox{if $i,j\geq q+1$, $i\neq j$, $(i,j)$: edge of $\Gamma$},\\
		&\psi(x_i)\neq\psi(x_j)\quad \mbox{if $i,j\leq q$, $i\neq j$, $(i,j)$: edge of $\Gamma$},\\
		&\psi(x_i)\neq p_1(x_j)\quad \mbox{if $i\leq q$, $j\geq q+1$, $(i,j)$: edge of $\Gamma$}\},
	\end{split} \]
\end{itemize}
Then the principal faces corresponding to the contractions in (\ref{eq:cntr}) are identified with the following spaces.
\begin{equation}\label{eq:type-space}
 \begin{array}{c|r}
	\mbox{type} & \mbox{space}\\\hline
	\mbox{(1)} & S^{n+1}\times S^{n-1}\times S^{n-1}\times C_{\delta_e\Gamma} \\
	\mbox{(2)} & S^{n-1}\times C_{\delta_e\Gamma}\\
	\mbox{(3)} & S^{n+1}\times S^{n-1}\times S^{n-1}\times C_{\delta_e\Gamma}\\
	\mbox{(4)} & S^{n-1}\times C_{\delta_e\Gamma}\\
	\mbox{(5)} & S^{n+1}\times S^{n-1}\times C_{\delta_e\Gamma}\\
	\mbox{(6)} &  S^{n-1}\times C_{\delta_e\Gamma}\\
	\mbox{(7)} &  S^{n-1}\times C_{\delta_e\Gamma}
	\end{array}
\end{equation}	

\begin{Lem}\label{lem:ori2}
Let $\Gamma_x$ be a quasi Jacobi diagram with one exceptional vertex $x$ and $\overline\Gamma_x$ be $\Gamma_x$ with its vertex orientation reversed. Then there exists a choice of an orientation $\Omega=\Omega(\Gamma_x)$ on $C_{\Gamma_x}(\psi)$ such that
\begin{equation}\label{eq:omega2}
\Omega(\overline\Gamma_x)=(-1)^n\Omega(\Gamma_x).
\end{equation}
\end{Lem}
\begin{proof}
If the $\theta$ part of $\Gamma_x$ is unitrivalent except for one internal tetravalent vertex $x$, where four edges $e_1, e_2, e_3, e_4$ with $\partial e_i=(v_i, x)$ for $i=1,2,3$, $\partial e_4=(x, v_4)$, $o_x(e_i)=i$ meet, set
\[	\begin{split}
	\Omega_x&\eqdef dX_x^1\wedge \cdots\wedge dX_x^{n+2}\\
		&\quad\wedge\omega_{n-1}(x^{(1)})\wedge\omega_{n-1}(x^{(2)})\wedge\omega_{n-1}(x^{(3)})\\
		&\quad\wedge\Omega_{\bar{e}_{1-}}\wedge\Omega_{\bar{e}_{2-}}\wedge\Omega_{\bar{e}_{3-}}\wedge\Omega_{\bar{e}_{4+}},\\
	\Omega(\Gamma_x)&\eqdef \Omega_x\wedge
		\bigwedge_{\tworows{e\in E_{\mathrm{\theta}}(\Gamma)\setminus}{\{e_1,e_2,e_3,e_4\}}}\Omega_e.
	\end{split} \]

If the $\theta$ part of $\Gamma_x$ is unitrivalent except for one external bivalent vertex $x$, where two edges $e_1, e_2$ with $\partial e_i=(v_i,x)$, $o_x(e_i)=i$ meet, set
\[	\begin{split}
	\Omega_x&\eqdef dX_x^1\wedge\cdots\wedge dX_x^n
		\wedge\Omega_{\bar{e}_{1-}}\wedge\Omega_{\bar{e}_{2-}},\\
	\Omega(\Gamma_x)&\eqdef \Omega_x\wedge\bigwedge_{\tworows{e\in E_{\mathrm{\theta}}(\Gamma)\setminus}{\{e_1,e_2\}}}\Omega_e.
	\end{split} \]

The cases of other edge orientations are similar. 

If $\Gamma_x$ has an exceptional vertex $x$ having 2 or 3 incident $\eta$-edges, we define
\[	\begin{split}
	\Omega_x&\eqdef dX_x^1\wedge\cdots\wedge dX_x^n,\\
	\Omega(\Gamma_x)&\eqdef \Omega_x\wedge\bigwedge_{e\in E_{\mathrm{\theta}}(\Gamma)}\Omega_e.
	\end{split} \]
Then the property (\ref{eq:omega2}) follows from the definition.
\end{proof}

\subsubsection{Canceling of the integral along principal faces}
\begin{Lem}\label{lem:ori-principal}
Let $\Gamma$ be a Jacobi diagram and $e$ be an edge of $\Gamma$. Then the orientation on the principal face induced from the orientation $\Omega(\Gamma)$ is as follows.
\[\begin{array}{c | r | r}
\mbox{type} & \mbox{space} & \mbox{induced orientation} \\\hline
\mbox{(1)} 
	& S^{n+1}\times S^{n-1}\times S^{n-1}\times C_{\delta_e\Gamma} 
	& -(-1)^n\omega_{n+1}\wedge\omega_{n-1}\wedge\omega_{n-1}\wedge\Omega(\delta_e\Gamma) \\
\mbox{(2)} 
	& S^{n-1}\times C_{\delta_e\Gamma}
	& (-1)^n\omega_{n-1}\wedge\Omega(\delta_e\Gamma)\\
\mbox{(3)} 
	& S^{n+1}\times S^{n-1}\times S^{n-1}\times C_{\delta_e\Gamma}
	& -(-1)^n\omega_{n+1}\wedge\omega_{n-1}\wedge\omega_{n-1}\wedge\Omega(\delta_e\Gamma)\\
\mbox{(4)}
	&S^{n-1}\times C_{\delta_e\Gamma}
	& \omega_{n-1}\wedge\Omega(\delta_e\Gamma)\\
\mbox{(5)} 
	& S^{n+1}\times S^{n-1}\times C_{\delta_e\Gamma}
	& (-1)^n\omega_{n+1}\wedge\omega_{n-1}\wedge\Omega(\delta_e\Gamma)\\
\mbox{(6)} 
	&  S^{n-1}\times C_{\delta_e\Gamma}
	& \omega_{n-1}\wedge\Omega(\delta_e\Gamma)\\
\mbox{(7)} 
	&  S^{n-1}\times C_{\delta_e\Gamma}
	& \omega_{n-1}\wedge\Omega(\delta_e\Gamma)\\
\end{array}\]

\end{Lem}
\begin{proof}
Let $\Gamma_1$ be as in the LHS of the case (1) and we consider the orientation induced on the principal face $\calS_{\Gamma_1,e_1}$ corresponding to the contraction of the edge $e_1=(w,y)$. The orientation on $C_{\Gamma_1}$ is 
\[ \begin{split}
	\Omega(\Gamma_1)&=(dX_y^1\omega_{n-1}^{(1)}d^nX_w)
		\wedge(dX_y^2\omega_{n-1}^{(2)}\Omega_{\bar{e}_v})
		\wedge(\Omega_{\bar{e}_u}dX_y^3\cdots dX_y^{n+2})
		\wedge\mbox{(the rest)}\\
	&=-(-1)^n dX_y^1\cdots dX_y^{n+2}\wedge\omega_{n-1}^{(1)}\wedge\omega_{n-1}^{(2)}
		\wedge d^nX_w\Omega_{\bar{e}_v}\Omega_{\bar{e}_u}
		\wedge\mbox{(the rest)}.
	\end{split} \]
Let $n_1$ is the outgoing unit normal vector field on $S^{n+1}$. Then the ingoing normal vector field on the principal face is given by $n_1$. The induced orientation on $\calS_{\Gamma_1,e_1}$ is then
\[ \begin{split}
	i_{n_1}\Omega(\Gamma_1)&=-(-1)^n\omega_{n+1}(\alpha)\wedge\omega_{n-1}^{(1)}\wedge\omega_{n-1}^{(2)}\wedge
		d^nX_w\Omega_{\bar{e}_v}\Omega_{\bar{e}_u}
		\wedge\mbox{(the rest)}\\
	&=-(-1)^n\omega_{n+1}(\alpha)\wedge\omega_{n-1}^{(1)}\wedge\omega_{n-1}^{(2)}\wedge\Omega(\delta_{e_1}\Gamma_1)
	\end{split} \]
where $\alpha=\frac{X_y-X_w}{|X_y-X_w|}$. For the cases (3), (5), the results follow from the following:
\[ \begin{split} i_{n_1}\Omega(\Gamma_3)&=-(-1)^n\omega_{n+1}\omega_{n-1}^{(1)}\omega_{n-1}^{(2)}
	d^nX_w\Omega_{\bar{e}_v}\Omega_{\bar{e}_u}\wedge\mbox{(the rest)}\\
	i_{n_1}\Omega(\Gamma_5)&=(-1)^{n}\omega_{n+1}\omega_{n-1}^{(4)}d^{n+2}X_x
		\omega_{n-1}^{(1)}\omega_{n-1}^{(2)}\omega_{n-1}^{(3)}
	\Omega_{\bar{e}_b}\Omega_{\bar{e}_a}\Omega_{\bar{e}_c}\Omega_{\bar{e}_d}\wedge\mbox{(the rest)}.
	\end{split} \]

Let $\Gamma_2$ be as in the LHS of the case (2) and we consider the orientation induced on the principal face $\calS_{\Gamma_2,e_2}$ corresponding to the contraction of the edge $e_2=(w,j)$. The orientation on $C_{\Gamma_2}$ is 
\[ \Omega(\Gamma_2)=d^nX_w\Omega_{\bar{e}_v}\Omega_{\bar{e}_u}d^nX_j\wedge\mbox{(the rest)}=d^nX_wd^nX_j\Omega_{\bar{e}_v}\Omega_{\bar{e}_u}\wedge\mbox{(the rest)}. \]
Let $n_2$ be the outgoing unit normal vector field on $S^{n-1}$. Then the induced orientation on $\calS_{\Gamma_2,e_2}$ is 
\[ i_{n_2}\Omega(\Gamma_2)=(-1)^n\omega_{n-1}(\beta)\wedge d^nX_w\Omega_{\bar{e}_u}\Omega_{\bar{e}_v}\wedge\mbox{(the rest)}
	=(-1)^n\omega_{n-1}(\beta)\wedge\Omega(\delta_{e_2}\Gamma_2) \]
where $\beta=\frac{X_j-X_w}{|X_j-X_w|}$. The case (4) is similar from the following:
\[ i_{n_2}\Omega(\Gamma_4)=\omega_{n-1}\Omega_{\bar{e}_w}\Omega_{\bar{e}_u}\Omega_{\bar{e}_v}\wedge\mbox{(the rest)}. \]

Let $\Gamma_6$ be as in the LHS of the case (6) and we consider the orientation induced on the principal face $\calS_{\Gamma_6,e_6}$ corresponding to the contraction of the edge $e_6=(i,j)$. The orientation on $C_{\Gamma_6}$ is 
\[ \Omega(\Gamma_6)=d^nX_jd^nX_i\wedge\mbox{(the rest)}. \]
The induced orientation on $\calS_{\Gamma_6,e_6}$ is 
\[ i_{n_2}\Omega(\Gamma_6)=\omega_{n-1}(\gamma)d^nX_i\wedge\mbox{(the rest)}
	=\omega_{n-1}(\gamma)\wedge\Omega(\delta_{e_6}\Gamma_6) \]
where $\beta=\frac{X_j-X_i}{|X_j-X_i|}$. The case (7) is similar.
\end{proof}

\begin{Lem}\label{lem:I(dG)}
Let $(\Gamma_i,e_i), (i=1,\ldots,7)$ be the pair of the Jacobi diagram and the edge to be contracted on it as in the LHS of the case ($i$) in (\ref{eq:cntr}). Then we have
\[ \begin{split}
	I(\Gamma_1,e_1)&=-I(\Gamma_2,e_2), \\
	I(\Gamma_3,e_3)&=-I(\Gamma_4,e_4), \\
	I(\Gamma_6,e_6)&=I(\Gamma_7,e_7).
	\end{split}
\]
\end{Lem}
\begin{proof}
By definition, 
\[ \begin{split}
	I(\Gamma_1,e_1)&=\int_{S^{n+1}\times S^{n-1}\times S^{n-1}\times C_{\delta_{e_1}\Gamma_1}}
		\theta_{e_1}\wedge\omega_{n-1}^{(1)}\wedge\omega_{n-1}^{(2)}\wedge\lambda_2(\delta_{e_1}\Gamma_1)\\
	&=\int_{S^{n+1}}\theta_{e_1}\int_{S^{n-1}}\omega_{n-1}^{(1)}\int_{S^{n-1}}\omega_{n-1}^{(2)}\int_{C_{\delta_{e_1}\Gamma_1}}\lambda_2(\delta_{e_1}\Gamma_1)
		=\int_{C_{\delta_{e_1}\Gamma_1}}\lambda_2(\delta_{e_1}\Gamma_1)
\end{split} \]
where $\theta_{e_1}$ denotes the pullback form of $\omega_{n+1}$ via the Gauss map $u$ with respect to $e_1$. Similarly we have
\[ I(\Gamma_2,e_2)
		=\int_{C_{\delta_{e_2}\Gamma_2}}\lambda_2(\delta_{e_2}\Gamma_2). \]
Since the induced orientations on $C_{\delta_{e_1}\Gamma_1}\cong C_{\delta_{e_2}\Gamma_2}$ are opposite by Lemma~\ref{lem:ori-principal}, it suffices to prove $\lambda_2(\delta_{e_1}\Gamma_1)=\lambda_2(\delta_{e_2}\Gamma_2)$. In terms of the labels of vertices in (\ref{eq:cntr}), we have
\[ \begin{split}
	\omega(\Gamma_1)&=\theta_{wy}\omega_{n-1}^{(1)}\wedge\theta_{vy}\omega_{n-1}^{(2)}\wedge\theta_{yu}\alpha\wedge\mbox{(the rest)}\\
		&=(-1)^{n-1}\theta_{wy}\omega_{n-1}^{(1)}\omega_{n-1}^{(2)}\wedge\theta_{vy}\theta_{yu}\wedge\mbox{(the rest)}
		\end{split} \]
for some $(n-1)$-form $\alpha$. On the principal face corresponding to the contraction of $e_1=(w,y)$, this form becomes $(-1)^{n-1}\theta_{e_1}\omega_{n-1}^{(1)}\omega_{n-1}^{(2)}\wedge\theta_{vw}\theta_{wu}\wedge\mbox{(the rest)}$ and thus $\lambda_2(\delta_{e_1}\Gamma_1)=(-1)^{n-1}\theta_{vw}\theta_{wu}\alpha\wedge\mbox{(the rest)}$. On the other hand, we have
\[ \omega(\Gamma_2)=(-1)^{n-1}\eta_{wj}\theta_{vw}\theta_{ju}\alpha\wedge\mbox{(the rest)}. \]
Thus $\lambda_2(\delta_{e_2}\Gamma_2)=(-1)^{n-1}\theta_{vw}\theta_{wu}\alpha\wedge\mbox{(the rest)}\equiv \lambda_2(\delta_{e_1}\Gamma_1)$.

The second and the third identity may be similarly proved from the facts: 
\[ \begin{split}
	\mbox{induced ori. on $C_{\delta_{e_3}\Gamma_3}\cong C_{\delta_{e_4}\Gamma_4}$ differ by $(-1)^{n-1}$},\ &\lambda_2(\delta_{e_3}\Gamma_3)=-(-1)^{n-1}\lambda_2(\delta_{e_4}\Gamma_4),\\
	\mbox{induced ori. on $C_{\delta_{e_6}\Gamma_7}\cong C_{\delta_{e_6}\Gamma_7}$ are the same},\ &\lambda_2(\delta_{e_6}\Gamma_6)=\lambda_2(\delta_{e_7}\Gamma_7).
	\end{split} \]
\end{proof}

\begin{proof}[Proof of Proposition~\ref{prop:principal}]
$\phantom{a}$\par

We consider the faces corresponding to the collapse of an edge of $\Gamma$. The proof is divided according to the cases appeared in the definition of $\delta_e$.

\begin{description}
\item[\underline{\bf (1), (2)}]{\it Appearance of the ST relation.} For the cases (1), (2), we will see that the contributions of the faces corresponding to the following pair of contractions of graphs cancel each other:
\[ \fig{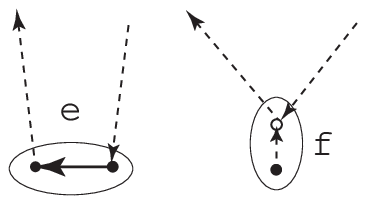} \]
which we denote by $(\Gamma^{\mathrm{a}},e), (\Gamma^{\mathrm{b}},f)$ respectively.

From (\ref{eq:type-space}), the face of $(\Gamma^{\rma},e)$ is identified with $S^{n-1}\times C_{\delta_e\Gamma^{\mathrm{a}}}$, where the factor $S^{n-1}$ is the relative configuration space of 2 points on $T_x\psi(\R^n)$. The face of $(\Gamma^\rmb,f)$ is identified with $S^{n+1}\times C_{\delta_f\Gamma^{\mathrm{b}}}$, where the factor $S^{n+1}$ is the relative configuration space of 2 points on $\R^{n+2}$.

Let $n_{\mathrm{a}}$ be the number of $\eta$-edges $e'$ on $\Gamma^{\rma}$ such that $\delta_e\Gamma^{\rma}=\delta_{e'}\Gamma^{\rma}$. $n_{\rmb}$ is similarly defined to be the number of $\eta$-edges $f'$ on $\Gamma^{\rmb}$ such that $\delta_f\Gamma^{\rmb}=\delta_{f'}\Gamma^{\rmb}$. Then we have
\[
\frac{|\Aut\Gamma^{\rma}|}{n_{\rma}}=\frac{|\Aut\Gamma^{\rmb}|}{n_{\rmb}}=|\Aut\delta_e\Gamma^{\rma}|=|\Aut\delta_f\Gamma^{\rmb}|.
\]
We assume that each automorphism in $\Aut\Gamma^{\rma}$ preserves the orientation of $C_{\Gamma^{\rma}}$ since if not $I(\Gamma^{\rma},e)=0$. The case of $\Gamma^{\rmb}$ is similar. Since the integration along the sphere factors contributes by 1, $z_k$ restricted to these faces is
\[ 	\frac{n_{\rma}I(\Gamma^{\mathrm{a}},e)w_k(\Gamma^{\mathrm{a}})}{|\Aut\Gamma^{\mathrm{a}}|}
	+\frac{n_{\rmb}I(\Gamma^{\mathrm{b}},f)w_k(\Gamma^{\mathrm{b}})}{|\Aut\Gamma^{\mathrm{b}}|}
	=\frac{I(\Gamma^{\mathrm{a}},e)}{|\Aut\delta_e\Gamma^{\mathrm{a}}|}(w_k(\Gamma^{\mathrm{a}})-w_k(\Gamma^{\mathrm{b}}))=0,
\]
by $I(\Gamma^{\mathrm{b}},f)=-I(\Gamma^{\mathrm{a}},e)$ by Lemma~\ref{lem:I(dG)} and the ST relation.

\item[\underline{\bf (3), (4)}]{\it Appearance of the STU relation.} For the cases (3), (4), we will see that the contributions of the faces corresponding to the following triple of graphs cancel each other:
\[ \fig{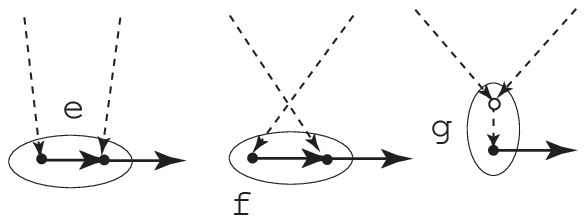} \]
which we denote by $(\Gamma^{\mathrm{c}},e), (\Gamma^{\mathrm{d}},f), (\Gamma^{\mathrm{e}},g)$ respectively.
From (\ref{eq:type-space}), the faces corresponding to $(\Gamma^{\mathrm{c}},e), (\Gamma^{\mathrm{d}},f), (\Gamma^{\mathrm{e}},g)$ can be identified with the following spaces:
\begin{center}
	\begin{tabular}{ll}
		$(\Gamma^{\mathrm{c}},e)$: & $S^{n-1}\times C_{\delta_e\Gamma^{\mathrm{c}}} ,$\\
		$(\Gamma^{\mathrm{d}},f)$: & $S^{n-1}\times C_{\delta_f\Gamma^{\mathrm{d}}},$\\
		$(\Gamma^{\mathrm{e}},g)$: & $S^{n+1}\times C_{\delta_g\Gamma^{\mathrm{e}}}$.
	\end{tabular}
\end{center}
Let $n_{\rmc}$ be the number of $\eta$-edges $e'$ on $\Gamma^{\rmc}$ such that $\delta_e\Gamma^{\rmc}=\delta_{e'}\Gamma^{\rmc}$. $n_{\rmd}, n_{\rme}$ are similarly defined. Then we have 
\[
	\frac{|\Aut\Gamma^{\rmc}|}{n_{\rmc}}=
	\frac{|\Aut\Gamma^{\rmd}|}{n_{\rmd}}=
	\frac{|\Aut\Gamma^{\rme}|}{n_{\rme}}=
	|\Aut\delta_e\Gamma^{\rmc}|=
	|\Aut\delta_f\Gamma^{\rmd}|=
	|\Aut\delta_g\Gamma^{\rme}|.
\]
We assume that each automorphism in $\Aut\Gamma^{\rmc}$ preserves the orientation of $C_{\Gamma^{\rmc}}$ since if not $I(\Gamma^{\rmc},e)=0$. The cases of $\Gamma^{\rmd}, \Gamma^{\rme}$ are similar.  Since the integration along the sphere factors contributes by 1, $z_k$ restricted to these faces is
\[ \begin{split}
	\frac{n_{\rmc}I(\Gamma^{\mathrm{c}},e)w_k(\Gamma^{\mathrm{c}})}{|\Aut\Gamma^{\mathrm{c}}|}
	&+\frac{n_{\rmd}I(\Gamma^{\mathrm{d}},f)w_k(\Gamma^{\mathrm{d}})}{|\Aut\Gamma^{\mathrm{d}}|}
	+\frac{n_{\rme}I(\Gamma^{\mathrm{e}},g)w_k(\Gamma^{\mathrm{e}})}{|\Aut\Gamma^{\mathrm{e}}|}\\
	=\frac{I(\Gamma^{\mathrm{c}},e)}{|\Aut\delta_e\Gamma^{\mathrm{c}}|}&(w_k(\Gamma^{\mathrm{c}})-w_k(\Gamma^{\mathrm{d}})-w_k(\Gamma^{\mathrm{e}}))=0,
	\end{split}
\]
by $I(\Gamma^{\mathrm{e}},g)=-I(\Gamma^{\mathrm{c}},e)$ by Lemma~\ref{lem:I(dG)}, and $I(\Gamma^{\mathrm{d}},f)=-I(\overline{\Gamma}^{\mathrm{e}},g)=I(\Gamma^{\mathrm{e}},g)$ because $\Omega(\delta_{g}\overline\Gamma^{\mathrm{e}})=(-1)^n\Omega(\delta_{g}\Gamma^{\mathrm{e}})$ and $\lambda_2(\delta_{g}\overline\Gamma^{\mathrm{e}})=(-1)^{n-1}\lambda_2(\delta_{g}\Gamma^{\mathrm{e}})$, and the STU relation.

\item[\underline{\bf (6),(7)}]{\it Appearance of the C relation.}
Then $z_k$ can be written without using $|\Aut\Gamma|$ as follows:
\[ z_k=\frac{1}{(2k)!}\sum_{\Gamma\mathrm{:labeled}}I(\Gamma)w_k(\Gamma), \]
since there are $(2k)!/|\Aut\Gamma|$ ways of labeling the $2k$ vertices with distinguished labels in $\{1,2,\cdots,2k\}$. There are just 4 labeled graphs $\Gamma^{\mathrm{L}}_1, \Gamma^{\mathrm{L}}_2, \Gamma^{\mathrm{R}}_1, \Gamma^{\mathrm{R}}_2$ as in Figure~\ref{fig:c1d}(a),(a'),(b),(b') respectively which are mapped by $\delta_e$ to the graph as in Figure~\ref{fig:c1d}(c). 
\begin{figure}
\cfig{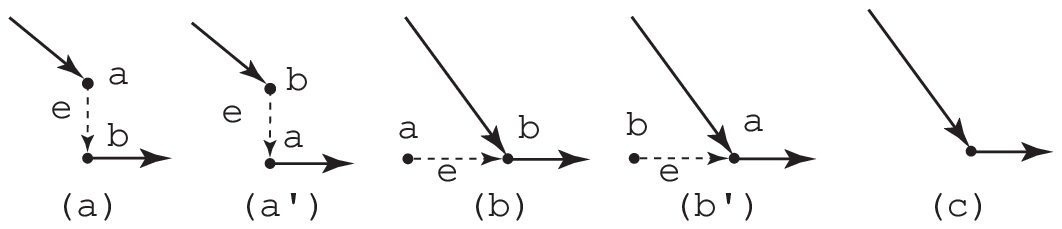}
\caption{}\label{fig:c1d}
\end{figure}
From (\ref{eq:type-space}), the corresponding principal faces can be identified with the following spaces via the tangent framing:
\begin{center}
		  $S^{n-1}\times C_{\delta_e\Gamma^{\mathrm{L}}_1}$,
		  $S^{n-1}\times C_{\delta_e\Gamma^{\mathrm{L}}_2}$,
		  $S^{n-1}\times C_{\delta_e\Gamma^{\mathrm{R}}_1}$,
		  $S^{n-1}\times C_{\delta_e\Gamma^{\mathrm{R}}_2}$.
\end{center}
Since the integration along $S^{n-1}$ fiber contributes by some common $2$-form $\mu$ on $I_n(\R^{n+2})$, the integrals restricted to these faces are\[ \begin{split}
	\frac{\mu}{(2k)!}&(I(\Gamma^{\mathrm{L}}_1,e)w_k(\Gamma^{\mathrm{L}}_1)
	+I(\Gamma^{\mathrm{L}}_2,e)w_k(\Gamma^{\mathrm{L}}_2)
	+I(\Gamma^{\mathrm{R}}_1,e)w_k(\Gamma^{\mathrm{R}}_1)
	+I(\Gamma^{\mathrm{R}}_2,e)w_k(\Gamma^{\mathrm{R}}_2))\\
	=\frac{1}{(2k)!}&
	(
	w_k(\Gamma^{\mathrm{L}}_1)+w_k(\Gamma^{\mathrm{L}}_2)
	+w_k(\Gamma^{\mathrm{R}}_1)+w_k(\Gamma^{\mathrm{R}}_2))
	\int_{C_{\delta_e\Gamma^{\mathrm{L}}_1}}\varphi^*\mu\wedge\lambda_2(\delta_e\Gamma^{\mathrm{L}}_1)=0,
	\end{split}
\]
by $I(\Gamma^{\mathrm{L}}_i,e)=I(\Gamma^{\mathrm{R}}_i,e)$ by Lemma~\ref{lem:I(dG)} and the C relation.

\item[\underline{(5)}]{\it Appearance of the IHX relation.} There are just 6 labeled graphs $\Gamma^{\mathrm{I}}_1, \Gamma^{\mathrm{I}}_2, \Gamma^{\mathrm{H}}_1, \Gamma^{\mathrm{H}}_2, \Gamma^{\mathrm{X}}_1, \Gamma^{\mathrm{X}}_2$ as in Figure~\ref{fig:ihxd}(a), (a'), (b), (b'), (c), (c') respectively which are mapped by $\delta_e$ to the graph as in Figure~\ref{fig:ihxd}(d) (the $-$ sign of (b) and (b') is because the cyclic permutation of 4 edges is an odd permutation). From (\ref{eq:type-space}), the corresponding principal faces can be identified with the following spaces:
\begin{center}
	\begin{tabular}{ll}
		  $S^{n+1}\times C_{\delta_e\Gamma^{\mathrm{I}}_1}$,
		 & $S^{n+1}\times C_{\delta_e\Gamma^{\mathrm{I}}_2}$,\\
		  $S^{n+1}\times C_{\delta_e\Gamma^{\mathrm{H}}_1}$,
		 & $S^{n+1}\times C_{\delta_e\Gamma^{\mathrm{H}}_2}$,\\
		  $S^{n+1}\times C_{\delta_e\Gamma^{\mathrm{X}}_1}$,
		 & $S^{n+1}\times C_{\delta_e\Gamma^{\mathrm{X}}_2}$.
	\end{tabular}
\end{center}
\begin{figure}
\cfig{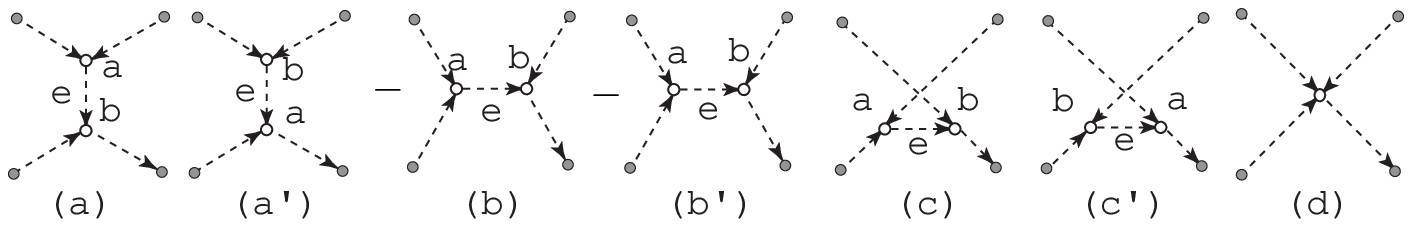}
\caption{}\label{fig:ihxd}
\end{figure}
Since the integration along the sphere factors contributes by 1, the integral restricted to these faces is
\[ 	\begin{split}
	\frac{1}{(2k)!}&(I(\Gamma^{\mathrm{I}}_1,e)w_k(\Gamma^{\mathrm{I}}_1)
	+I(\Gamma^{\mathrm{I}}_2,e)w_k(\Gamma^{\mathrm{I}}_2)
	+I(\Gamma^{\mathrm{H}}_1,e)w_k(\Gamma^{\mathrm{H}}_1)\\
	&+I(\Gamma^{\mathrm{H}}_2,e)w_k(\Gamma^{\mathrm{H}}_2)
	+I(\Gamma^{\mathrm{X}}_1,e)w_k(\Gamma^{\mathrm{X}}_1)
	+I(\Gamma^{\mathrm{X}}_2,e)w_k(\Gamma^{\mathrm{X}}_2))\\
	&=\frac{I(\Gamma^{\mathrm{I}}_1,e)}{(2k)!}
	(w_k(\Gamma^{\mathrm{I}}_1)+w_k(\Gamma^{\mathrm{I}}_2)-w_k(\Gamma^{\mathrm{H}}_1)-w_k(\Gamma^{\mathrm{H}}_2)+w_k(\Gamma^{\mathrm{X}}_1)+w_k(\Gamma^{\mathrm{X}}_2))=0,
	\end{split}
\]
by the IHX relation.
\end{description}
\end{proof}
\subsection{Hidden faces}
Now we prove the vanishing of $I(\Gamma)w_k(\Gamma)$ restricted to the hidden faces. Let $A\subset V(\Gamma), |A|\geq 3$ be a set of points to be collapsed having external edges (edges connecting points in $A$ and points not in $A$) and let $\Gamma_A$ be the subgraph of $\Gamma$ consisting of edges connecting two vertices in $A$. The goal of this subsection is to prove the following proposition.
\begin{Prop}\label{prop:hiddenall}
The term $I(\Gamma)w_k(\Gamma)$ restricted to any face corresponding to the diagonal where the points in $A\subset V(\Gamma)$ with $A\neq V(\Gamma)$, $|A|\geq 3$ coincide, vanishes.
\end{Prop}
The proof of Proposition~\ref{prop:hiddenall} is almostly done by using the Vanishing Lemmata (Lemma~\ref{lem:involution}, \ref{lem:hiddenonknot}, \ref{lem:hiddenchords} below) proved in \cite{R} and \cite{CCL}.
\begin{Lem}\label{lem:hiddenexternaledge}
If $|A|\geq 3$ and $A$ has external edges, then the integral restricted to the corresponding face vanishes.
\end{Lem}
\begin{proof}
By Lemma~\ref{lem:notconnected} below, we need only to consider the cases in which $\Gamma_A$ is connected.

Let $a$ be a vertex in $\Gamma_A$ which is an end of an external edge. The possible cases in which $\Gamma_A$ are connected are as in the following picture:
\cfig{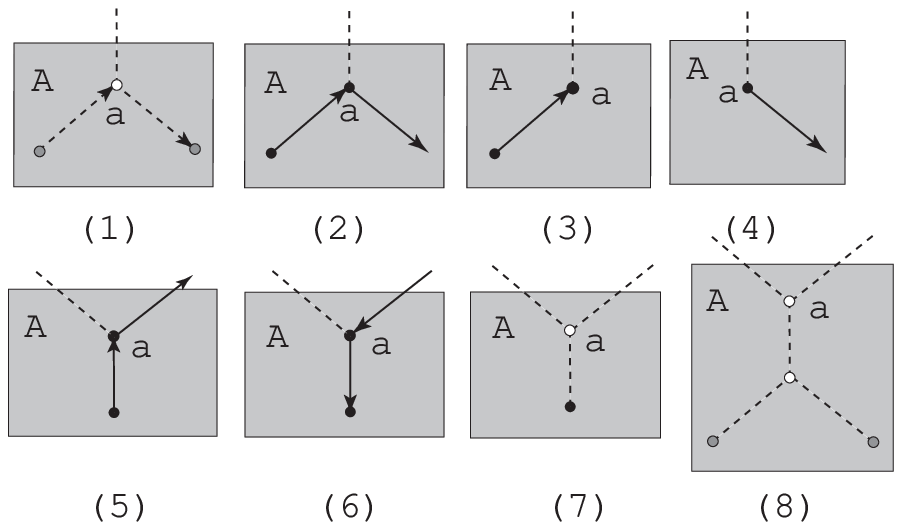}
The reasons for the vanishing of the integral $I(\Gamma)$ along $\calS_{\Gamma,A}$ are as follows:
\begin{description}
\item[(1),(2)] By Lemma~\ref{lem:involution} below.
\item[(3),(4),(5),(6),(7),(8)] By Lemma~\ref{lem:hiddenonknot} below.
\end{description}
\end{proof} 
\begin{proof}[Proof of Proposition~\ref{prop:hiddenall}]
By Lemma~\ref{lem:hiddenexternaledge} and Proposition~\ref{prop:A_n}, for a non vanishing face the $\theta$ part of $\Gamma_A$ must form a disjoint union of chords, trees and wheels entirely included in $\Gamma_A$. So we assume this condition in the rest of the proof.

If there does not exists an $\eta$-edge going from a point in $A$ to a point not in $A$, then $w_k(\Gamma)=0$ by Y relation or L relation.

Suppose that there exists an $\eta$-edge going from a point in $A$ to a point not in $A$. Then $\Gamma_A$ consists of a sequence of subgraphs as in (\ref{eq:hairy-strut})(b) joined by $\eta$-edges. If $\Gamma_A$ has a tree component with at least one internal vertex, then consider the automorphism on the fiber of $\widehat{B}_A$ defined as in Figure~\ref{fig:inv-tree}, which reverses (resp. preserves) the orientation of the fiber and preserves (resp. reverses) the sign of the form if $n$ is odd (resp. even). Hence the integral vanishes. If $\Gamma$ has only chords in its $\theta$ part, then by the assumption, there has to be at least two chords connected by an $\eta$-edge. But this case does not contribute by Lemma~\ref{lem:hiddenchords}.
\begin{figure}
\fig{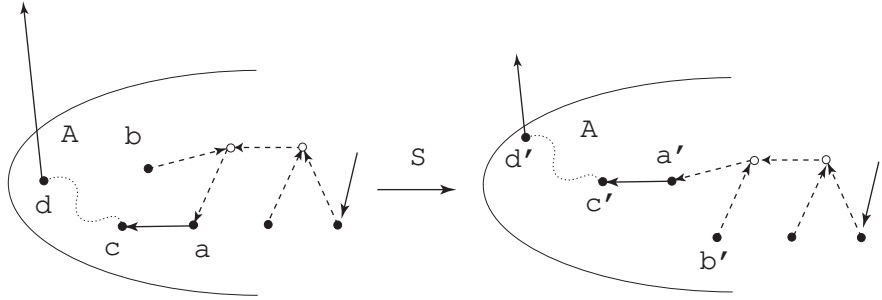}
\caption{$(x_d,\ldots,x_c,x_a,x_b,\ldots)\atop\mapsto(x_d+(x_b-x_a),\ldots,x_c+(x_b-x_a),x_b,x_a,\ldots)$}\label{fig:inv-tree}
\end{figure}
\end{proof}

\begin{Lem}\label{lem:notconnected}
If $\Gamma_A$ is not connected, then the integral restricted to $\calS_{\Gamma,A}$ vanishes.
\end{Lem}
\begin{proof}
By the observations in \S\ref{ss:face}, the face $\calS_{\Gamma,A}$ is a bundle $\calS_{\Gamma_A}\to C_{\Gamma/A}$ with fiber $F_A$ the configuration space associated to the subgraph $\Gamma_A$ modulo translations and dilations. Since $\Gamma_A$ is not connected, we can consider the action $T$ on $F_A$ which translates one connected component in $\Gamma_A$ fixing all other components. It is easy to see from $|A|\geq 3$ that this action is non trivial on $F_A$ and the quotient map $F_A\to F_A/T$ gives rise to a bundle $\pi^T:\calS_{\Gamma_A}\to\calS_{\Gamma_A}/T$ with $n$ or $(n+2)$ dimensional fiber, depending on whether the translated component has external vertices or not. Since the $T$-action does not affect to the integrand form $\omega(\Gamma)$ extended to $\calS_{\Gamma,A}$, namely the integrand form can be written as $(\pi^T)^*\omega(\Gamma)^T$ for some form $\omega(\Gamma)^T$ on $\calS_{\Gamma,A}/T$. Thus 
\[ \int_{\calS_{\Gamma,A}}\omega(\Gamma)=\int_{\calS_{\Gamma}}(\pi^T)^*\omega(\Gamma)^T
	=\int_{\calS_{\Gamma}/T}\omega(\Gamma)^T \]
where the integration is the fiber integration. The last integral is over the space with codimension at least $(n+1)$ and so it vanishes.
\end{proof}
The following lemma is an analogue of a lemma in \cite{Kon}, which is proved in \cite{R}.
\begin{Lem}\label{lem:involution}
If $\Gamma_A$ has a subgraph of one of the following forms:
\cfig{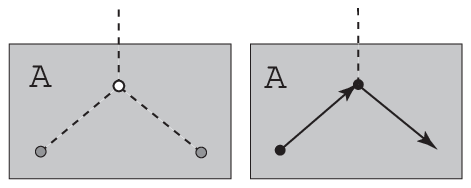}
where the gray vertices may both be internal and external, then the integral restricted to $\calS_{\Gamma,A}$ vanishes.
\end{Lem}

The following two lemmas are proved in \cite{R}.
\begin{Lem}\label{lem:hiddenonknot}
If $|A|\geq 3$ and $\Gamma_A$ has a subgraph of one of the following forms:
\cfig{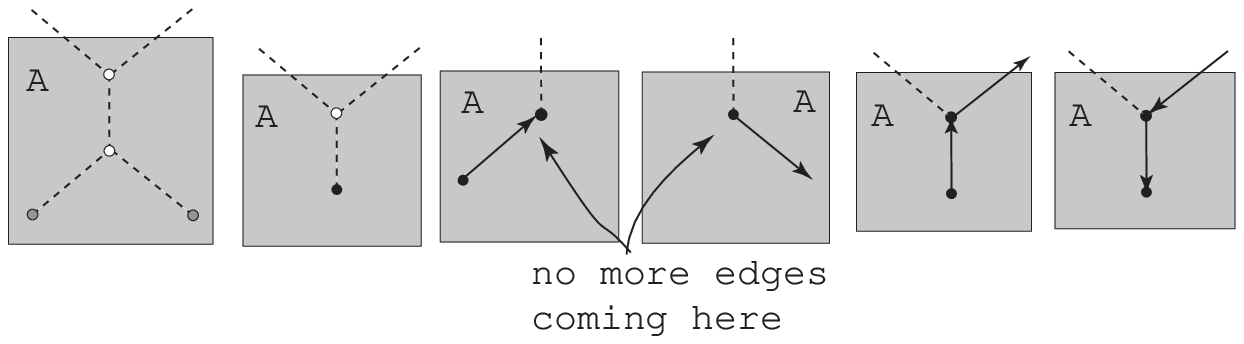}
then the integral restricted to $\calS_{\Gamma,A}$ vanishes.
\end{Lem}

\begin{Lem}\label{lem:hiddenchords}
If $\Gamma_A$ has a subgraph of the form:
\cfig{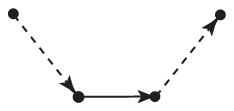}
then the integral restricted to $\calS_{\Gamma,A}$ vanishes.
\end{Lem}

\subsection{Infinite and anomalous faces}

For the infinite faces, the following proposition is proved in \cite{R}.
\begin{Prop}\label{prop:infvanish}
The integral $I(\Gamma)$ restricted to the infinite faces for $\infty\in S^{n+2}$ vanishes.
\end{Prop}

\subsubsection{Anomaly correction in the case $n=2$}\label{ss:anomaly-n=2}
In the case $n=2$ and $k=3$, the following proposition is proved in \cite{CR,R}.
\begin{Prop}\label{prop:anomaly-n=2}
If $n=2$, then there exists a 2-form $\hat{\rho}(\Gamma)$ on $I_2(\R^4)$ such that
\[ d\int_{C_1}\varphi^*\hat{\rho}(\Gamma)=-I(\Gamma,V(\Gamma)). \]
\end{Prop}

\subsubsection{Anomaly vanishing in the case $n$ odd}\label{ss:anomaly-n-odd}
We show the following proposition as already announced in \cite{R}. Recall that $I(\Gamma,V(\Gamma))=\int_{C_1}\varphi^*\pi_*^{\partial}\lambda_1(\Gamma)$.
\begin{Prop}\label{prop:anomaly-n-odd}
$\pi_*^{\partial}\lambda_1(\Gamma)$ vanishes for any $\Gamma$ in the case $n$ is odd.
\end{Prop}
\begin{proof}
Note first that each $\Gamma$ has even number of vertices and even number of edges. Let $x_1,\ldots,x_{j}$ be external vertices of $\Gamma$ mapped into an underlying $n$ dimensional plane with a tangent frame $f\in I_n(\R^{n+2})$ and $y_1,\ldots,y_{k}$ be internal vertices of $\Gamma$ in the ambient space $\R^{n+2}$ so that $(f;x_1,\ldots,x_{j},y_1,\ldots,y_{k})\in\widehat B_{V(\Gamma)}$. Then consider the involution
\[ \begin{split}
	S:&(f;x_1,x_2,\ldots,x_{j},y_1,\ldots,y_{k})\\
	&\mapsto(f;x_1,2x_1-x_2,\ldots,2x_1-x_{j},2f(x_1)-y_1,\ldots,2f(x_1)-y_{k})
	\end{split} \]
on the fiber. Since $j+k-1$ is odd if $n$ is odd, $S$ reverses the orientation of the fiber of $\widehat B_{V(\Gamma)}$ and preserves the form $\lambda_1(\Gamma)$. Therefore the integration along the fiber vanishes as in the proof of Lemma~\ref{lem:involution}.
\end{proof}

\end{appendix}

\end{document}